\documentclass[12pt,twoside,reqno]{amsart}
\usepackage[colorlinks=true,citecolor=blue]{hyperref}
\usepackage{mathptmx, amsmath, amssymb, amsfonts, amsthm, mathptmx, enumerate, color}
\usepackage{graphicx}
\usepackage{epstopdf}
\usepackage{framed}
\usepackage{booktabs} 
\usepackage{epstopdf}
\usepackage{algorithm}
\usepackage{algorithmic}
\usepackage{url}

\newtheorem{theorem}{Theorem}[section]
\newtheorem{corollary}{Corollary}[section]
\newtheorem{lemma}{Lemma}[section]

\theoremstyle{definition}

\newtheorem{remark}{Remark}[section]

\newtheorem{assumption}{Assumption}[section]
\numberwithin{equation}{section}

\newcommand{\R}{\mathbb{R}}
\newcommand{\Rext}{\mathbb{R}\cup\{+\infty\}}
\newcommand{\abs}[1]{\left\vert#1\right\vert}
\newcommand{\set}[1]{\left\{#1\right\}}
\newcommand{\norm}[1]{\left\Vert#1\right\Vert}
\newcommand{\Eproof}{\hfill $\square$}

\newcommand{\argmin}{\mathrm{arg}\min}
\newcommand{\dom}[1]{\mathrm{dom}(#1)}
\newcommand{\Ac}{\mathcal{A}}
\newcommand{\Pc}{\mathcal{P}}
\newcommand{\Rc}{\mathcal{R}}
\newcommand{\Xc}{\mathcal{X}}
\newcommand{\Sc}{\mathcal{S}}
\newcommand{\Lc}{\mathcal{L}}
\newcommand{\Nc}{\mathcal{N}}
\newcommand{\iprods}[1]{\langle #1\rangle}
\newcommand{\Id}{\mathbb{I}}
\newcommand{\rank}[1]{\mathrm{rank}\left(#1\right)}
\newcommand{\trace}[1]{\mathrm{trace}\left(#1\right)}
\renewcommand{\vec}[1]{\mathrm{vec}\left(#1\right)}
\newcommand{\mat}[1]{\mathrm{mat}\left(#1\right)}

\begin{document}
\setcounter{page}{1}

\vspace*{1.0cm}
\title[Gauss-Newton-Based Algorithms for A Class of Matrix  Optimization Problems]{
Extended Gauss-Newton and ADMM-Gauss-Newton algorithms for low-rank matrix optimization
}
\author[Q. Tran-Dinh]{Quoc Tran-Dinh$^{*}$}
\maketitle
\vspace*{-0.6cm}

\begin{center}
{\footnotesize {\it
Department of Statistics and Operations Research\\
The University of North Carolina at Chapel Hill, Chapel Hill, NC 27599.

}}\end{center}

\vskip 4mm {\small \noindent {\bf Abstract.}
In this paper, we develop a variant of the well-known Gauss-Newton (GN) method to solve a class of nonconvex optimization problems involving low-rank matrix variables.
As opposed to standard GN method, our algorithm allows one to handle general smooth convex objective function. 
We show, under mild conditions, that the proposed algorithm globally and locally converges to a stationary point of the original problem. 
We also show empirically that the GN algorithm achieves higher accurate solutions than the alternating minimization algorithm (AMA).
Then, we specify our GN scheme to handle the symmetric case and prove its convergence, where AMA is not applicable.
Next, we incorporate our GN scheme into the alternating direction method of multipliers (ADMM) to develop an  ADMM-GN algorithm.
We prove that, under mild conditions and a proper choice of the penalty parameter, our ADMM-GN  globally converges to a  stationary point of the original problem.
Finally, we provide several numerical experiments to illustrate the  proposed algorithms.
Our results show that the new algorithms have encouraging performance compared to existing methods.

}

\renewcommand{\thefootnote}{}
\footnotetext{ $^*$Corresponding author.
\par
E-mail addresses: \texttt{quoctd@email.unc.edu} (Quoc Tran-Dinh).
\par
The first version was on Arxiv on June 10, 2016.
}

\section{Introduction}\label{sec:intro}
\paragraph{\textbf{Problem statement:}}
In this paper, we consider the following class of low-rank matrix nonconvex optimization problems:
\begin{equation}\label{eq:LMA_prob}
\Phi^{\star} := \min_{U, V}\Big\{ \Phi(U,V) := \phi\left( \Ac(UV^{\top}) - B\right) + \Rc(U,V) \ : \  U\in\R^{m\times r}, \ V\in\R^{n\times r} \Big\},
\end{equation}
where $\Ac(Z) := [\trace{A_1^{\top}Z},  \trace{A_2^{\top}Z},  \cdots, \trace{A_l^{\top}Z}]$ for $l$ matrices $A_1,\cdots, A_l$ in $\R^{m\times n}$ is a linear operator;  $\phi : \R^l\to\Rext$ is a proper, closed and convex function; and  $B\in\R^l$ is an observed vector.
The function $\Rc$ is often referred to as a regularizer, which can be chosen as $\Rc(U,V) := \frac{1}{4}\Vert U^{\top}U - V^{\top}V\Vert_F^2$ as suggested in \cite{tu2015low}.  
Clearly, \eqref{eq:LMA_prob} is  nonconvex  due to the bilinear term $UV^{\top}$. 
Hence, it is NP-hard \cite{natarajan1995sparse}, and numerical methods for solving \eqref{eq:LMA_prob} aim at obtaining a local optimum or a stationary point of \eqref{eq:LMA_prob}.
In this paper, we are interested in the low-rank case, where $r \ll \min\set{m, n}$. 

Problem \eqref{eq:LMA_prob} covers various practical models in low-rank embedded problems, function learning, matrix completion in recommender systems, inpainting and compression in image processing, robust principal component analysis in statistics, and semidefinite programming relaxations in combinatorial optimization, see, e.g., \cite{Candes2012b,Candes2011a,Esser2010a,goldfarb2011convergence,kyrillidis2014matrix,liu2015efficient,wen2012solving}. 
Among these applications, the following  problems have been recently attracted a great attention.
The most common case is when  $\phi(\cdot) := (1/2)\Vert\cdot\Vert_2^2$, where \eqref{eq:LMA_prob} becomes a least-squares low-rank  approximation problem in compressive sensing (see, e.g., \cite{kyrillidis2014matrix}):
\begin{equation}\label{eq:compressive_sensing}
\min_{U, V}\big\{ (1/2)\Vert\Ac(UV^{\top}) - B \Vert_2^2 \ : \ U\in\R^{m\times r}, \ V\in\R^{n\times r} \big\}.
\end{equation} 
Here, the linear operator $\Ac$ is often assumed to satisfy a restricted isometric property (RIP) \cite{Candes2006} that allows us to recover an exact solution from a few number of observations in $B$.
In particular, if  $\Ac = \mathcal{P}_{\Omega}$, the projection on a given index subset $\Omega$, then \eqref{eq:compressive_sensing} covers the matrix completion model:
\begin{equation}\label{eq:mc}
\min_{U, V}\big\{ (1/2)\Vert \mathcal{P}_{\Omega}(UV^{\top}) - B_{\Omega}\Vert_F^2 \ : \ U\in\R^{m\times r}, \ V\in\R^{n\times r} \big\},
\end{equation} 
where $B_{\Omega}$ is the observed entries in $\Omega$.
If  $\Ac$ is an identity operator and $B\in\R^{m\times n}$ is a given, then \eqref{eq:compressive_sensing} becomes  a low-rank matrix factorization problem
\begin{equation}\label{eq:LMA_prob2}
\Phi^{\star} := \min_{U, V}\big\{ \Phi(U,V)  = (1/2)\Vert UV^{\top} - B\Vert^2_F ~:~ U\in\R^{m\times r}, ~V\in\R^{n\times r} \big\}.
\end{equation}
Especially, if  $U = V$ and  $B$ is symmetric positive definite, then \eqref{eq:LMA_prob2} reduces to 
\begin{equation}\label{eq:LMA_prob3}
\Phi^{\star} := \min_{U}\big\{ \Phi(U) := (1/2)\Vert UU^{\top} - B\Vert_F^2 \ : \ U\in\R^{n\times r} \big\},
\end{equation}
which is studied in \cite{liu2015efficient}.
Alternatively, if we choose $\Phi(U) := (1/2)\Vert \Ac(UU^{\top}) - B\Vert_F^2$ in \eqref{eq:compressive_sensing}, then \eqref{eq:LMA_prob} reduces to the case investigated in \cite{Bhojanapallli2015}.
While both special cases, \eqref{eq:LMA_prob2} and \eqref{eq:LMA_prob3}, possess a closed form solution via a truncated SVD and an eigenvalue decomposition, respectively, GN methods can also be applied to solve these problems. 
In \cite{liu2015efficient}, the authors  demonstrated the advantages of a GN method for solving \eqref{eq:LMA_prob3}  with significantly encouraging performance.

\vspace{0.5ex}
\paragraph{\textbf{Related work:}}
The low-rank structure is key to recast many existing problems into new frameworks or to design new models by means of regularizers to promote solution structures in various applications such as matrix completion (MC) \cite{Candes2012b}, robust principal component analysis (RPCA) \cite{Candes2011}, and their variants. 
Hitherto, extensions to group structured sparsity, low-rankness, tree models, and tensor representation have attracted a great attention in recent years, see, e.g., \cite{fazel2002matrix,grasedyck2013literature,huang2011learning,kyrillidis2015structured,Recht2010,Signoretto2014,Yu2014}.
A majority of research for low-rank models focuses on estimating sample complexity results for specific instances of \eqref{eq:LMA_prob}, while numerous recent papers revolve around the RPCA settings, MC, and their extensions \cite{Candes2011,Candes2012b,johnson1990matrix}.

Along with modeling, solution methods have also been extensively developed for solving concrete instances of \eqref{eq:LMA_prob} in low-rank matrix completion and recovery settings.  
Among various approaches, convex optimization is perhaps one of the most powerful tools to solve several instances of \eqref{eq:LMA_prob}, including MC, RPCA, their variants, and extensions. 
Unfortunately, convex models only provide an approximation to the low-rank model \eqref{eq:LMA_prob} by convex relaxations using, e.g., nuclear or max norms, which may not adequately approximate the desired rank.
Alternatively, nonconvex  as well as  discrete optimization methods have also been considered for solving \eqref{eq:LMA_prob}, see, e.g., \cite{burer2003nonlinear,kyrillidis2014matrix,Lin2009,Shen2012,wen2012solving,yu2014parallel}.
While these approaches work directly on the original problem \eqref{eq:LMA_prob}, they can only find a local optimum or a critical point, and strongly depend on the priori knowledge of problems, the initial points of algorithm, and predicted ranks. 
However, recent empirical evidence has been provided to support these approaches, and surprisingly, in many cases, they outperform the convex optimization approach in terms of ``accuracy'' to the original model,  and the overall computational time \cite{kyrillidis2014matrix,Lin2009,wen2012solving}.
Other approaches such as stochastic gradient descent, Riemann manifold descent,  greedy methods, parallel and distributed algorithms have also recently been studied for solving \eqref{eq:LMA_prob},  see, e.g., \cite{bottou2010large,johnson1990matrix,keshavan2009gradient,vandereycken2013low,yu2014parallel}.

\vspace{0.5ex}
\paragraph{\textbf{Motivation:}}
Gauss-Newton (GN) methods work extremely well for nonlinear least-squares problems \cite{Bjorck1996}. 
When $\phi$ is quadratic and the residual term $\Ac(UV^{\top})-B$ in \eqref{eq:LMA_prob} is small or zero at solutions, they can achieve local superlinear and even quadratic convergence rate. 
With a ``good'' initial point (i.e., close to the set of stationary points), GN methods often reach a stationary point within a few iterations \cite{Deuflhard2006}.
Such a ``good'' initial point can be obtained using priori knowledge of the problem and the underlying algorithm (e.g., steady states of dynamical systems, or previous iterations of the algorithm) as a warm-start strategy.

As in classical GN methods, we develop an iterative scheme  for solving \eqref{eq:LMA_prob} by using a linearization of $\Ac(UV^{\top}) - B$ and a quadratic surrogate of $\phi$. 
At each iteration, it requires to solve a simple convex problem to form a GN direction and then incorporates with a globalization strategy to update the next iteration. 
In our setting, computing GN  direction reduces to solving a linear least-squares problem.  
Comparing to the alternating minimization method (AMA) \cite{wen2012solving} that alternatively solves for each $U$ and $V$, GN simultaneously solves for $U$ and $V$ using the linearization of $UV^{\top}$.
We have observed that ({\emph{cf}} Subsection \ref{subsec:experiment0}) GN uses a linearization of $UV^{\top}$ providing a good local approximate model to $UV^{\top}$ compared to the alternating form $U\bar{V}^{\top}$ (or $\bar{U}V^{\top}$), when $U - \bar{U}$ (or $V - \bar{V})$ is relatively large. This makes AMA saturated and does not significantly improve the objective values.
In addition, without regularization, AMA may fail to converge as indicated by a counterexample in \cite{gonzalez2005accelerating}.
Moreover, AMA is not applicable to solving the symmetric case of \eqref{eq:LMA_prob} as shown in Section~\ref{sec:symmetric_case}, but the GN method is. 

While GN methods often use  in nonlinear least squares \cite{Nocedal2006}, they have not widely been exploited for  matrix optimization.
Our aim in this paper is to extend the GN method for solving a class of problems \eqref{eq:LMA_prob} with a general smooth convex objective function $\phi$ and low-rank matrix variables.
This paper is also inspired by a recent work \cite{liu2015efficient}, where the authors proposed a simple symmetric GN scheme to solve \eqref{eq:LMA_prob3}, and demonstrated its very encouraging performance.

\vspace{0.5ex}
\paragraph{\textbf{Contribution:}}
Our contribution in this paper can be summarized as follows:
\begin{itemize}
\item[$\mathrm{(a)}$] We extend the GN method to solve the low-rank matrix optimization problem \eqref{eq:LMA_prob} with smooth convex objective function $\phi$. 
We prove the existence of a GN direction and provide a closed form formulation to compute it.
We empirically show that our GN method can achieve higher accurate solutions than the well-known AMA scheme within the same number of iterations in certain cases. 

\item[$\mathrm{(b)}$] We show that there exists an explicit step-size to guarantee a descent property of the GN direction, 
which  allows us to perform a  backtracking linesearch procedure. 
We specify our framework to the symmetric case. 
Under mild conditions, we prove a global convergence of the proposed methods.

\item[$\mathrm{(c)}$] We prove a local linear and quadratic convergence rate of the full-step GN variant under standard assumptions imposed on \eqref{eq:LMA_prob} at its solution set.

\item[$\mathrm{(d)}$] We also combine an alternating direction method of multipliers (ADMM) and the GN method to obtain a new algorithm for handling \eqref{eq:LMA_prob}. 
Under standard assumptions  on \eqref{eq:LMA_prob}, we prove global convergence of the proposed algorithm.
\end{itemize}
Unlike AMA whose only achieves a sublinear convergence even with good initial points, GN methods may require additional computation for GN directions, but they can achieve a fast local linear or quadratic convergence rate, which is key for online and real-time implementations by using warm-start. 
Alternatively, gradient descent-based methods can achieve local linear convergence but often require much strong assumptions imposed on \eqref{eq:LMA_prob}. 
In contrast, GN methods work with the ``small residual'' setting under mild assumptions, and can easily achieve high accuracy solutions within a small number of iterations.

\vspace{0.5ex}
\paragraph{\textbf{Paper outline:}}
The rest of this paper is organized as follows. 
We first review basic concepts related to problem \eqref{eq:LMA_prob} in Section \ref{sec:preliminary}.
Section \ref{sec:GN_method} presents a linesearch GN method for solving \eqref{eq:LMA_prob} and its convergence guarantees.
Section \ref{sec:admm_gn_method} develops an ADMM-Gauss-Newton algorithm to solve \eqref{eq:LMA_prob} and investigates its global convergence.
Section \ref{sec:symmetric_case} specifies the GN algorithm to the symmetric case and proves its convergence.
Section \ref{sec:impl_aspects} discusses the implementation aspects of our algorithms and their extension to the nonsmooth objective function case.
Numerical experiments are conducted in Section \ref{sec:num_experiments} with several examples in different fields.
For the sake of presentation, we move all the  proofs in the main text  to the appendix.

\section{Basic notation and optimality condition}\label{sec:preliminary}
We briefly describe basic notation,  the optimality condition of \eqref{eq:LMA_prob}, and our assumptions.

\subsection{Basic notation and concepts}
For a matrix $X$,  $\sigma_{\min}(X)$ and $\sigma_{\max}(X)$ denote its positive smallest and largest singular values, respectively.
If $X$ is symmetric, then $\lambda_{\min}(X)$ and $\lambda_{\max}(X)$ denote its smallest and largest eigenvalues, respectively.
We use $X = P\Sigma Q^{\top}$ for SVD and $X = U\Lambda U^{-1}$ for eigenvalue  decomposition.
$X^{\dagger}$ denotes the Moore-Penrose pseudo-inverse of $X$. When $X$ is full-column rank, $X^{\dagger} = (X^{\top}X)^{-1}X^{\top}$. 
We define $P_{X} := XX^{\dagger}$ the projection onto the range space of $X$, and $P_{X}^{\perp} := \mathbb{I} - P_{X}$ the orthogonal projection of $P_X$, i.e., $P_XP^{\perp}_X = P_X^{\perp}P_X = 0$, where $\Id$ is the identity matrix. Clearly, $P_X^{\perp}X = 0$.
We define $\vec{X} = (X_{11}, \cdots X_{m1}, \cdots, X_{1n}, \cdots, X_{mn})^{\top}$ the vectorization of $X$, and $\mathrm{mat}$ the inverse mapping of $\mathrm{vec}$, i.e., $\mat{\vec{X}} = X$.
$X\otimes Y$ denotes the Kronecker product of $X$ and $Y$. We have $\vec{AXB} = (B^{\top}\otimes A)\vec{X}$ and $(A \otimes B)(C \otimes D) = AC\otimes BD$.
$\Ac^{*}$ denotes the adjoint of a linear operator $\Ac$. 
We say that a continuously differentiable function $f$ is $L_f$-smooth if there exists a constant $L_f \in [0, +\infty)$ such that $\Vert\nabla{f}(x) - \nabla{f}(y)\Vert_2 \leq L_f\Vert x - y\Vert_2$ for all $x, y\in\dom{f}$.
Here, $L_f$ is called a Lipschitz constant of $f$.
A function $f$ is said to be $\mu_f$-strongly convex if $f(\cdot) - \frac{\mu_f}{2}\norm{\cdot}^2_2$ remains convex.
If $\mu_f = 0$, then $f$ is just convex.

\subsection{Optimality condition and basic assumptions}
We define  $X := [U, V]$ as the joint variable of $U$ and $V$.  
We assume that $\phi$ in \eqref{eq:LMA_prob} is smooth. 
The \textit{optimality condition} of \eqref{eq:LMA_prob} can be written as follows:
\begin{equation}\label{eq:opt_cond}
\left\{\begin{array}{ll}
U_{\star}^{\top}\Ac^{*}\left(\nabla{\phi}(\Ac(U_{\star}V_{\star}^{\top}) - B)\right) &= 0, \vspace{0.75ex}\\
\Ac^{*}\left(\nabla{\phi}(\Ac(U_{\star}V_{\star}^{\top}) - B)\right)V_{\star} &= 0.
\end{array}\right.
\end{equation}
Any $X_{\star} = [U_{\star}, V_{\star}]$ satisfying \eqref{eq:opt_cond} is called a \textit{stationary point} of \eqref{eq:LMA_prob}.  
We denote by $\Xc_{\star}$ the set of stationary points of \eqref{eq:LMA_prob}.
Since $r \leq \min\set{m,n}$, the solution  of \eqref{eq:opt_cond} is generally nonunique.
Our aim is to design algorithms for generating a sequence $\set{X_k}$ converging to  $X_{\star}\in\Xc_{\star}$ under the following assumptions.

\begin{assumption}\label{as:A1}
Problem \eqref{eq:LMA_prob} satisfies the following conditions:
\begin{itemize}
\item[$\mathrm{(a)}$] The set $\Xc_{\star}$ of stationary points of \eqref{eq:LMA_prob} is nonempty, and $\Phi^{\star} > -\infty$ in \eqref{eq:LMA_prob}.
\item[$\mathrm{(b)}$]~$\phi$ is $L_{\phi}$-smooth and $\mu_{\phi}$-convex with $0 \leq \mu_{\phi} \leq L_{\phi} < +\infty$.
\end{itemize}
\end{assumption}
We allow $\mu_{\phi} = 0$, which also covers the non-strongly convex case.
Since $\phi$ is smooth and $\Ac$ is linear, $\Phi$ in \eqref{eq:LMA_prob} is also smooth. 
Moreover, as shown in \cite{Nesterov2004}, $\phi$ satisfies
\begin{equation}\label{eq:key_ineq1}
\frac{\mu_{\phi}}{2}\norm{y - x}_2^2 \leq \phi(y) - \phi(x) - \iprods{\nabla{\phi}(x), y-x} \leq\frac{L_{\phi}}{2}\Vert y - x\Vert_2^2, \quad \forall x, y\in\dom{\phi}.
\end{equation}
Note that Assumption A.\ref{as:A1}(b) covers a wide range of applications, including logistic loss, Huber loss,  and entropy function in statistics and machine learning \cite{Boyd2011}.

\section{Linesearch Gauss-Newton method}\label{sec:GN_method}
In this section, we develop a linesearch  Gauss-Newton (Ls-GN) algorithm for solving \eqref{eq:LMA_prob}. 

\subsection{Forming a surrogate of the objective}
By Assumption A.\ref{as:A1}, it follows from \eqref{eq:key_ineq1} and $\Phi(U,V) := \phi(\Ac(UV^{\top}) - B)$ that
\begin{align}\label{eq:surrogate}
\Phi(\hat{U},\hat{V}) \leq \Phi(U,V)  + \frac{L_{\Phi}}{2}\Vert \hat{U}\hat{V}^{\top} - (UV^{\top} - L_{\Phi}^{-1}\Phi^{\prime}(UV^{\top}))\Vert_F^2 - \frac{1}{2L_{\Phi}}\Vert \Phi^{\prime}(UV^{\top}) \Vert_F^2,
\end{align}
for any $U$, $V$, $\hat{U}$, and $\hat{V}$, 
where $\Phi^{\prime}(UV^{\top}) := \Ac^{*}\nabla{\phi}(\Ac(UV^{\top}) - B)$, and $L_{\Phi} := L_{\phi}\Vert\Ac\Vert^2$ is the Lipschitz constant of the gradient of $\phi(\Ac(\cdot)-B)$.

Gradient descent-type methods rely on finding a descent direction of $\Phi$ by  approximately minimizing the right-hand side surrogate of $\Phi$ in \eqref{eq:surrogate}. 
Unfortunately, this surrogate remains nonconvex due to the bilinear term $\hat{U}\hat{V}^{\top}$. 
Our next step is to linearize this   term around a given point $[U, V]$ as follows:
\begin{equation}\label{eq:bilinear_approx}
\hat{U}\hat{V}^{\top} \approx UV^{\top} + U(\hat{V} - V)^{\top} + (\hat{U} - U)V^{\top}.
\end{equation}
Then, the minimization of the right-hand side of \eqref{eq:surrogate} is approximated by 
\begin{equation}\label{eq:subprob1}
\min_{\hat{U},\hat{V}}\set{ (1/2)\Vert U(\hat{V} - V)^{\top} + (\hat{U} - U)V^{\top} + L_{\Phi}^{-1}\Phi^{\prime}(UV^{\top}) \Vert_F^2 }.
\end{equation}
This is  a linear least-squares problem, and can be solved by standard linear algebra routines.

\subsection{Computing Gauss-Newton direction}
Let us define
\begin{equation*}
D_U := \hat{U} - U, \quad D_V := \hat{V} - V, \quad \text{and} \quad Z :=  -L_{\Phi}^{-1}\Ac^{\ast}\left(\nabla{\phi}(\Ac(UV^{\top}) - B)\right).
\end{equation*} 
Then, we rewrite \eqref{eq:subprob1} as
\begin{equation}\label{eq:subprob2}
\min_{D_U, D_V}\set{ (1/2)\Vert UD_V^{\top} +D_UV^{\top} - Z\Vert_F^2 \ : \ D_U\in\R^{m\times r}, D_V \in\R^{m\times r}}.
\end{equation}
The optimality condition of \eqref{eq:subprob2} becomes
\begin{equation}\label{eq:opt_cond_subprob2}
\arraycolsep=0.2em
\left\{\begin{array}{llcl}
& U^{\top}UD_V^{\top} + U^{\top}D_UV^{\top} & = & U^{\top}Z, \vspace{1ex}\\
& UD_V^{\top}V + D_UV^{\top}V &= & ZV.
\end{array}\right.
\end{equation}
As usual, we  refer to \eqref{eq:opt_cond_subprob2} as the normal equation of \eqref{eq:subprob2}.
We will construct a closed form solution of \eqref{eq:opt_cond_subprob2}  in Lemma \ref{le:gauss_newton_dir}, whose proof is in Appendix \ref{apdx:le:gauss_newton_dir}.

\begin{lemma}\label{le:gauss_newton_dir}
The rank of the square linear system \eqref{eq:opt_cond_subprob2} does not exceed $r(m+n-r)$.
In addition, \eqref{eq:opt_cond_subprob2} has a solution.
If  $\rank{U} = \rank{V} = r \leq \min\set{m, n}$, then the solution of  \eqref{eq:opt_cond_subprob2} is given explicitly by 
\begin{equation}\label{eq:DuDv}
\left\{\begin{array}{lcl}
D_U &= & P^{\perp}_UZ(V^{\dagger})^{\top} + U\hat{D}_r, \vspace{1ex}\\
D_V^{\top} & = & U^{\dagger}Z -  \hat{D}_rV^{\top},
\end{array}\right.
\end{equation}
which forms a linear subspace in $\R^{r\times r}$, and $\hat{D}_r\in\R^{r\times r}$ is an arbitrary matrix.

In particular, if we choose $\hat{D}_r := 0.5 U^{\dagger}Z(V^{\dagger})^{\top} \in\R^{r\times r}$, then
\begin{equation}\label{eq:gn_dir}
D_U = \big(\Id_m - 0.5 P_U\big)Z(V^{\dagger})^{\top} \quad \text{and} \quad D_V^{\top} = U^{\dagger}Z\big(\Id_n -  0.5 P_V\big).
\end{equation}
Moreover, the optimal value of \eqref{eq:subprob2} is  $(1/2)\Vert P^{\perp}_UZP^{\perp}_V\Vert_F^2$.
\end{lemma}

Lemma \ref{le:gauss_newton_dir} also shows that if either $Z$ is in the null space of $P_U$ or $Z^{\top}$ is in the null space of $P_V$, then $\Vert P^{\perp}_UZP^{\perp}_V\Vert_F^2 = 0$. 
Since  \eqref{eq:gn_dir} only gives us one choice for $D_{X} := [D_U, D_V]$, if $\hat{D}_r = \boldsymbol{0}^r$,   we  obtain another simple GN  search direction. 

\begin{remark}\label{re:symmetric_case}
Let $m = n$. If we assume that $U = V$, then $D_U = D_V$ and 
\begin{equation*}
D_U =  P^{\perp}_UZ(U^{\dagger})^{\top} + U\hat{D}_r,  \quad \text{where} \quad \hat{D}_r \in \Sc_r := \set{ \hat{D}_r\in\R^{r\times r} \ : \ \hat{D}_r + \hat{D}_r^{\top} = U^{\dagger}Z(U^{\dagger})^{\top} }.
\end{equation*}
Clearly, $\Sc_r$ is a linear subspace, and its dimension is $r(r+1)/2$.
\end{remark}

\subsection{The damped-step Gauss-Newton scheme}
Using  Lemma \ref{le:gauss_newton_dir}, we can form a damped step GN scheme as follows:
\begin{equation}\label{eq:gau_newton_scheme}
\left\{\begin{array}{lcl}
U_{+} & := & U + \alpha D_U, \vspace{1ex}\\
V_{+} & :=  & V + \alpha D_V,
\end{array}\right.
\end{equation}
where $D_U$ and $D_V$ defined in \eqref{eq:gn_dir} is a GN direction, and $\alpha > 0$ is a given step-size determined in the next lemma.

Since the GN direction computed from \eqref{eq:subprob2} is not unique, we need to choose an appropriate $D_{X}$ such that it is a descent direction of $\Phi$ at $X$.
We prove in Lemma~\ref{le:descent_dir} that \eqref{eq:gn_dir} indeed gives a descent direction of  $\Phi$ at $X$. The proof of this lemma is deferred to Appendix \ref{apdx:le:descent_dir}.

\begin{lemma}\label{le:descent_dir}
Let $X := [U, V]$ be a non-stationary point of \eqref{eq:LMA_prob} and $D_X := [D_U, D_V]$ be given by  \eqref{eq:gn_dir}. 
If $D_X\neq 0$ and  $\alpha$ is chosen in $0 < \alpha \leq \underline{\alpha}$ where
\begin{equation}\label{eq:step_size_min}
\underline{\alpha} := \min\set{1,  \frac{L_{\Phi}\sigma_{\min}^3}{2\Vert\nabla{\Phi}(U,V)\Vert_F}, \frac{3\sigma_{\min}^4}{32\sigma_{\max}^2\Vert \Phi^{\prime}(UV^{\top}) \Vert_F}} \in (0, 1],
\end{equation}
then we have
\begin{equation}\label{eq:descent_property}
\Phi(U_{+}, V_{+}) \leq \Phi(U, V) - \frac{\alpha\sigma_{\min}^2}{128L_{\Phi}\sigma_{\max}^4}\Vert\nabla{\Phi}(U,V)\Vert^2,
\end{equation}
where $\Phi^{\prime}(\cdot) = \Ac^{*}\nabla{\phi}(\Ac(\cdot) - B)$, $L_{\Phi} := L_{\phi}\norm{\Ac}^2$, $\sigma_{\min} := \min\set{\sigma_{\min}(U), \sigma_{\min}(V)}$, and $\sigma_{\max} := \max\set{\sigma_{\max}(U), \sigma_{\max}(V)}$.
Hence,  $D_X$ is a descent direction of  $\Phi$.
\end{lemma}

Lemma \ref{le:descent_dir} shows that if the residual term $ \Ac^{*}\nabla{\phi}(\Ac(UV^{\top}) - B)$ is sufficient small near $\Xc_{\star}$, then we obtain a full-step size $\alpha = 1$. 

The existence of the GN direction in Lemma \ref{le:gauss_newton_dir} requires $U$ and $V$ to be full-rank. We prove in   Appendix \ref{apdx:le:full_rank_uv} the following lemma.

\begin{lemma}\label{le:full_rank_uv}
If $\rank{U} = \rank{V} = r$, then $X_{+} := [U_{+}, V_{+}]$ updated by \eqref{eq:gau_newton_scheme} using the step-size $\underline{\alpha}$ in \eqref{eq:step_size_min} satisfies
\begin{equation}\label{eq:rank_preserve}
\sigma_{\min}(U_{+}) \geq 0.5\sigma_{\min}(U) \quad \text{and} \quad \sigma_{\min}(V_{+}) \geq 0.5\sigma_{\min}(V).
\end{equation} 
Hence, \eqref{eq:gau_newton_scheme} preserves the rank of $U_{+}$ and $V_{+}$, i.e., $\rank{U_{+}} = \rank{V_{+}} = r$.
\end{lemma}

\subsection{The algorithmic template and its global convergence}
Theoretically, we can use the step-size $\underline{\alpha}$ in Lemma \ref{le:descent_dir} for \eqref{eq:gau_newton_scheme}. 
However, in practice, computing $\underline{\alpha}$ requires a high computational cost. 
We instead incorporate the GN scheme \eqref{eq:gau_newton_scheme} with an Armijo  backtracking linesearch to find an appropriate step-size $\alpha \geq \beta\underline{\alpha}$ for a given $\beta \in (0, 1)$.
\begin{framed}
\vspace{-1ex}
\noindent Find the smallest integer number $i_k\geq 0$ such that $\alpha := \beta^{i_k}\alpha_0 \geq \underline{\alpha}$ and 
\begin{equation}\label{eq:ls_condition}
\Phi(U + \alpha D_U, V + \alpha D_V) \leq \Phi(U, V) - 0.5c_1\alpha\Vert\nabla{\Phi}(U, V)\Vert_F^2,
\end{equation}
where $\alpha_0 > 0$, $c_1 > 0$, and $\beta\in (0, 1)$ are given (e.g., $c_1 := 0.5$ and $\beta := \frac{\sqrt{5}-1}{\sqrt{5} + 1}$).
\vspace{-1ex}
\end{framed}
By Lemma \ref{le:descent_dir}, this   procedure is terminated after a finite number of iterations $i_k$ such that
\begin{equation}\label{eq:i_k}
0 \leq i_k \leq \lfloor \log_{\beta}(\underline{\alpha}/\alpha_0) \rfloor + 1,
\end{equation}
where $\underline{\alpha}$ is given by \eqref{eq:step_size_min}.
Now, we describe the complete linesearch GN algorithm for approximating a stationary point of \eqref{eq:LMA_prob} as in Algorithm \ref{alg:A1}.

\begin{algorithm}[!ht]\caption{(\textit{Linesearch Gauss-Newton Algorithm} (Ls-GN))}\label{alg:A1}
\begin{algorithmic}[1]
   \STATE {\bfseries Initialization:} Given a tolerance $\varepsilon > 0$. Choose $X_0 := [U_0, V_0]$. Set $c_1 := 0.5$ and $\alpha_0 := 1$.
   \FOR{$k = 0$ {\bfseries to} $k_{\max}$}
	\STATE\label{step:gn_dir} \textit{GN direction}: Let $Z_k := -L_{\Phi}^{-1}\Phi^{\prime}(U_kV_k^{\top})$. 
	Compute  $D_{X_k} := [D_{U_k}, D_{V_k}]$:
	\begin{equation*}
	D_{U_k} := \left(\Id_m - 0.5 P_{U_k}\right)Z_k(V_k^{\dagger})^{\top} \quad \text{and} \quad D_{V_k} = \left(\Id_n - 0.5 P_{V_k}\right)Z_k^{\top}(U_k^{\dagger})^{\top}.
	\vspace{-2ex}
	\end{equation*}
	\STATE \textit{Stopping criterion}: If \texttt{stopping\_criterion}, then TERMINATE.
	\STATE\label{step:linesearch} \textit{Backtracking linesearch}:  Find the smallest integer number $i_k \geq 0$ such that
	\begin{equation*} 
		\Phi(U_k+ \alpha_kD_{U_k}, V_k + \alpha_kD_{V_k}) \leq \Phi(U_k, V_k) - 0.5c_1\alpha_k\Vert\nabla{\Phi}(U_k, V_k)\Vert_F^2,
	\end{equation*}
	where $\alpha_k := \alpha_0\beta^{i_k}$.
	\STATE Update  $X_{k+1} := [U_{k+1}, V_{k+1}]$ as $U_{k+1} := U_k + \alpha_kD_{U_k}$ and $V_{k+1} := V_k + \alpha_kD_{V_k}$.
   \ENDFOR
\end{algorithmic}
\end{algorithm}

\paragraph{\textbf{Per-teration complexity:}}
The main steps of Algorithm \ref{alg:A1} are Steps \ref{step:gn_dir} and \ref{step:linesearch}, i.e. computing $D_{X_k}$ and performing the linesearch routine, respectively. 
\begin{itemize}
\item[$\mathrm{(a)}$]~Computing $D_{X_k}$ requires two inverses $(U^{\top}U)^{-1}$ and $(V^{\top}V)^{-1}$ of the size $r\times r$, and two matrix-matrix multiplications (of the size $m\times r$ or $n\times r$). 

\item[$\mathrm{(b)}$]~Evaluating $\Phi^{\prime}$ requires one matrix-matrix multiplication $UV^{\top}$ and one evaluation of the form $\Ac^{*}\nabla{\phi}(\Ac(\cdot) - B)$. 
When $\Ac$ is a subset projection $\mathcal{P}_{\Omega}$ (e.g., in matrix completion), we can compute $(UV^{\top})_{(i,j)\in\Omega}$ instead of the full matrix $UV^{\top}$. 

\item[$\mathrm{(c)}$]~Each step of the linesearch needs one matrix-matrix multiplication $UV^{\top}$ and one evaluation of $\Phi$. 
It requires at most $\lfloor \log_{\beta}(\underline{\alpha}/\alpha_0) \rfloor + 1$ linesearch iterations.
However, we observe that $i_k$ often varies from $1$ to $2$ on average in our experiments in Section~\ref{sec:num_experiments}. 
\end{itemize}

\paragraph{\textbf{Global convergence:}}
Since  \eqref{eq:LMA_prob} is nonconvex, we only expect $\set{X_k}$ generated by  Algorithm \ref{alg:A1} to converge to a stationary point $X_{\star}\in\Xc_{\star}$.
However, Lemma \ref{le:full_rank_uv} only guarantees the full-rankness of $U_k$ and $V_k$ at each iteration, but we may have 
$\lim\limits_{k\to\infty}\sigma_{\min}(U_k) = 0$ or $\lim\limits_{k\to\infty}\sigma_{\min}(V_k) = 0$. 
In order to prove a global convergence of Algorithm \ref{alg:A1}, we require one additional condition:
There exists   $\underline{\sigma} > 0$ such that:
\begin{equation}\label{eq:bounded_assumption}
\sigma_{\min}(U_k) \geq \underline{\sigma} \quad \text{and} \quad \sigma_{\min}(V_k) \geq \underline{\sigma} \quad \text{for all}~k\geq 0.
\end{equation}
Under Assumption A.\ref{as:A1}, the following sublevel set of $\Phi$:
\begin{equation*}
\mathcal{L}_{\Phi}(\gamma) := \set{ [U, V]\in\dom{\Phi} \ : \ \Phi(U, V) \leq \gamma}
\end{equation*}
 is bounded for a given $\gamma > 0$.
We  prove in Appendix \ref{apdx:th:convergence_guarantee1} a global convergence of Algorithm \ref{alg:A1} stated in the following theorem.

\begin{theorem}\label{th:convergence_guarantee1}
Let $\set{X_k}$ with $X_k := [U_k, V_k]$  be generated by Algorithm \ref{alg:A1}. 
Then, under Assumption A.\ref{as:A1}, we have
\begin{equation}\label{eq:sum_inf}
\sum_{k=0}^{\infty}\alpha_k\Vert \nabla{\Phi}(U_k, V_k)\Vert_F^2 <+\infty, ~~\text{and}~~ \lim_{k\to\infty} \alpha_k\Vert \nabla{\Phi}(U_k, V_k)\Vert_F = 0.
\end{equation}
If, in addition, the condition \eqref{eq:bounded_assumption} holds and $\set{X_k}$ is bounded, then 
\begin{equation}\label{eq:lim_inf}
\lim_{k\to\infty} \Vert \nabla{\Phi}(U_k, V_k)\Vert_F = 0.
\end{equation}
There exists a limit point $X_{\star}$ of $\set{X_k}$, and any limit point $X_{\star}$ is in $\Xc_{\star}$.
\end{theorem}

\subsection{Local linear convergence without strong convexity}
We prove a local convergence of the full-step Gauss-Newton scheme \eqref{eq:gau_newton_scheme} when $\alpha = 1$.
Generally, problem \eqref{eq:LMA_prob} does not satisfy the regularity assumption: the Jacobian $J_R(X) = A[V \otimes \Id_m, \Id_n\otimes U]\in\R^{l\times (m+n)r}$ of the objective residual $R(X) := \Ac(UV^{\top}) - B$ in \eqref{eq:LMA_prob} is not full-column rank, where $A$ is the matrix form of the linear operator $\Ac$.
However, we can still guarantee a fast local convergence under the following conditions:

\begin{assumption}\label{as:A2} Problem \eqref{eq:LMA_prob} satisfies the following conditions:
\begin{itemize}
\item[$\mathrm{(a)}$]~$\phi$ is twice continuously differentiable on a neighborhood $\Nc(Z_{\star})$ of $Z_{\star} = \Ac(U_{\star}V^{\top}_{\star}) -B$, and its Hessian $\nabla^2{\phi}$ is Lipschitz continuous in $\Nc(Z_{\star})$ with the constant $L_{\phi^{''}}$.

\item[$\mathrm{(b)}$]~The Hessian $\nabla^2{\Phi}(X_{\star})$ of $\Phi(X) := \phi(\Ac(UV^{\top}) - B)$ at $X_{\star}\in\Xc_{\star}$ satisfies
\begin{equation}\label{eq:small_residual}
\Vert \left[\Id - L_{\Phi}^{-1}H(X)^{\dagger}\nabla^2{\Phi}(X_{\star})\right](X - X_{\star})\Vert_F \leq \kappa(X_{\star})\Vert X - X_{\star}\Vert_F, \quad \forall X\in\Nc(X_{\star}),
\end{equation}
where $H(X) := \begin{bmatrix}V^{\top}\otimes U & V^{\top}V \otimes \Id_m \\ \Id_n\otimes U^{\top}U & V\otimes U^{\top}\end{bmatrix}$, $L := L_{\phi}\Vert\Ac\Vert^2$, and $0 \leq \kappa(X_{\star}) \leq \bar{\kappa} < 1$.
\end{itemize}
\end{assumption}

Assumption~A.\ref{as:A2}(b) relates to a ``small residual condition''.
For instance, if  $\phi(\cdot) = (1/2)\Vert\cdot\Vert_2^2$,  and $\Ac = \Id$, the identity operator, then the residual term becomes $R(X) = UV^{\top}-B$, and $\Phi(X) = (1/2)\Vert R(X)\Vert^2_F$.
In this case, condition \eqref{eq:small_residual} holds if $\Vert R(X_{\star})\Vert_F \leq \kappa(X_{\star}) < 1$ (i.e., we have a ``small residual'' case).

Now, we prove in Appendix \ref{apdx:th:local_convergence} a local convergence of the full-step GN variant.  

\begin{theorem}\label{th:local_convergence}
Let $\set{X_k}$ be generated by \eqref{eq:gau_newton_scheme} with a full step-size $\alpha_k = 1$, and $X_{\star} := [U_{\star}, V_{\star}]\in\Xc_{\star}$ be a given stationary point of \eqref{eq:LMA_prob} such that $\rank{U_{\star}} = \rank{V_{\star}} = r$. Assume that Assumptions A.\ref{as:A1} and A.\ref{as:A2} hold.
Then, there exists a neighborhood $\Nc(X_{\star})$ of $X_{\star}$ and a constant $K_1 > 0$ independent of $X_k$ such that 
\begin{equation}\label{eq:gn_local_est1}
\Vert X_{k+1} - X_{\star}\Vert_F \leq \big(\bar{\kappa} + 0.5K_1\Vert X_k - X_{\star}\Vert_F \big)\Vert X_k - X_{\star}\Vert_F, \quad \forall X_k\in\Nc(X_{\star}).
\end{equation}
Consequently, if $H(X_{\star})^{\dagger}\nabla^2{\Phi}(X_{\star}) = L_{\Phi}\Id$ in \eqref{eq:small_residual}  $($i.e., zero residual$)$, then there exists a constant $K_2 > K_1$ such that the sequence $\set{X_k}$ generated by our full-step GN algorithm starting from $X_0\in\Nc(X_{\star})$ with $\Vert X_0-X_{\star}\Vert_F < 2K_2^{-1}$ quadratically converges  to $X_{\star}\in\Xc_{\star}$.

If $\kappa(X_{\star}) \in (0, 1)$ in \eqref{eq:small_residual} $($i.e., small residual$)$, then, for any $X_0\in\Nc(X_{\star})$ such that $\Vert X_0-X_{\star}\Vert_F \leq \bar{r}_0 < 2K_1^{-1}(1-\bar{\kappa})$, $\set{X_k}$ linearly converges to $X_{\star}$.
\end{theorem}

\section{ADMM-Gauss-Newton Algorithm}\label{sec:admm_gn_method}
The GN method only works well and has a fast local  convergence for the ``small residual'' case. 
In general, it may converge very slowly or even fails to converge. 
In this section, we propose to combine the GN scheme \eqref{eq:gau_newton_scheme} and the alternating direction method of multipliers (ADMM) to develop a  new algorithm for solving \eqref{eq:LMA_prob} called GN-ADMM. 
The ADMM can be viewed as a variant of augmented Lagrangian-based methods in nonlinear optimization \cite{Bertsekas1996a,Hestenes1969,polyak2009local}.
It can also be derived from Douglas-Rachford's method in convex optimization.

\subsection{The augmented Lagrangian function and  ADMM scheme}\label{subsec:AL_scheme}
We  introduce  $W = \Ac(UV^{\top}) - B$ and rewrite  \eqref{eq:LMA_prob} as the following  problem:
\begin{equation}\label{eq:constr_LMA_prob2}
\Phi_{\star} := \min_{U, V, W}\set{ \phi(W) \ : \ \Ac(UV^{\top})  - W = B }.
\end{equation}
We can define the augmented Lagrangian function associated  with \eqref{eq:constr_LMA_prob2} as
\begin{equation}\label{eq:aug_Lagrangian}
\begin{array}{lcl}
\Lc_{\rho}(U, V, W, \Lambda) & := & \phi(W) + \iprods{\Lambda, \Ac(UV^{\top}) - W-B} + \frac{\rho}{2}\Vert \Ac(UV^{\top}) - W-B\Vert_2^2{} \vspace{1ex}\\
& = &  \phi(W) + \frac{\rho}{2}\Vert \Ac(UV^{\top}) - W-B + \rho^{-1}\Lambda \Vert_2^2 - \frac{1}{2\rho}\Vert \Lambda\Vert_2^2,
\end{array}{}
\end{equation}
where $\rho > 0$ is a penalty parameter and $\Lambda$ is a Lagrange multiplier.

Next, we apply the standard ADMM scheme to \eqref{eq:constr_LMA_prob2} which leads to the following $3$ steps:
\begin{subequations}\label{eq:aug_method2}
\begin{eqnarray}
& &(U_{k+1}, V_{k+1})   := \displaystyle\argmin_{U, V}\set{ \Vert \Ac(UV^{\top}) - W_k - B + \rho^{-1}\Lambda_k \Vert_2^2 },\label{eq:aug_method2_a}\\
& &W_{k+1}   :=  \displaystyle\argmin_{W}\Big\{  \phi(W) + (\rho/2)\Vert W - \big(\Ac(U_{k+1}V_{k+1}^{\top}) - B + \rho^{-1}\Lambda_k \big) \Vert_2^2\Big\},{~~~~~~~}\label{eq:aug_method2_b}\\
& &\Lambda_{k+1}    :=  \Lambda_k + \rho(\Ac(U_{k+1}V_{k+1}^{\top}) - W_{k+1}-B).\label{eq:aug_method2_c}
\end{eqnarray}
\end{subequations}
Obviously, both subproblems \eqref{eq:aug_method2_a} and \eqref{eq:aug_method2_b} remain computationally expensive.
While \eqref{eq:aug_method2_a} is nonconvex,  \eqref{eq:aug_method2_b} is smooth and convex.
Without any further step applying to \eqref{eq:aug_method2}, convergence theory for this nonconvex ADMM scheme can be found in several recent papers including \cite{li2015global,wang2015global,wen2012solving}.
However, \eqref{eq:aug_method2} remains impractical since  \eqref{eq:aug_method2_a} and \eqref{eq:aug_method2_b} cannot be solved with a closed form or  a highly accurate solution. 
We approximately solve these subproblems.

\subsection{Approximation of the alternating steps}
We apply the GN scheme to approximate \eqref{eq:aug_method2_a} and a linearization to approximate \eqref{eq:aug_method2_b}  in our ADMM scheme above.

\paragraph{\textbf{Gauss-Newton step for the $UV$-subproblem \eqref{eq:aug_method2_a}:}}
We first apply on step of \eqref{eq:gau_newton_scheme} to solve \eqref{eq:aug_method2_a} as follows.
We first approximate $\Vert \Ac(UV^{\top}) -  W_k - B  +  \rho^{-1}\Lambda_k\Vert^2_2$ by using the quadratic surrogate of $\Ac(\cdot)$ and the linearization $U_kV_k^{\top} + U_kD_V^{\top} + D_UV_k^{\top}$ of $UV^{\top}$ with $D_U := U - U_k$ and $D_V := V - V_k$ as new variables.
By letting $Z_k := -L_{\Ac}^{-1}\Ac^{\ast}\left(\Ac(U_kV_k^{\top}) - W_k - B + \rho^{-1}\Lambda_k\right)$ with $L_{\Ac} :=\Vert\Ac\Vert^2$, we solve
\begin{equation}\label{eq:linearizeD_UV_subprob}
[D_{U_k}, D_{V_k}] := \mathrm{arg}\min_{D_U, D_V}\Big\{ \mathcal{Q}_k(D_U, D_V) := \frac{1}{2}\Vert U_kD_V^{\top} + D_UV_k^{\top} - Z_k\Vert_F^2 \Big\}.
\end{equation}
Here, the Lipschitz constant $L_{\Ac} :=\Vert\Ac\Vert^2$ can be computed by a power method \cite{Golub1996}.
Using Lemma \ref{le:gauss_newton_dir}, we can compute $[D_{U_k}, D_{V_k}]$ as
\begin{equation}\label{eq:gn_dir2}
\left\{\begin{array}{lcl}
D_{U_k} &:= & \big(\Id_m - 0.5P_{U_k}\big)Z_k(V_k^{\dagger})^{\top} \vspace{1ex}\\
D_{V_k}^{\top} &:= & U_k^{\dagger}Z_k\big(\Id_n - 0.5P_{V_k}\big).
\end{array}\right.
\end{equation}
The corresponding objective value is  $\mathcal{Q}_k(D_{U_k}, D_{V_k}) := (1/2)\Vert P^{\perp}_{u_k}Z_kP^{\perp}_{V_k}\Vert_F^2$.
Then, we update $X_{k+1} := [U_{k+1}, V_{k+1}]$ as  
\begin{equation}\label{eq:uv_update}
U_{k+1} := U_k + \alpha_kD_{U_k} \quad \text{and} \quad V_{k+1} := V_k + \alpha_kD_{V_k},
\end{equation}
where $\alpha_k > 0$ is a step-size computed by a linesearch procedure as in \eqref{eq:ls_condition}.

\vspace{0.5ex}
\paragraph{\textbf{Gradient step for the $W$-subproblem \eqref{eq:aug_method2_b}:}}
If $\phi$ does not have a tractably proximal operator (i.e., its proximal operator cannot be computed in a closed form, or with a low-order polynomial-time algorithm), we approximate \eqref{eq:aug_method2_b} by using one gradient step as
\begin{equation}\label{eq:aug_method2b}
W_{k+1} := \displaystyle\argmin_{W}\Big\{  \frac{L_{\phi}}{2}\Vert W - (W_k - L_{\phi}^{-1}\nabla{\phi}(W_k))\Vert_2^2 + \frac{\rho}{2}\Vert W - E_k \Vert_2^2\Big\},
\end{equation}
where $E_k := \Ac(U_{k+1}V_{k+1}^{\top})-B + \rho^{-1}\Lambda_k$.
Solve  \eqref{eq:aug_method2b} directly, we get
\begin{equation}\label{eq:second_subprob_sol}
W_{k+1} :=   (\rho + L_{\phi})^{-1}\left(L_{\phi}W_k - \nabla{\phi}(W_k) +  (\Lambda_k + \rho (\Ac(U_{k+1}^{\top}V_{k+1}) - B)  \right).
\end{equation}

\subsection{The ADMM-Gauss-Newton algorithm and its global convergence}\label{subsec:admm_gn_alg}
Putting \eqref{eq:uv_update}, \eqref{eq:aug_method2_c}, and  \eqref{eq:aug_method2_b} or \eqref{eq:second_subprob_sol} together, we obtain the following ADMM-GN scheme with two options:
\begin{equation}\label{eq:aug_method3}
\left\{\begin{array}{ll}
Z_k & := -L_{\Ac}^{-1}\Ac^{\ast}\left(\Ac(U_kV_k^{\top})  - W_k - B + \rho^{-1}\Lambda_k\right), \vspace{1ex}\\
U_{k+1} &:=  U_k + \alpha_k\big(\Id_m - 0.5P_{U_k}\big)Z_k(V_k^{\dagger})^{\top}, \vspace{1ex}\\
V_{k+1} &:=  V_k + \alpha_k\big(\Id_n - 0.5P_{V_k}\big)Z_k(U_k^{\dagger})^{\top}, \vspace{1ex}\\
W_{k+1} &~\text{is computed by \eqref{eq:aug_method2_b} for \textbf{Option 1}, or by \eqref{eq:aug_method2b} for \textbf{Option 2}},\vspace{1ex}\\
\Lambda_{k+1} & := \Lambda_k + \rho(\Ac(U_{k+1}V_{k+1}^{\top}) - W_{k+1} - B).
\end{array}\right.
\end{equation}
Clearly, computing $[U_{k+1}, V_{k+1}]$ in \eqref{eq:aug_method3} using the step-size in Lemma \ref{le:descent_dir} is impractical. 
Similar to Algorithm \ref{alg:A1}, we find an appropriate $\alpha_k$ by  a backtracking linesearch on $\mathcal{Q}_k(U, V) := (1/2)\norm{\Ac(UV^{\top}) - W_k - B + \rho^{-1}\Lambda_k}_2^2$ as
\begin{equation}\label{eq:ls_cond2}
\mathcal{Q}(U_k + \alpha D_{U_k}, V_k + \alpha_kV_k) \leq \mathcal{Q}(U_k, V_k) -  0.5c_1\alpha_k\Delta_k^2, 
\end{equation}
where $\Delta_k^2 := \Vert U_k^{\top}\Ac^{\ast}(E_k -W_k)\Vert_F^2 + \Vert \Ac^{\ast}(E_k-W_k)V_k \Vert_F^2$ and $\alpha_k := \beta^{i_k}\alpha_0$  with $\alpha_0 > 0$ and $\beta := (\sqrt{5}-1)/(\sqrt{5}+1) \in (0, 1)$ given a priori.
Obviously, by Lemma  \ref{le:descent_dir}, this   procedure terminates after a finite number of linesearch steps $i_k$ satisfying \eqref{eq:i_k}. 
In addition, $D_{X_k} := [D_{U_k}, D_{V_k}]$ is a descent direction of the quadratic objective $\mathcal{Q}_k$ at $X_k$.

Now, we expand  \eqref{eq:aug_method3} algorithmically  as in Algorithm \ref{alg:A2}.

\begin{algorithm}[!ht]\caption{(\textit{ADMM-Gauss-Newton Algorithm} (ADMM-GN))}\label{alg:A2}
\begin{algorithmic}[1]
   \STATE {\bfseries Initialization:} Given $\varepsilon > 0$, choose $\rho > 0$ and  $X_0 := [U_0,V_0]$. 
    \STATE\hspace{2ex} Set $W_0 := U_0V_0^{\top}$ and $\Lambda_0 := \boldsymbol{0}^{m\times n}$.
   \FOR{$k = 0$ {\bfseries to} $k_{\max}$}
	\STATE\label{step:gnadmm_gn} \textit{Gauss-Newton step}: Compute a GN direction $D_{X_k} := [D_{U_k}, D_{V_k}]$ by \eqref{eq:gn_dir2}.
	\STATE\label{step:gnadmm_ls}  \textit{Linesearch step}: Find $\alpha_k > 0$ from the linesearch condition \eqref{eq:ls_cond2} and update
	\vspace{-1ex}
	\begin{align*}
	U_{k+1} := U_k + \alpha_kD_{U_k} \quad \text{and} \quad V_{k+1} := V_k + \alpha_k D_{V_k}.
	\end{align*}
	\vspace{-3ex}
	\STATE\label{step:gnadmm_grad} \textit{Gradient step}: Evaluate $Y_{k+1} :=  \Ac(U_{k+1}V_{k+1}^{\top})-B$, and $\Phi^{\prime}(W_k)$,  and
	\begin{equation*} 
		\text{\textbf{Option 1}: update $W_{k+1}$ by \eqref{eq:aug_method2_b}} \quad \text{or} \quad \text{\textbf{Option 2}: update $W_{k+1}$ by \eqref{eq:aug_method2b}} 
	\end{equation*}
	\vspace{-3ex}
	\STATE  If \texttt{stopping\_criterion}, then TERMINATE.
	\STATE\label{step:gnadmm_ud} Update $\Lambda_{k+1}  := \Lambda_k + \rho(Y_{k+1} - W_{k+1})$.
   \ENDFOR
\end{algorithmic}
\end{algorithm}

\paragraph{\textbf{Per-iteration complexity:}}
The main steps of Algorithm \ref{alg:A2} remain at Steps \ref{step:gnadmm_gn} and \ref{step:gnadmm_ls}, where they require to compute $D_{X_k} := [D_{U_k}, D_{V_k}]$ and to perform a linesearch procedure, respectively. 
Steps \ref{step:gnadmm_grad} and \ref{step:gnadmm_ud} only require matrix-matrix additions which have the complexity of $\mathcal{O}(m\times n)$.
Overall, the per-iteration complexity of Algorithm \ref{alg:A2} is higher than of Algorithm \ref{alg:A1}, but as we can see from Section \ref{sec:num_experiments} that  we can simply use the full-step GN scheme at Step~\ref{step:gnadmm_gn} without linesearch, and Algorithm \ref{alg:A2} often requires a fewer number of iterations than Algorithm \ref{alg:A1}. 
Moreover, Algorithm \ref{alg:A2} seems working well for the ``large residual case'', i.e., $\Ac^{\ast}\nabla{\phi}(\Ac(U_{\star}V_{\star}^{\top}) - B)$ is large.

\vspace{0.5ex}
\paragraph{\textbf{Global convergence analysis:}}\label{subsec:convergence_analysis}
We first write the optimality condition (or the KKT condition) for \eqref{eq:constr_LMA_prob2} as follows:
\begin{equation}\label{eq:kkt_cond}
\nabla{\phi}(W_{\star}) - \Ac^{\ast}(\Lambda_{\star}) = 0,~~U_{\star}^{\top}\Ac^{\ast}(\Lambda_{\star}) = 0, ~~\Ac^{\ast}(\Lambda_{\star}) V_{\star} = 0, ~\text{and}~\Ac(U_{\star}V_{\star}^{\top}) - W_{\star} = B.{}
\end{equation}
This condition can be rewritten as \eqref{eq:opt_cond} by eliminating $W_{\star}$ and the multiplier $\Lambda_{\star}$. 
Hence, if $[U_{\star}, V_{\star}, W_{\star}, \Lambda_{\star}]$ satisfies \eqref{eq:kkt_cond}, then $X_{\star} := [U_{\star}, V_{\star}]\in\Xc_{\star}$.

The following lemma provides a key step to prove the convergence of  Algorithm \ref{alg:A2},  whose proof is given in Appendix \ref{apdx:le:descent2}.

\begin{lemma}\label{le:descent2}
Let $\set{[U_k, V_k, W_k, \Lambda_k]}$ be generated by Algorithm \ref{alg:A2}. 
Suppose that Assumption~A.\ref{as:A1} holds and $\set{(U_k, V_k)}$ is bounded.
Then, the following statements hold:
\begin{itemize}
\item[$\mathrm{(a)}$]~The sequence $\set{ (W_k, \Lambda_k) }$ is bounded. In addition, for $k\geq 1$, we have
\begin{equation}\label{eq:bound_lambda}
\begin{array}{llll}
&\Vert\Lambda_{k+1} - \Lambda_k\Vert_2 &\leq L_{\phi}\Vert W_{k+1} - W_k\Vert_2 &{} \text{for \textbf{Option 1}},\vspace{1ex}\\
\text{or} & \Vert\Lambda_{k+1} - \Lambda_k\Vert_2 &\leq L_{\phi}\big(\Vert W_k - W_{k-1}\Vert_2 + \Vert W_{k+1} - W_{k-1}\Vert_2\big) &\text{for \textbf{Option 2}}.
\end{array}
\vspace{1ex}
\end{equation}
\item[$\mathrm{(b)}$]~Let $\Lc_{\rho}$ be defined by \eqref{eq:aug_Lagrangian}.
Then, for any $\rho > 0$, we have
\begin{equation}\label{eq:descent2}
\begin{array}{ll}
&\Lc_{\rho}(U_{k+1}, V_{k+1}, W_{k+1}, \Lambda_{k+1}) \leq \Lc_{\rho}(U_k,V_k,W_k, \Lambda_k) - 
\frac{\eta_1}{2}\Vert W_{k+1} - W_k\Vert_2^2 \vspace{1ex}\\
&\qquad\quad + {~} \frac{\eta_0}{2}\Vert W_k - W_{k-1}\Vert_2^2 -  \frac{c_1\rho\alpha_k}{2}\left[\Vert U_k^{\top}\Ac^{\ast}(E_k - W_k)\Vert_F^2 + \Vert \Ac^{\ast}(E_k -W_k)V_k \Vert_F^2\right],
\end{array}
\end{equation}
where $E_k := \Ac(U_kV_k^{\top})-B + \rho^{-1}\Lambda_k$, and 
\begin{equation}\label{eq:descent2b}
\begin{array}{llllll}
&\eta_1 &:= \rho^{-1}\big(\rho^2 + \mu_{\phi}\rho - 2L_{\phi}^2\big) ~~&\text{and} &\eta_0 := 0 &\text{for \textbf{Option 1}},\vspace{1ex}\\
\text{or}~&\eta_1 &:= \rho^{-1}\big(\rho^2 + L_{\phi}\rho - 4L_{\phi}^2\big)~~&\text{and} &\eta_0 := 8\rho^{-1}L_{\phi}^2 &\text{for \textbf{Option 2}}.
\end{array}
\end{equation}
\end{itemize}
\end{lemma}
Similar to Algorithm \ref{alg:A1}, we prove a global convergence of Algorithm \ref{alg:A2} in the following theorem, whose proof is deferred to  Appendix \ref{apdx:th:global_convergence2}.

\begin{theorem}\label{th:global_convergence2}
Under Assumption A.\ref{as:A1} and condition \eqref{eq:bounded_assumption}, let  $\set{[U_k, V_k]}$   generated by Algorithm~\ref{alg:A2} be bounded. Then, if we choose $\rho$ such that
\begin{equation}\label{eq:choice_of_rho}
\left\{\begin{array}{lll}
\rho &> 0.5\big(\big(\mu_{\phi} + 8L_{\phi}^2\big)^{1/2} + \mu_{\phi}\big) &\text{for \textbf{Option 1}},\vspace{1ex}\\
\rho &> 3L_{\phi} &\text{for \textbf{Option 2}}, 
\end{array}\right.
\end{equation}
then
\begin{equation}\label{eq:norm_grad_limit}
\lim_{k\to\infty}\Vert \nabla{\Phi}(U_k, V_k)\Vert_F = 0.
\end{equation}
Consequently,  there exists a limit point $X_{\star} := [U_{\star}, V_{\star}]$ of $\set{[U_k, V_k]}$ and $X_{\star} \in \Xc_{\star}$.
\end{theorem}

\section{Symmetric low-rank matrix optimization}\label{sec:symmetric_case}
In this section, we develop a symmetric GN variant of Algorithm~\ref{alg:A1} for solving the following special symmetric setting of \eqref{eq:LMA_prob} when $U = V$:
\begin{equation}\label{eq:sym_LMA_prob}
\Phi^{\star} := \min_{U}\set{ \Phi(U) := \phi\big( \Ac(UU^{\top}) - B \big) \ : \ U\in\R^{m\times r}}.
\end{equation}
Clearly,  \eqref{eq:bounded_assumption} is a generalization of the least-squares problem in \cite{liu2015efficient}. 
In addition, we cannot directly apply alternating scheme to solve \eqref{eq:sym_LMA_prob} without reformulating it into other form.
The optimality condition of \eqref{eq:sym_LMA_prob} is written as
\begin{equation}\label{eq:sym_opt_cond}
U^{\top}\Ac^{*}\big(\nabla{\phi}(\Ac(UU^{\top}) - B)\big) = 0.
\end{equation}
Any $U_{\star}$ satisfying this condition is called a \textit{stationary point} of \eqref{eq:sym_LMA_prob}. 
We again assume that the set of stationary points $\mathcal{U}_{\star}$ of \eqref{eq:sym_LMA_prob} is nonempty.

We now customize Algorithm \ref{alg:A1} to find a stationary point of \eqref{eq:sym_LMA_prob}.
Since $U = V$, the symmetric GN direction can be computed from Remark \ref{re:symmetric_case} as
\begin{equation*}
D_{U} = \left(\Id - 0.5P_{U}\right)Z(U^{\dagger})^{\top}, \quad \textrm{where} \quad Z =   -L_{\Phi}^{-1}\Ac^{\ast}\big(\nabla{\phi}(\Ac(UU^{\top}) - B) \big).
\end{equation*}
Combining this step and modifying the linesearch procedure   \eqref{eq:ls_condition}, we can describe a new variant of Algorithm \ref{alg:A1} for solving \eqref{eq:sym_LMA_prob} as in  Algorithm \ref{alg:A1b}.

\begin{algorithm}[!ht]\caption{(\textit{Symmetric linesearch Gauss-Newton algorithm} (SLs-GN))}\label{alg:A1b}
\begin{algorithmic}[1]
   \STATE {\bfseries Initialization:} Given a tolerance $\varepsilon > 0$. 
    Choose  $U_0\in\R^{m\times r}$.
    Set $\alpha_0  := 1$ and  $c_1 := 0.5$.
    \FOR{$k = 0$ {\bfseries to} $k_{\max}$}
	\STATE\label{step:symgn_gn}  \textit{Gauss-Newton direction}: 
	Evaluate $Z_k := -L_{\Phi}^{-1}\Ac^{\ast}\nabla{\phi}(\Ac(U_kU_k^{\top}) - B)$ and compute
	\begin{equation*}
	D_{U_k} := \left(\Id_m - 0.5P_{U_k}\right)Z_k(U_k^{\dagger})^{\top}.
	\vspace{-2ex}
	\end{equation*}
	\STATE  If $\Vert D_{U_k} \Vert_F \leq \varepsilon\max\set{1, \Vert U_k\Vert_F }$, then TERMINATE.
	\STATE\label{step:symgn_ls} \textit{Linesearch}: Find the smallest number $i_k \geq 0$ such that $\alpha_{i_k} := \beta^{i_k}\alpha_0$ and
	\begin{equation*} 
		\Phi(U_k+ \alpha_{i_k}D_{U_k}) \leq \Phi(U_k) - 0.5c_1\alpha_{i_k}\Vert\nabla{\Phi}(U_k)\Vert_F^2.
	\vspace{-2ex}
	\end{equation*}
	\STATE Update $U_{k+1} := U_k + \alpha_kD_{U_k}$.
   \ENDFOR
\end{algorithmic}
\end{algorithm}

\vspace{0.5ex}
\paragraph{\textbf{Per-iteration complexity:}}
Computing $U^{\dagger}$ requires one QR-factorization of an $m\times r$ matrix to get $[Q, R] = \texttt{qr}(U)$. 
Then, we   form $U^{\dagger} = R^{\dagger}Q^{T}$, where $R^{\dagger}$ is obtained by solving an upper triangle linear system. $P_{U_k}$ is computed by $P_{U_k} = U_kU^{\dagger}_k$.
Computing $Z_k$ at Step \ref{step:symgn_gn} requires $U_kU_k^{\top}$, one linear operator $\Ac$ and one adjoint $\Ac^{*}$. 
The linesearch routine at Step \ref{step:symgn_ls} requires $i_k$ function evaluations as indicated in \eqref{eq:i_k}. 
Each linesearch step needs one $U_kU_k^{\top}$ and one $\Ac(\cdot)$.

The following corollary summarizes the convergence properties of Algorithm \ref{alg:A1b}, which is a direct consequence of Lemma~\ref{le:descent_dir} and Theorem~\ref{th:convergence_guarantee1}.

\begin{corollary}\label{co:sym_congergence_guarantee1}
Let $\set{U_k}$ be  generated by Algorithm \ref{alg:A1b}. Then, under  Assumption A.\ref{as:A1}:
\begin{itemize}
\item[$\mathrm{(a)}$] There exists $\underline{\alpha}_k := \min\set{1, \frac{L_{\Phi}\sigma_{\min}^3(U_k)}{2\Vert \nabla{\Phi}(U_k)\Vert_F}, \frac{3\sigma_{\min}(U_k)^4}{32\sigma_{\max}(U_k)^2\Phi'(U_k)}} \in (0, 1]$ such that 
\begin{equation}\label{eq:sym_descent}
\Phi(U_{k} +  \alpha_k D_{U_k}) \leq \Phi(U_k) - \frac{\alpha_k}{128L_{\Phi}}\frac{\sigma_{\min}^2(U_k)}{\sigma_{\max}^4(U_k)}\Vert \nabla{\Phi}(U_k)\Vert_F^2, \quad\forall \alpha_k \in (0, \underline{\alpha}_k].
\end{equation}
Consequently, the linesearch procedure at Step 5 is well-defined $($i.e., it terminates after a finite number of iterations $i_k$$)$. 

\item[$\mathrm{(b)}$] If there exists $\underline{\sigma} > 0$ such that $\sigma_{\min}(U_k) \geq \underline{\sigma}$ for all $k\geq 0$ and $\set{U_k}$ is bounded, then  $\lim_{k\to\infty}\Vert \nabla{\Phi}(U_k)\Vert_F = 0$, and any limit point of $\set{U_k}$ is in $\mathcal{U}_{\star}$.
\end{itemize}
\end{corollary}

The results in Corollary~\ref{co:sym_congergence_guarantee1} is fundamentally different from \cite{liu2015efficient}, even when $\phi(\cdot) := (1/2)\Vert\cdot\Vert_2^2$ and $\Ac$ is identical, since  $B$ is not positive definite.
We note that Algorithm~\ref{alg:A2} can be specified to handle the symmetric case \eqref{eq:sym_LMA_prob} by substituting Steps \ref{step:gnadmm_gn} and \ref{step:gnadmm_ls} by Steps \ref{step:symgn_gn} and \ref{step:symgn_ls} in Algorithm \ref{alg:A1b}, respectively. 
We omit the details of this specification.

\section{Numerical experiments}\label{sec:num_experiments}
In this section, we first discuss some implementation remarks.
Next, we compare the full-step GN scheme and AMA. 
Then, we test Algorithm~\ref{alg:A1} on a low-rank matrix approximation problem and compare it with standard SVDs.
Finally, we apply Algorithms \ref{alg:A1}, \ref{alg:A2} and \ref{alg:A1b} to solve three problems: matrix completion,  matrix recovery, and robust low-rank matrix recovery.

\subsection{Implementation remarks}\label{sec:impl_aspects}
The following aspects are implemented in our experiments.

\vspace{0.5ex}
\paragraph{\textbf{Computing initial points:}}
Since \eqref{eq:LMA_prob} is nonconvex, the performance of the above algorithms strongly depends on an initial point.
Principally, these algorithms still converge from any initial point. 
However, we propose to use the following simple procedure for finding an initial point:
We first form a matrix $M\in\R^{m\times n}$ such that $\Ac(M) = B$. 
Then, we compute the $r$-truncated  SVD of $M$ as $[U_f, \Sigma_f, V_f]$ and form
\begin{equation*}
U_0 := U_f(:,1:r)\Sigma_f(1:r)^{1/2}~~\text{and}~~ V_0 := V_f(:,1:r)\Sigma_f(1:r)^{1/2}. 
\end{equation*}
In Algorithm \ref{alg:A2}, given $[U_0, V_0]$, we set $W_0 := \Ac(U_0V_0^{\top})-B$ and $\Lambda_0 := \boldsymbol{0}^l$.

\vspace{0.5ex}
\paragraph{\textbf{Stopping criterions:}}
We can implement  different stopping criterions for Algorithms \ref{alg:A1} and \ref{alg:A2}. 
The first criterion is based on the optimality condition \eqref{eq:opt_cond}:
\begin{equation}\label{eq:stop_cond1}
\max\big\{ \Vert U^{\top}_k\Phi^{\prime}(U_kV_k^{\top})\Vert_F, \Vert\Phi^{\prime}(U_kV_k^{\top})V_k\Vert_F \big\} \leq \varepsilon_1\max\set{1, \norm{B}_F},
\end{equation}
where $\Phi^{\prime}(UV^{\top}) := \Ac^{\ast}\left(\nabla{\phi}(\Ac(UV^{\top}) - B)\right)$.
We can terminate Algorithm \ref{alg:A1} if
\begin{equation}\label{eq:stop_cond2}
\max\set{ \Vert D_{U_k}\Vert_F, \Vert D_{V_k}\Vert_F} \leq \varepsilon_1\max\set{1, \norm{B}_F}.
\end{equation}
We can add to Algorithm \ref{alg:A2} the following condition for feasibility in \eqref{eq:constr_LMA_prob2}:
\begin{equation}\label{eq:stop_cond3}
\Vert U_kV_k^{\top}- W_k\Vert_F \leq \varepsilon_1\max\set{1, \norm{B}_F}.
\end{equation}
When  $\phi(\cdot) := (1/2)\Vert\cdot\Vert_F^2$ and the optimal value is zero, we also use
\begin{equation}\label{eq:stop_cond4}
\Vert \Ac(U_kV_k^{\top}) - B\Vert_F \leq \varepsilon_2\max\set{1, \norm{B}_F}.
\end{equation}
Similar stopping criterions are applied to Algorithm \ref{alg:A1b}.

\vspace{0.5ex}
\paragraph{\textbf{Penalty parameter update:}} 
Theoretically, we can fix any parameter $\rho$ as indicated in \eqref{eq:choice_of_rho}.
However, in Section \ref{sec:num_experiments}, we follow the update rule used in  \cite{Shen2012} but with different parameters.
We also use the full-step GN scheme at Step~\ref{step:gnadmm_gn}. 
%
\subsection{Comparison of Gauss-Newton and Alternating Minimization Algorithm}\label{subsec:experiment0}
In order to observe the advantage of the GN scheme over AMA (also called alternating direction method)  for solving \eqref{eq:LMA_prob}, we compare these algorithms on the following special case of \eqref{eq:LMA_prob}:
\begin{equation}\label{eq:qp_LMA_prob}
\Phi^{\star} := \min_{U\in\R^{m\times r}, V\in\R^{n\times r}}\set{ \Phi(U,V) := (1/2)\Vert \Ac(UV^{\top}) - B\Vert_2^2 }.
\end{equation}
Since $\Ac$ is nonidentical, we upper bound $(1/2)\Vert\Ac(\cdot) - B\Vert_2^2$ as
\begin{equation*}
\begin{array}{lcl}
\frac{1}{2}\Vert\Ac(UV^{\top}) - B\Vert_2^2 &\leq & \frac{1}{2}\Vert\Ac(U_kV^{\top}_k) - B\Vert_2^2 + \frac{1}{2}\Vert UV^{\top} - (U_kV^{\top}_k - L^{-1}\Ac^{\ast}(\Ac(U_kV_k^{\top}) - B))) \Vert_2^2\vspace{0.75ex}\\
&& - {~} \frac{1}{2L}\Vert \Ac^{\ast}(\Ac(U_kV_k^{\top}) - B)) \Vert_2^2,
\end{array}
\end{equation*}
where $L := \Vert\Ac\Vert^2$ is the Lipschitz constant of the gradient of $(1/2)\Vert\Ac(\cdot) - B\Vert^2_2$.

Let $Z_k := L^{-1}\Ac^{\ast}(\Ac(U_kV_k^{\top}) - B))$. We can write  AMA as
\begin{equation}\label{eq:adm_scheme}
{}\left\{ \begin{array}{lcl}
U_{k+1} &:= & \displaystyle\argmin_{U}\set{ (1/2)\Vert UV_k^{\top} - (U_kV_k^{\top} - Z_k)\Vert_2^2},\vspace{1ex}\\
V_{k+1} &:= & \displaystyle\argmin_{V}\set{  (1/2)\Vert U_{k+1}V^{\top} - (U_kV_k^{\top} - Z_k)\Vert_2^2}.
\end{array}\right.{}\tag{AMA}
\end{equation} 
We compare this algorithm and the following full-step GN scheme of \eqref{eq:gau_newton_scheme}:
\begin{equation}\label{eq:fsGN_scheme}
(U_{k+1}, V_{k+1}) := \argmin_{U,V}\set{ (1/2)\Vert U_kV^{\top} + UV_k^{\top} - (U_kV_k^{\top} + Z_k)\Vert_2^2}. \tag{FsGN}
\end{equation}
Clearly, \ref{eq:adm_scheme} alternates between $U$ and $V$ and solves for them separately, while \ref{eq:fsGN_scheme} linearizes $UV^{\top}$ and solves   for $U_{k+1}$ and $V_{k+1}$ simultaneously.

We implement these  schemes in Matlab and running on a MacBook laptop with a 2.6 GHz Intel Core i7 processor and 16GB memory.
The input data is generated as follows.
For $\Ac$, we generate an $(mn\times mn)$-matrix from either a fast Fourier transform (fft) or a standard Gaussian distribution, and take $l$ random sub-samples from the rows of this matrix to form $\Ac$, where $l \leq mn$.  
We generate $B = \Ac( U^{\natural}(V^{\natural})^{\top}) + \Nc(0,\sigma^2\Id)$, where $U^{\natural}\in\R^{m\times r}$ and $V^{\natural}\in\R^{n\times r}$ are given matrices, and $\Nc(0,\sigma^2\Id)$ is i.i.d. Gaussian noise of variance $\sigma^2$. 
We consider two cases: the underdetermined case with $l<r(m+n)$, and the overdetermined case with $l > r(m+n)$.
In the first case, problem \eqref{eq:qp_LMA_prob} always has a solution with zero residual.
We choose  $(U_0, V_0)$ randomly, which may not  be in the local convergence region of the GN method.

\begin{figure}[!ht]
\begin{center}
\vspace{-0ex}
\includegraphics[width = 1.0\textwidth]{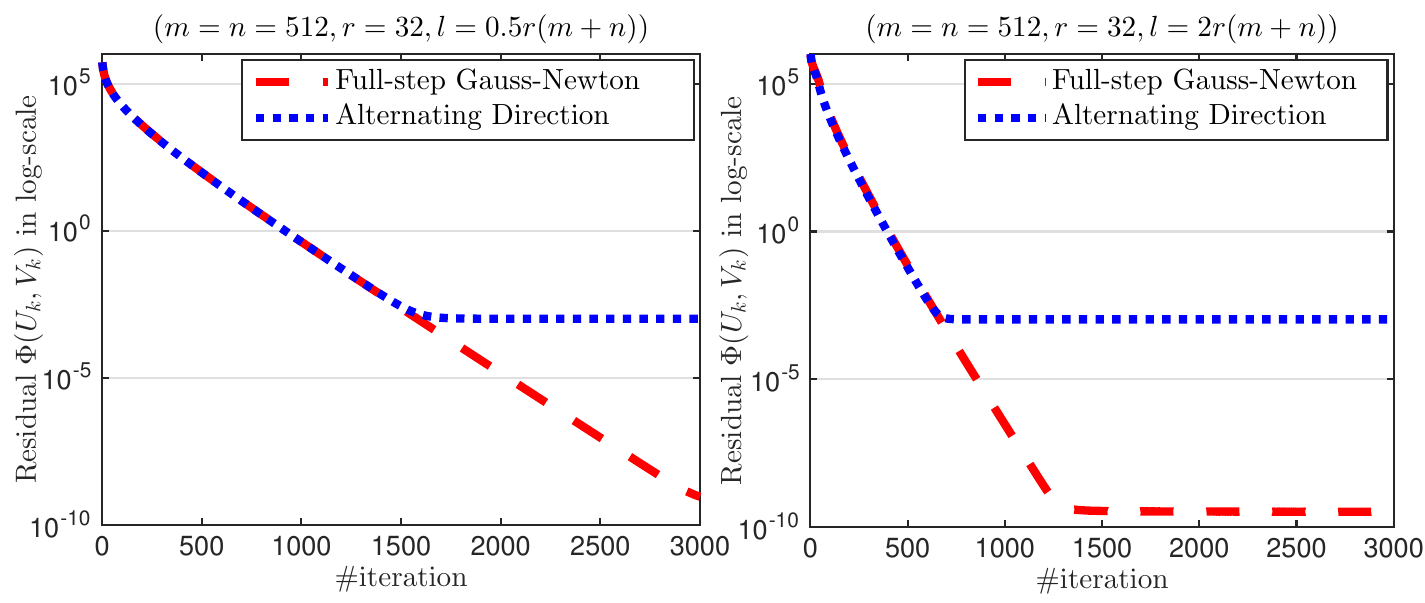}
\vspace{-3ex}
\caption{A comparison between \ref{eq:fsGN_scheme} (\textit{Legend}: Full-step Gauss-Newton)  and \ref{eq:adm_scheme} (\textit{Legend:} Alternating Direction). 
Left: The underdetermined case -- $l = 0.5r(m+n)$. Right: The overdetermined case -- $l=2r(m+n)$.}
\label{fig:compare_adm_gn}
\vspace{-1ex}
\end{center}
\end{figure}

Figure \ref{fig:compare_adm_gn} shows the convergence behavior  of the two algorithms.
The right plot  is  $l = 2r(m+n)$, and  the left one is   $l = 0.5r(m+n)$, where $m = n =512$ and $r = 32$.

\begin{figure}[ht!]
\begin{center}
\includegraphics[width = 1.0\textwidth]{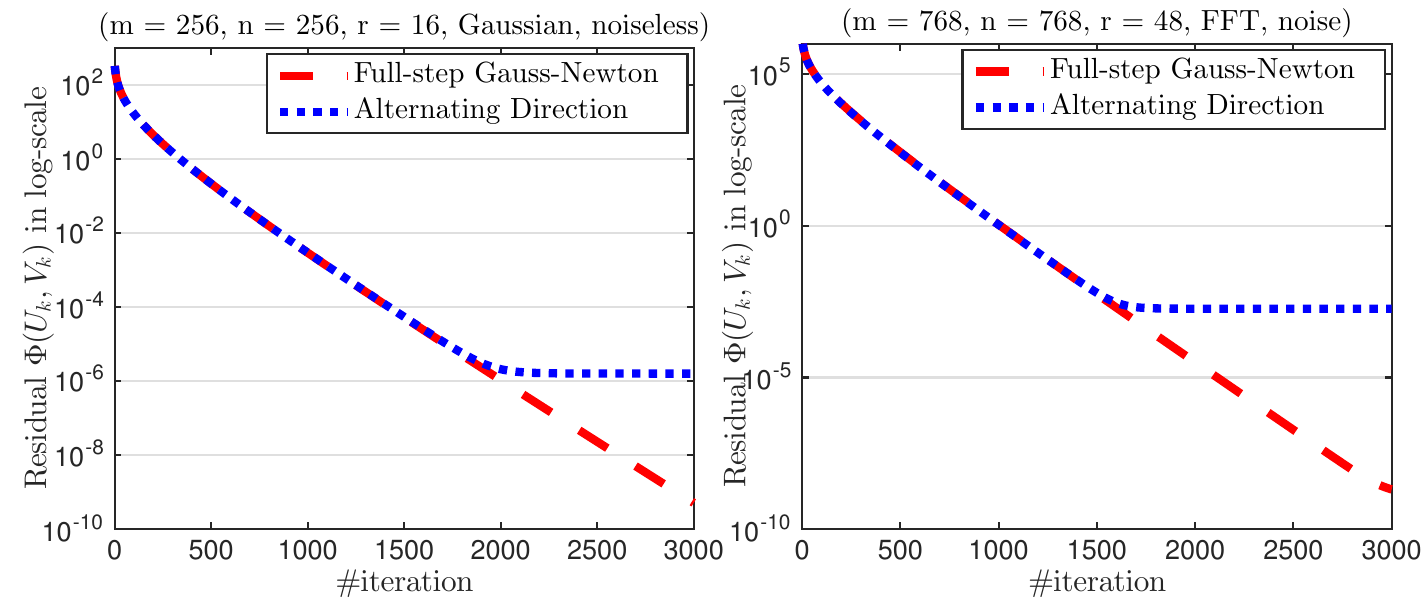}
\vspace{-3ex}
\caption{A comparison between \ref{eq:fsGN_scheme} and \ref{eq:adm_scheme}. 
\textit{Left:} sparse  Gaussian operator without noise. 
\textit{Right:} subsampling FFT linear operator with noise.}
\label{fig:compare_adm_gn2}
\vspace{-1ex}
\end{center}
\end{figure}

We can see from Figure \ref{fig:compare_adm_gn} that  both algorithms perform very similarly in early iterations, but then \ref{eq:fsGN_scheme} gives better result in terms of accuracy (terminated around $10^{-9}$ in the overdetermined case due to the nonzero objective residual), while \ref{eq:adm_scheme} is saturated at a certain level, and does not improve the objective values.
In addition, Figures \ref{fig:compare_adm_gn} and \ref{fig:compare_adm_gn2} show that the full-step Gauss-Newton scheme has a local linear convergence rate for the underdetermined case.
However, as a compensation, \ref{eq:fsGN_scheme} requires one $(r\times m)$-matrix multiplication $U^{\top}U$ and one $(r\times r)$-inverse compared to \ref{eq:adm_scheme}. 
This suggests that we can perform \ref{eq:adm_scheme} in early iterations and switch to \ref{eq:fsGN_scheme} if \ref{eq:adm_scheme} does not make significant progress to improve the objective values.

We test the underdetermined case by choosing a Gaussian operator $\Ac$ generated as $\Ac = \frac{1}{\sqrt{l}}\texttt{sprandn}(l,mn,0.05)$. 
The convergence of two algorithms on this dataset is plotted in Figure \ref{fig:compare_adm_gn2} (left).
Finally, we consider the effect of noise to both algorithms by adding a Gaussian noise with $\sigma^2 = 10^{-3}$.
The performance of these algorithms is plotted in Figure \ref{fig:compare_adm_gn2} (right).

We can observe from Figure \ref{fig:compare_adm_gn2} the same behavior as in the previous test. 
Our \ref{eq:fsGN_scheme} still maintains a local linear convergence even with noise, while \ref{eq:adm_scheme} is saturated at a certain level of the objective values.

\subsection{Low-rank matrix factorization and linear subspace selection}
We consider a special case of \eqref{eq:LMA_prob} by taking  $\phi(\cdot) := (1/2)\Vert\cdot \Vert_F^2$ and $\Ac = \Id$ as
\begin{equation}\label{eq:LMA_prob_exam1}
\Phi^{\star} := \min_{U, V}\set{\Phi(U,V) := (1/2)\Vert UV^{\top} - B\Vert_F^2 \ : \ U\in\R^{m\times r}, V\in\R^{n\times r} }.
\end{equation}
Although this problem has a closed form solution by truncated SVD, our objective is to compare  the full-step GN variant of Algorithm \ref{alg:A1} with standard Matlab singular value decomposition routines: \texttt{svds} and \texttt{lansvd}. 
The full-step GN scheme for \eqref{eq:LMA_prob_exam1} is presented as
\begin{equation}\label{eq:fs_GN}
V_{k+1}^{\top} := U_k^{\dagger}B  \quad \text{and} \quad U_{k+1} := U_k + (B - U_k^{\top}V_{k+1})(V_k^{\dagger})^{\top}. 
\end{equation}
At each iteration, \eqref{eq:fs_GN} requires two $(r\times r)$-matrix inverses $U_k^{\top}U_k$ and $V_k^{\top}V_k$, and three $(m\times r)$- or $(n\times r)$- matrix - $(r\times r)$-matrix multiplications.
We  compute these two inverses by Cholesky decomposition. 
We note that we do not form  the $(m\times n)$-matrix $U_kV_k^{\top}$ at each iteration, but we can occasionally compute it to check the objective value if required.
We choose  $U_0 := [\Id_r, \boldsymbol{0}_{(m-r)\times r}^{\top}]^{\top}$ and $V_0 := [\boldsymbol{0}_{(n-r)\times r}^{\top}, \Id_r]^{\top}$ as a starting point, where $\Id_r$ is the identity matrix.

Scheme \eqref{eq:fs_GN} generates two low-rank matrices $U_k$ and $V_k$ so that $U_kV_k^{\top} \approx B$.
We can perform a Rayleigh--Ritz (RR) routine to orthonormalize $U_k$ and $V_k$, 
\begin{itemize}
\item Compute $[Q_u, R_u] = \texttt{qr}(U_k, 0)$ and $[Q_v, R_v] = \texttt{qr}(V_k, 0)$, the two economic QR-factorizations of size $r$.
\item Compute $[U_r, \Sigma_r, V_r] = \texttt{svd}(Q_u^{\top}BQ_v)$ the singular value decomposition of the $r\times r$ matrix $Q_u^{\top}BQ_v$.
\item Then form $U = Q_uU_r$ and $V = Q_vV_r$ to obtain two orthogonal matrices $U$ and $V$ of the size $m\times r$ and $n\times r$ so that $[U, \Sigma, V] = \texttt{svds}(B, r)$.
\end{itemize}
Here, \eqref{eq:fs_GN} works on a symmetric positive definite matrix compared to \cite{liu2015efficient}.

Now, we test \eqref{eq:fs_GN} in combining with the Rayleigh--Ritz procedure, and compare it with \texttt{svds} and \texttt{lansvd}. 
We generate an input matrix $B$ of size $m\times n$ with rank $r$.
Once $m$ is chosen, we set $n = m$ and either $r = 0.01 \times m$ or $r = 0.05 \times m$ (which is either $1\%$ or $5\%$ of the problem size, respectively).
Then, we generate $B\in\R^{m\times n}$ using the following Matlab code:
\begin{small}
\begin{framed}
\begin{verbatim}
min_mn      = min(m, n);
nnz_sig_vec = [1:1:r].^(-0.01);
sig_vec     = [nnz_sig_vec(:); zeros(min_mn-r, 1)];
n_sig_vec   = sqrt(length(sig_vec))/norm(sig_vec(:), 2)*sig_vec;
B           = gallery('randcolu', n_sig_vec, max(m, n), 1);    
G           = sprandn(m, n, nnz(B)/(m*n));
M_mat       = B + 0.1*norm(B, 'fro')*G/norm(G, 'fro');
\end{verbatim}
\end{framed}
\end{small}
Clearly, the singular values $\sigma_i$ of $B$ are clustered into two parts: $\sigma_i = i^{-0.01}$ for $i=1,\cdots, r$, and $\sigma_i=0$ for $i=r+1,\cdots,\min\set{m,n}$.
In addition, an i.i.d. Gaussian noise $\frac{\norm{B}_F}{10\norm{G}_F}G$ is added to $B$, where $G = \mathcal{N}(0, \sigma\Id)$, with $\sigma$ being the sparsity of $B$.
We terminate \eqref{eq:fs_GN} using either  \eqref{eq:stop_cond1} or  \eqref{eq:stop_cond2} with $\varepsilon_1 =  10^{-6}$ or $\varepsilon_2 = 10^{-4}$, respectively. 
We also terminate \texttt{svds} and \texttt{lansvd} using $\texttt{tol} = 10^{-4}$, which is a moderate accuracy.

\begin{figure}[!ht]
\begin{center}
\includegraphics[width = 1.0\textwidth]{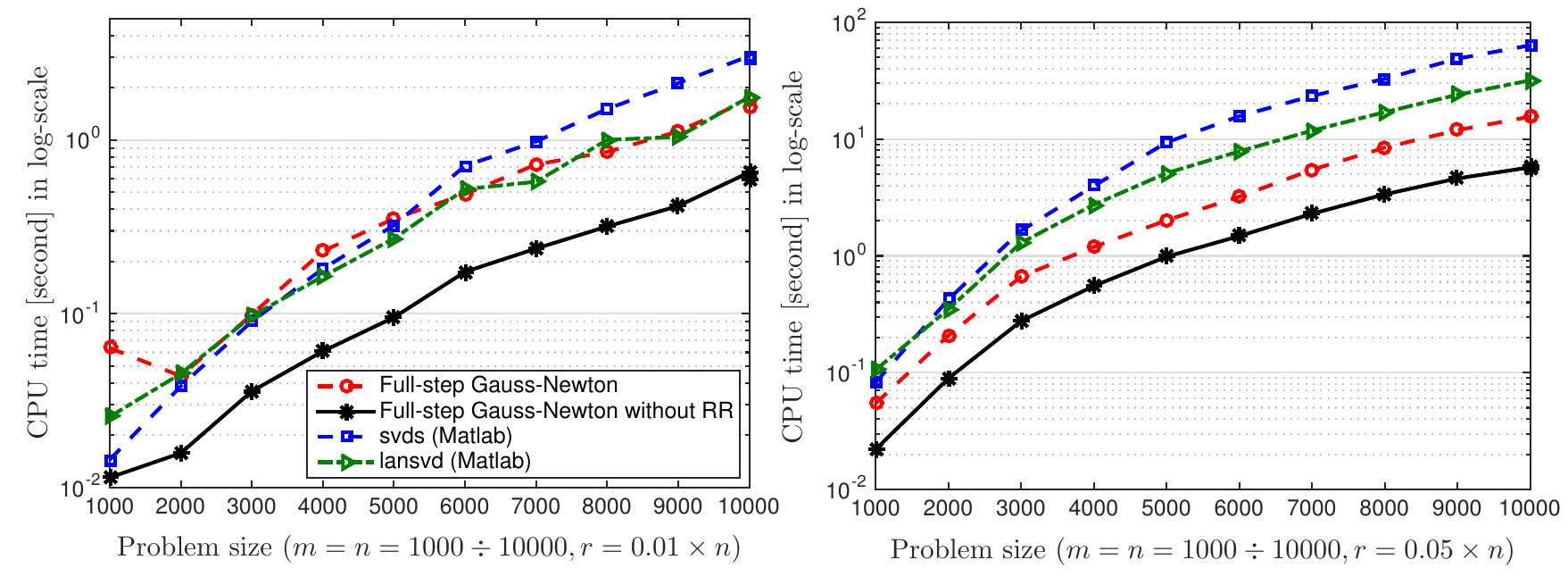}
\vspace{-3ex}
\caption{A comparison between the full-step GN scheme, $\mathtt{Matlab~SVDS}$, and $\mathtt{Matlab~LanSVD}$ on $10$ problem sizes $($from $1,000$ to $10,000$$)$, and two different ranks. The result is on the average of $10$ random runs for each problem size.}
\label{fig:LMA_time_vs_size}
\vspace{-1ex}
\end{center}
\end{figure}

The performance of three algorithms in terms of computational time vs. problem size is plotted in Figure \ref{fig:LMA_time_vs_size} for $10$ problems from $m=n=1,000$ to $m=n=10,000$, carried out on a MacBook laptop with a 2.6 GHz Intel Core i7 processor and 16GB memory.
We run each problem size $10$ times and compute the averaging computational time. The abbreviation \texttt{Full-step Gauss-Newton} indicates the time of both scheme \eqref{eq:fs_GN} and  Rayleigh-Ritz procedure, while \texttt{Full-step Gauss-Newton without RR} only counts for the time of \eqref{eq:fs_GN}.
Figure \ref{fig:LMA_time_vs_size} (left) shows the performance with $r = 0.01 \times m$, while Figure~\ref{fig:LMA_time_vs_size} (right) reveals the case $r = 0.05\times m$.

When the rank $r$ is about $1\%$ of the problem size, \eqref{eq:fs_GN} is comparable to \texttt{lansvd} while it is slightly better than \texttt{svds}.
However, when the rank $r$ is increased to $5\%$ of the problem size, \eqref{eq:fs_GN} clearly outperforms  both \texttt{lansvd} and \texttt{svds}.

\subsection{Recovery with Pauli measurements in quantum tomography}\label{subsec:quantum}
We consider a $d$ spin-1/2 system with unknown state $S$ as described in \cite{gross2010quantum}. 
A $d$-qubit  Pauli matrix is given by the form $w = \otimes_{i=1}^dw_i$, where $w_i \in \set{\boldsymbol{1}, \sigma^x, \sigma^y, \sigma^z}$ is a given set of elements. 
There are $n^2$, $n = 2^d$, such matrices denoted by $w(s)$ with $s\in\set{1,\cdots, n^2}$. 
A compressive sensing  procedure takes $m$ integer numbers $s_1,\cdots, s_m\in\set{1,\cdots, n^2}$ randomly and measures the expected values $\trace{S w(s_i)}$. 
Then, it solves the following convex problem to construct the unknown states:
\begin{equation}\label{eq:QT}
\trace{X} = 1, \quad \trace{Xw(s_i)} = \trace{w_i S} \quad (i=1,\cdots, m).
\end{equation}
From \cite{gross2010quantum}, the number of measurement $m$ to reconstruct the quantum states can be estimated as $m = c n r \log^2 n \ll n^2$ for some constant $c$ and the rank $r$.

Given that $X$ characterizes a density matrix, which is positive semidefinite Hermitian, we instead consider the following least-squares formulation of \eqref{eq:QT}:
\begin{equation}\label{eq:sdp}
\min_{X\in\mathcal{H}^{n}_{+}}\set{ (1/2)\Vert B - \mathcal{A}(X)\Vert_F^2  \ : \ \mathrm{trace}(X) = 1},
\end{equation}
where $\mathcal{H}^n_{+}$ is the set of positive semidefinite Hermitian matrices of size $n$, and $\Ac$ and $B$ are the measurement operator and observed measurements obtained from \eqref{eq:QT}, respectively.
Assume that $X = UU^{\top}$, where $U\in\mathbb{C}^{n\times 1}$, we can write \eqref{eq:sdp} into
\begin{equation}\label{eq:sdp_rank}
\min_{U\in \mathbb{C}^{n\times 1}}\set{ (1/2)\Vert B - \mathcal{A}(UU^{\top})\Vert_F^2 },
\end{equation}
where $\mathbb{C}^{n\times 1}$ is the set of $(n\times 1)$ - complex matrices.
Clearly, problem \eqref{eq:sdp_rank} falls into the special form \eqref{eq:sym_LMA_prob} of \eqref{eq:LMA_prob} which can be solved by Algorithm \ref{alg:A1b}.

\begin{table}[ht!]
\newcommand{\cell}[1]{#1}
\newcommand{\cella}[1]{{}#1{}}
\newcommand{\cellbf}[1]{{}{\color{blue}#1}{}}
\begin{center}
\begin{scriptsize}
\caption{Numerical results of three algorithms on noiseless and noisy data}\label{tbl:QT_exam}
\begin{tabular}{ c c c  | rrr | rrr | rrr}
\toprule
\multicolumn{3}{c}{} & \multicolumn{3}{c}{ Algorithm \ref{alg:A1b} } & \multicolumn{3}{c}{ Frank-Wolfe without LS }  & \multicolumn{3}{c}{ Frank-Wolfe with LS }\\  \midrule
\cell{\#qubits} & $m$ & $n$  & \cell{iter} & \cell{time[s]} & \cell{$\frac{\norm{B - \Ac(X)}_F}{\norm{B}_F}$} & \cell{iter} & \cell{time[s]} & \cell{$\frac{\norm{B - \Ac(X)}_F}{\norm{B}_F}$} & \cell{iter} & \cell{time[s]} & \cell{$\frac{\norm{B - \Ac(X)}_F}{\norm{B}_F}$} \\ \midrule
\multicolumn{12}{c}{The noiseless case} \\ \midrule
\cell{10} & \cell{14196} & \cell{1024} & \cell{26} & \cellbf{12.25} & \cellbf{3.21e-06} & \cell{1707} & \cell{664.90} & \cell{1.65e-03} & \cell{322} & \cell{129.35} & \cell{1.62e-03} \\ 
\cell{11} & \cell{ 31231} & \cell{2048} & \cell{25} & \cellbf{71.60} & \cellbf{2.64e-06} & \cell{1654} & \cell{2803.18} & \cell{1.61e-03} &\cell{370} & \cell{593.56} & \cell{1.54e-03} \\
\cell{12} & \cell{68140} & \cell{4096} & \cell{25} & \cellbf{696.27} & \cellbf{1.78e-06} & \cell{1577} & \cell{17990.98} & \cell{1.56e-03} &\cell{254} & \cell{1741.19} & \cell{1.54e-03} \\
\cell{13} & \cell{147635} & \cell{8192} & \cell{27} & \cellbf{1516.97} & \cellbf{1.73e-06} & \cell{648} & \cell{20574.13} & \cell{3.68e-03} &\cell{303} & \cell{9654.69} & \cell{1.52e-03} \\  \midrule
\multicolumn{12}{c}{The depolarizing noisy case (1\%) } \\ \midrule
\cell{10} & \cell{14196} & \cell{1024} &  \cell{24} & \cellbf{16.07} & \cellbf{8.99e-06} & \cell{1711} & \cell{692.16} & \cell{1.66e-03} &\cell{238} & \cell{78.22} & \cell{1.66e-03} \\
\cell{11} & \cell{ 31231} & \cell{2048} & \cell{23} & \cellbf{94.98} & \cellbf{8.80e-06} & \cell{1663} & \cell{2683.90} & \cell{1.62e-03} &\cell{258} & \cell{423.27} & \cell{1.61e-03} \\
\cell{12} & \cell{68140} & \cell{4096} & \cell{23} & \cellbf{589.73} & \cellbf{6.03e-06} & \cell{1585} & \cell{12146.01} & \cell{1.56e-03} &\cell{247} & \cell{1892.66} & \cell{1.56e-03} \\
\cell{13} & \cell{147635} & \cell{8192} & \cell{24} & \cellbf{3684.57} & \cellbf{8.76e-06} & \cell{648} & \cell{20537.15} & \cell{3.70e-03} &\cell{292} & \cell{8691.90} & \cell{1.53e-03} \\
\bottomrule
\end{tabular}
\end{scriptsize}
\end{center}
\end{table}
We test Algorithm \ref{alg:A1b} and compared it with Frank-Wolfe's method proposed in \cite{Jaggi2013}. 
We use both the standard Frank-Wolfe and its linesearch variant.  
We generate $U_0 := [\Id_r, \boldsymbol{0}_{(n-r)\times r}^{\top}]^{\top}$ and terminate Algorithm \ref{alg:A1b} using either \eqref{eq:stop_cond1}, \eqref{eq:stop_cond2}, or \eqref{eq:stop_cond4} with $\varepsilon_1  = 10^{-9}$ and $\varepsilon_2 = 10^{-6}$, respectively.
We generate $\Ac$ and $B$ using the procedures in \cite{gross2010quantum}.
We perform two cases: noise and noiseless. 
In the noisy case, we set $S$ to be $0.99 S + 0.01\Id_n/n$ before computing the observed measurement $B$.
Since Frank-Wolfe's algorithms take long time to reach a high accuracy, we terminate them if $\norm{\Ac(X) - B}_F\leq 10^{-3}\sqrt{2}\norm{B}_F$ which is different from Algorithm \ref{alg:A1b}.

\begin{figure}[!ht]
\begin{center}
\includegraphics[width = 1.0\textwidth]{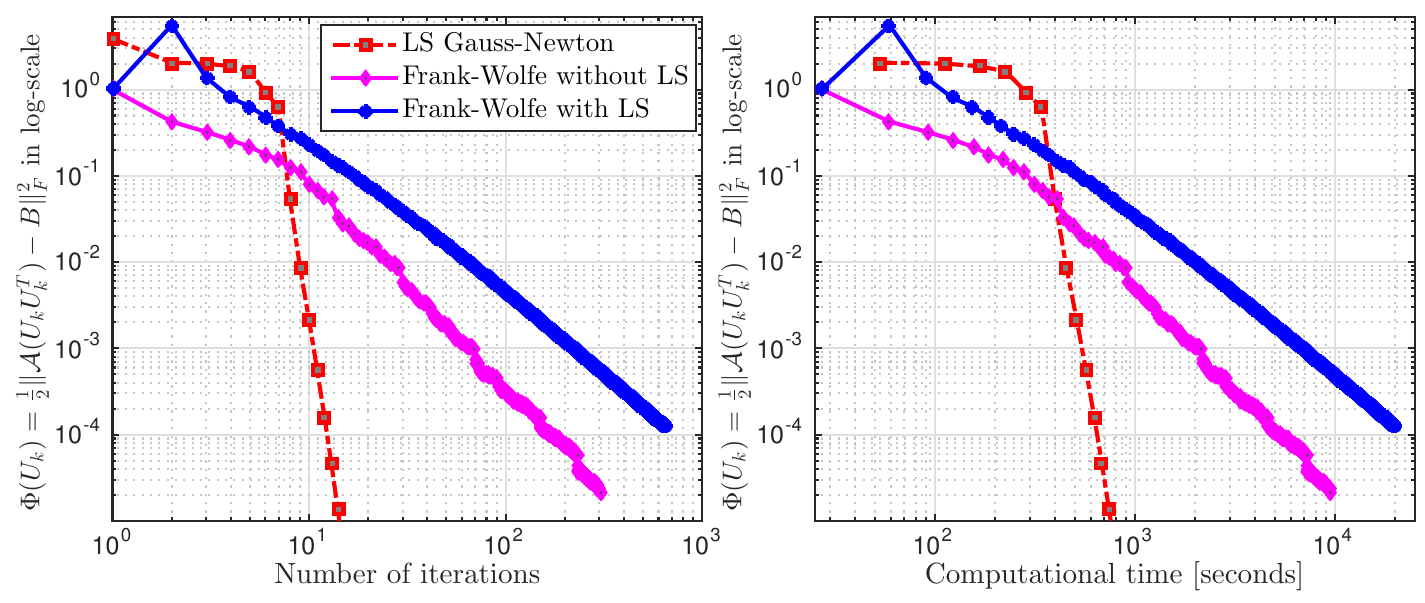}
\includegraphics[width = 1.0\textwidth]{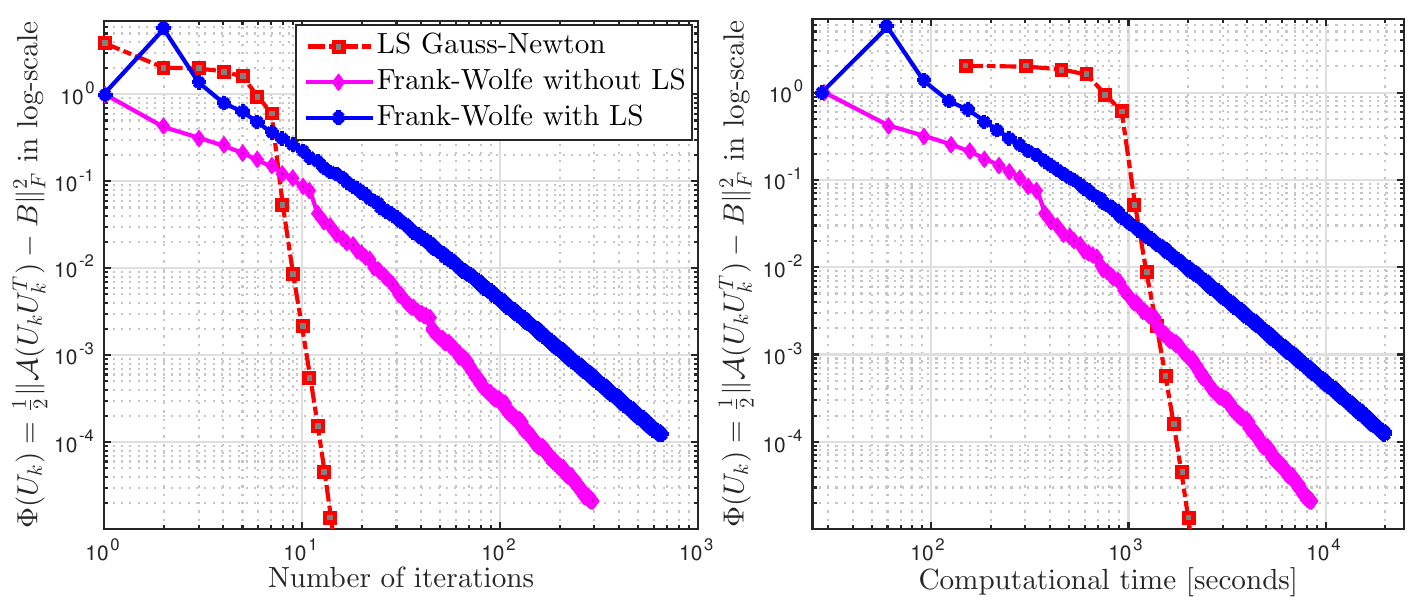}
\vspace{-3ex}
\caption{A comparison of three algorithms for the noiseless case (first row) and for the $0.01$ depolarizing noisy case (second row).}
\label{fig:QT_convg_noiseless}
\vspace{-1ex}
\end{center}
\end{figure}

We test on 4 problems of the size $d$ with $d \in \set{10, 11,12, 13}$ being the number of qubits running one a single node of an Intel(R) Xeon(R) 2.67GHz cluster with 4GB memory, but can share up to 320GB RAM.
The results and performance of three algorithms are reported in Table \ref{tbl:QT_exam}, where $m$ is the number of measurements, $n = 2^d$, \textrm{iter} is the number of iterations, \textrm{time[s]} is the computational time in seconds. 
The convergence behavior of three algorithms for both noiseless and noisy cases with $d = 13$ is also plotted in Figure \ref{fig:QT_convg_noiseless}.

We can observe from our results that Algorithm \ref{alg:A1b} highly outperforms the two Frank-Wolfe variants. 
It also reaches a highly  accurate solution after a few iterations. 
However, each iteration of Algorithm \ref{alg:A1b} is more expensive than that of Frank-Wolfe's algorithms.
As can be seen from Figure \ref{fig:QT_convg_noiseless}, Algorithm \ref{alg:A1b} behaves like super-linearly convergent.

\subsection{Matrix completion}
Our next experiment is solve the well-known matrix completion (MC) widely used in recommender systems \cite{Candes2012b,goldfarb2011convergence,wen2012solving}.
This problem is a special case of \eqref{eq:LMA_prob} and can be written as follows:
\vspace{-0.75ex}
\begin{equation}\label{eq:mc_prob}
\min_{U, V}\set{ \Vert (1/2)\mathcal{P}_{\Omega}(UV^{\top}) - B \Vert_F^2 \ : \ U \in \R^{r\times m}, V\in\R^{r\times n} },
\vspace{-0.75ex}
\end{equation} 
where $\mathcal{P}_{\Omega}$ is a selection operator on an index subset $\Omega$, and $B$ is the set of observed entries.

There are two major approaches to solve \eqref{eq:mc_prob}. The first one is using a convex relaxation for the rank constraint via nuclear or max norms. 
Methods based on this approach have been widely developed, including SVT \cite{cai1956singular}, and [accelerated] gradient descent \cite{goldfarb2011convergence,toh2010accelerated}.
The second approach is using nonconvex optimization, including, e.g., OpenSpace \cite{keshavan2009gradient} and LMaFit \cite{Lin2009,wen2012solving}.

In this experiment, we select the most efficient algorithms for our comparison: the over-relaxation alternating direction method (LMaFit)  in \cite{wen2012solving}, and the accelerated proximal gradient method (APGL) in \cite{toh2010accelerated}. 
We will test the four algorithms on synthetic datasets and the three first algorithms on some real datasets.

\vspace{0.5ex}
\paragraph{\textbf{Synthetic datasets:}}
Since data in rating systems is often integer, our synthetic dataset is generated as follows. 
We randomly generate two integer matrices $U$ and $V$ whose entries are in $\set{1,\cdots, 5}$ of the size $m\times r$ and $n\times r$, respectively. 
Then, we form $M = UV^{\top}$. Finally, we randomly take either $50\%$ or $30\%$ entries of $M$ as an output matrix $B$. We can also add a standard Gaussian noise to $B$ if necessary.
A Matlab script for generating such a dataset is given below.
\begin{small}
\begin{framed}
\begin{verbatim}
U_org    = randi(5, m, r);
V_org    = randi(5, n, r);
M_org    = U_org*V_org';
s        = round(0.5*m*n); 
Omega    = randsample(m*n, s); 
M_omega  = M_org(Omega);
B        = M_omega + sigma*randn(size(M_omega));
\end{verbatim}
\end{framed}
\end{small}
We first test these algorithms with a fixed rank $r$ and $50\%$ randomly observed entries, which is relative dense. 
We terminate Algorithms \ref{alg:A1} and \ref{alg:A2} using the conditions given in Section \ref{sec:impl_aspects} with $\varepsilon_1 = 10^{-6}$ and $\varepsilon_2 = 10^{-4}$, respectively.
We also terminate LMaFit and APGL with the same tolerance $\mathtt{tol} = 10^{-4}$. 
The initial point is computed by a truncated SVD as   in Section \ref{sec:impl_aspects}. 

\begin{table}[ht!]
\newcommand{\cell}[1]{{}#1{}}
\newcommand{\cella}[1]{{}#1{}}
\newcommand{\cellb}[1]{{}{\color{blue}#1}{}}
\newcommand{\cellm}[1]{{}{\color{magenta}#1}{}}
\begin{center}
\begin{scriptsize}
\caption{Comparison of four algorithms on synthetic integer data without noise}\label{tbl:mc_syn_data1}
\begin{tabular}{ rrr | rrrrr | rrrrr}
\toprule
\multicolumn{3}{c}{} & \multicolumn{5}{c}{ Algorithm \ref{alg:A1} } & \multicolumn{5}{c}{ Algorithm \ref{alg:A2} }  \\  \midrule
$m$ & $n$ & $r$  & \cell{iter} & \cell{time[s]} & \cell{$\delta{f}_k$} & \cell{NMAE } & \cell{rank} & \cell{iter} & \cell{time[s]} & \cell{$\delta{f}_k$} & \cell{NMAE} & \cell{rank} \\ 
\cell{1000} & \cell{2000} & \cell{10} &\cell{15.9} & \cell{2.11} & \cellb{4.15e-05} & \cellb{1.39e-05} & \cell{10} & \cell{30.0} & \cell{2.07} & \cell{8.34e-05} & \cell{3.55e-05} & \cell{10} \\ 
\cell{1000} & \cell{2000} & \cell{50} &\cell{20.8} & \cellb{4.34} & \cellb{6.91e-05} & \cellb{5.78e-05} & \cell{50} & \cell{37.6} & \cell{4.43} & \cell{9.61e-05} & \cell{7.83e-05} & \cell{50} \\ 
\midrule
\cell{2500} & \cell{2500} & \cell{25} &\cell{14.3} & \cell{8.05} & \cellb{5.31e-05} & \cellb{2.97e-05} & \cell{25} & \cell{31.1} & \cell{9.18} & \cell{1.04e-04} & \cell{6.18e-05} & \cell{25} \\ 
\cell{2500} & \cell{2500} & \cell{125} &\cell{26.3} & \cell{28.94} & \cellb{7.35e-05} & \cellb{8.18e-05} & \cell{125} & \cell{35.2} & \cellb{23.68} & \cell{1.05e-04} & \cell{1.32e-04} & \cell{125} \\ 
\midrule
\cell{5000} & \cell{5000} & \cell{50} &\cell{15.7} & \cell{53.40} & \cellb{5.42e-05} & \cellb{3.90e-05} & \cell{50} & \cell{32.0} & \cell{56.71} & \cell{9.87e-05} & \cell{7.80e-05} & \cell{50} \\ 
\cell{5000} & \cell{5000} & \cell{250} &\cell{23.6} & \cell{180.70} & \cell{7.94e-05} & \cell{1.35e-04} & \cell{250} & \cell{35.0} & \cellb{165.89} & \cell{1.10e-04} & \cell{1.94e-04} & \cell{250} \\ 
\midrule
\cell{5000} & \cell{7500} & \cell{50} &\cell{14.9} & \cell{76.93} & \cell{5.10e-05} & \cell{3.64e-05} & \cell{50} & \cell{32.0} & \cell{85.72} & \cell{8.51e-05} & \cell{6.65e-05} & \cell{50} \\ 
\cell{5000} & \cell{7500} & \cell{250} &\cell{23.7} & \cell{273.30} & \cell{7.61e-05} & \cell{1.22e-04} & \cell{250} & \cell{35.0} & \cellb{245.93} & \cell{9.97e-05} & \cell{1.68e-04} & \cell{250} \\ 
\midrule
\cell{10000} & \cell{10000} & \cell{100} &\cell{16.2} & \cell{289.14} & \cell{5.99e-05} & \cell{5.86e-05} & \cell{100} & \cell{32.2} & \cell{319.92} & \cell{1.10e-04} & \cell{1.16e-04} & \cell{100} \\ 
\cell{10000} & \cell{10000} & \cell{500} &\cell{24.8} & \cell{1303.01} & \cell{8.02e-05} & \cell{1.75e-04} & \cell{500} & \cell{35.0} & \cellb{1173.38} & \cell{1.14e-04} & \cell{2.60e-04} & \cell{500} \\ 
\midrule
\multicolumn{3}{c}{} & \multicolumn{5}{c}{ LMaFit \cite{wen2012solving}} & \multicolumn{5}{c}{ APGL \cite{toh2010accelerated}}  \\  \midrule
$m$ & $n$ & $r$  & \cell{iter} & \cell{time[s]} & \cell{$\delta{f}_k$} & \cell{NMAE } & \cell{rank} & \cell{iter} & \cell{time[s]} & \cell{$\delta{f}_k$} & \cell{NMAE} & \cell{rank} \\ 
\midrule
\cell{1000} & \cell{2000} & \cell{10} &\cell{13.3} & \cellb{0.94} & \cell{4.74e-05} & \cell{1.43e-05} & \cell{10} & \cell{28.0} & \cell{4.29} & \cell{3.31e-04} & \cell{1.40e-04} & \cell{10} \\ 
\cell{1000} & \cell{2000} & \cell{50} &\cell{109.8} & \cell{12.13} & \cell{2.71e-04} & \cell{1.51e-04} & \cell{50} & \cell{28.6} & \cell{6.79} & \cell{1.06e-02} & \cell{9.02e-03} & \cell{41.6} \\ 
\midrule
\cell{2500} & \cell{2500} & \cell{25} &\cell{10.0} & \cellb{3.10} & \cell{6.98e-05} & \cell{3.86e-05} & \cell{25} & \cell{29.0} & \cell{16.64} & \cell{5.06e-04} & \cell{3.07e-04} & \cell{25} \\ 
\cell{2500} & \cell{2500} & \cell{125} &\cell{135.3} & \cell{85.78} & \cell{2.99e-04} & \cell{2.59e-04} & \cell{125} & \cell{31.6} & \cellb{20.00} & \cell{1.92e-02} & \cell{2.43e-02} & \cell{5} \\ 
\midrule
\cell{5000} & \cell{5000} & \cell{50} &\cell{10.0} & \cellb{18.43} & \cell{5.57e-05} & \cell{4.30e-05} & \cell{50} & \cell{30.6} & \cell{81.69} & \cell{8.48e-04} & \cell{6.73e-04} & \cell{50.2} \\ 
\cell{5000} & \cell{5000} & \cell{250} &\cell{140.9} & \cell{631.58} & \cellb{6.70e-05} & \cellb{8.58e-05} & \cell{250} & \cell{30.4} & \cellb{60.75} & \cell{1.38e-02} & \cell{2.44e-02} & \cell{5} \\ 
\midrule
\cell{5000} & \cell{7500} & \cell{50} &\cell{10.1} & \cellb{28.21} & \cellb{4.38e-05} & \cellb{2.92e-05} & \cell{50} & \cell{30.8} & \cell{122.47} & \cell{7.35e-04} & \cell{5.76e-04} & \cell{50} \\ 
\cell{5000} & \cell{7500} & \cell{250} &\cell{126.8} & \cell{845.90} & \cellb{7.02e-05} & \cellb{7.30e-05} & \cell{250} & \cell{30.5} & \cellb{90.74} & \cell{1.38e-02} & \cell{2.33e-02} & \cell{5} \\ 
\midrule
\cell{10000} & \cell{10000} & \cell{100} &\cell{11.0} & \cellb{112.82} & \cellb{3.10e-05} & \cellb{3.24e-05} & \cell{100} & \cell{32.7} & \cell{266.16} & \cell{2.15e-02} & \cell{2.28e-02} & \cell{5} \\ 
\cell{10000} & \cell{10000} & \cell{500} &\cell{120.4} & \cell{3818.05} & \cellb{6.25e-05} & \cellb{8.63e-05} & \cell{500} & \cell{30.3} & \cellb{206.86} & \cell{9.85e-03} & \cell{2.26e-02} & \cell{5} \\ 
\bottomrule
\end{tabular}
\end{scriptsize}
\end{center}
\end{table}

The test is conducted on $10$ problems of different sizes running on a single node of an Intel(R) Xeon(R) 2.67GHz cluster with 4GB memory, but can share up to 100GB RAM.
We run each problem size $10$ times and compute the average result and performance. 
The problem sizes and results are reported in Table \ref{tbl:mc_syn_data1} for two different ranks. 
The rank $r$ is chosen as $r = 0.01\times m$,  and $r = 0.05\times m$, which correspond to $1\%$, and $5\%$ of the problem size.
Here, \textrm{iter} and \textrm{time[s]} are the number of iterations and the computational time in seconds, respectively; \textrm{rank} is the rank of $U_kV_k^{\top}$ given by the algorithms; and
\begin{equation*}
\delta{f}_k := \Vert\Pc_{\Omega}(U_kV_k^{\top}) - B\Vert_F/\norm{B}_F \quad \text{and} \quad \mathrm{NMAE} := C^{-1}\sum_{(i,j)\in\Omega}\big\vert  (U_kV_k^{\top})_{ij} - B_{ij} \big\vert,
\end{equation*}
are the relative objective residual; and the Normalized Mean Absolute Error, respectively, where $C :=  (\max_{i,j}  B_{ij} -\min_{ij} B_{ij})\vert\Omega\vert$.

The results in Table \ref{tbl:mc_syn_data1} show that both Algorithms \ref{alg:A1} and \ref{alg:A2} produce similar results as LMaFit in terms of the relative objective residual and NMAE.
When the rank is small (i.e., $1\%$ of problem size), Algorithm \ref{alg:A1} and LMaFit have similar number of iterations, but LMaFit has better computational time.
When the rank is increasing up to $5\%$ of the problem size, both Algorithm~\ref{alg:A1} and Algorithm~\ref{alg:A2} require a fewer iterations than LMaFit, and outperform this solver in terms of computational time.
In this experiment, the number of iterations in Algorithm~\ref{alg:A2} is very similar in all the test cases, from $30$ to $38$ iterations, and similar to APGL.
Note that we fix the rank in the first three algorithms, since APGL uses a convex approach, it cannot predict well an approximate rank if it is $5\%$ of the problem size, or when the problem size is increasing.

\begin{figure}[ht!]
\begin{center}
\includegraphics[width = 1.0\textwidth]{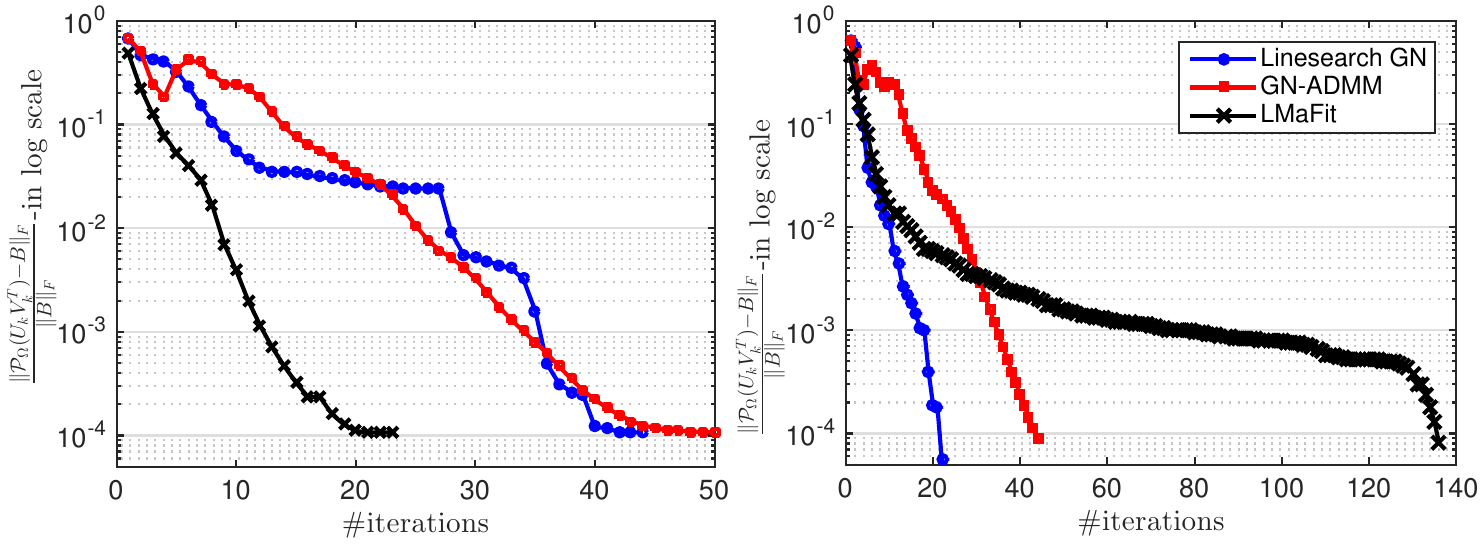} 
\vspace{-3ex}
\caption{The convergence behavior of three algorithms ($m = n = 5000$) with noise ($\sigma = 0.01$) and $30\%$ known entries (Left: $r = 0.01m$, Right: $r = 0.025m$).}
\label{fig:mc_convergence}
\vspace{-0ex}
\end{center}
\end{figure}

Now, we add i.i.d. Gaussian noise $\Nc(0, \sigma \Id)$ with $\sigma = 0.01$ to $B$ as $B := B^{\natural} + 5\times\Nc(0, \sigma \Id)$, and only randomly take $30\%$ observed entries.
The convergence behavior of three algorithms for one problem instance with $m=n=5000$ is plotted in Figure \ref{fig:mc_convergence}.
When the rank $r = 0.01m$ (i.e., $1\%$ of the problem size), LMaFit outperforms Algorithms \ref{alg:A1} and \ref{alg:A2} in terms of iterations, but when the rank $r = 0.025m$ (i.e., $2.5\%$ of the problem size), Algorithms \ref{alg:A1} and \ref{alg:A2} are much better than LMaFit. Algorithm~\ref{alg:A1} works really well in the second case, and takes only $22$ iterations.
We also observe the monotone decrease in Algorithm~\ref{alg:A1} as guaranteed by our theory, but not in Algorithm~\ref{alg:A2}.

\begin{table}[ht!]
\newcommand{\cell}[1]{{}#1{}}
\newcommand{\cella}[1]{{}#1{}}
\newcommand{\cellb}[1]{{}{\color{blue}#1}{}}
\newcommand{\cellm}[1]{{}{\color{magenta}#1}{}}
\begin{center}
\begin{scriptsize}
\caption{{}Comparison of the four algorithms on synthetic integer datasets  with noise.{}}\label{tbl:mc_syn_data2}{}
\begin{tabular}{ lrr| rrrrr |rrrrr}
\toprule
\multicolumn{3}{c}{} & \multicolumn{5}{c}{ Algorithm \ref{alg:A1} } & \multicolumn{5}{c}{ Algorithm \ref{alg:A2} }  \\  \midrule
$m$ & $n$ & $r$  & \cell{iter} & \cell{time[s]} & \cell{$\delta{f}_k$} & \cell{NMAE } & \cell{rank} & \cell{iter} & \cell{time[s]} & \cell{$\delta{f}_k$} & \cell{NMAE} & \cell{rank} \\ 
\midrule
\cell{1000} & \cell{2000} & \cell{10} &\cell{40.8} & \cell{3.56} & \cell{5.33e-04} & \cell{2.27e-04} & \cell{10} & \cell{40.5} & \cell{2.24} & \cell{5.32e-04} & \cell{2.26e-04} & \cell{10} \\ 
\cell{1000} & \cell{2000} & \cell{25} &\cell{49.1} & \cell{5.50} & \cell{6.05e-04} & \cell{2.76e-04} & \cell{25} & \cell{62.0} & \cellb{4.62} & \cell{2.09e-04} & \cell{1.31e-04} & \cell{25} \\ 
\midrule
\cell{5000} & \cell{5000} & \cell{50} &\cell{19.5} & \cell{39.82} & \cell{1.10e-04} & \cell{8.71e-05} & \cell{50} & \cell{45.0} & \cell{67.70} & \cell{1.09e-04} & \cell{8.63e-05} & \cell{50} \\ 
\cell{5000} & \cell{5000} & \cell{125} &\cell{16.6} & \cellb{56.70} & \cell{7.90e-05} & \cell{9.58e-05} & \cell{125} & \cell{40.0} & \cell{103.45} & \cell{9.50e-05} & \cell{1.19e-04} & \cell{125} \\ 
\midrule
\multicolumn{3}{c}{} & \multicolumn{5}{c}{ LMaFit \cite{wen2012solving}} & \multicolumn{5}{c}{ APGL \cite{toh2010accelerated}}  \\  \midrule
$m$ & $m$ & $r$  & \cell{iter} & \cell{time[s]} & \cell{$\delta{f}_k$} & \cell{NMAE } & \cell{rank} & \cell{iter} & \cell{time[s]} & \cell{$\delta{f}_k$} & \cell{NMAE} & \cell{rank} \\ 
\midrule
\cell{1000} & \cell{2000} & \cell{10} &\cell{31.6} & \cellb{1.71} & \cell{5.31e-04} & \cell{2.26e-04} & \cell{10} & \cell{28.0} & \cell{3.08} & \cell{6.34e-04} & \cell{2.70e-04} & \cell{10} \\ 
\cell{1000} & \cell{2000} & \cell{25} &\cell{121.0} & \cell{8.55} & \cell{2.08e-04} & \cell{1.31e-04} & \cell{25} & \cell{30.7} & \cell{5.86} & \cell{9.87e-04} & \cell{5.93e-04} & \cell{25} \\ 
\midrule
\cell{5000} & \cell{5000} & \cell{50} &\cell{20.0} & \cellb{30.39} & \cell{1.07e-04} & \cell{8.49e-05} & \cell{50} & \cell{28.5} & \cell{52.96} & \cell{4.97e-03} & \cell{3.82e-03} & \cell{48.2} \\ 
\cell{5000} & \cell{5000} & \cell{125} &\cell{48.0} & \cell{121.11} & \cell{7.19e-05} & \cell{7.34e-05} & \cell{125} & \cell{31.3} & \cellb{46.99} & \cell{1.92e-02} & \cell{2.42e-02} & \cell{5} \\ 
\bottomrule
\end{tabular}
\end{scriptsize}
\end{center}
\end{table}

Finally, we test three first algorithms on two problem instances with $30\%$ observed entries in $B$ and with i.i.d. Gaussian noise $\Nc(0, 0.01\Id)$). The results of this test is reported in Table \ref{tbl:mc_syn_data2}. 
LMaFit remains working well for then low-rank cases, while getting slower when the rank $r$ increases.
Algorithms \ref{alg:A1} and \ref{alg:A2} have similar performance in this case. 

\vspace{0.5ex}
\paragraph{\textbf{Real datasets:}} 
Now, we test three algorithms: Algorithms~\ref{alg:A1} and~\ref{alg:A2}, and LMaFit on  MovieLens and Jester jokes datasets available on \url{http://grouplens.org/datasets/movielens/}. 
For the MovieLens dataset, we test our algorithms on $5$ problems: ``movie-lens-latest (small)'', ``movie-lens'' 100k, 1M, 10M, and 20M, which we abbreviate by ``movie(s)", and ``moviexM" in Table~\ref{tbl:mc_real_data2}, respectively.
We also test all problems in Jester joke dataset: ``jester-1'', ``jester-2'', ``jester-3'', and ``jester-all''.

\begin{table}[ht!]
\newcommand{\cell}[1]{{}#1{}}
\newcommand{\cells}[1]{{}#1{}}
\newcommand{\cella}[1]{{}#1{}}
\newcommand{\cellb}[1]{{}{\color{blue}#1}{}}
\begin{center}
\caption{Summary of results of four the algorithms for MC on ``real'' datasets}\label{tbl:mc_real_data2}
\begin{scriptsize}
\begin{tabular}{ lrr | rrr | rrr | rrr}
\toprule
\multicolumn{3}{c}{} & \multicolumn{3}{c}{ Algorithm \ref{alg:A1} } & \multicolumn{3}{c}{ Algorithm \ref{alg:A2} }  &  \multicolumn{3}{c}{ LMaFit \cite{wen2012solving} } \\  \midrule
Name & $m$ & $n$  & \cell{iter} & \cell{time[s]} & \cell{$\delta{f}_k$} & \cell{iter} & \cell{time[s]} &  \cell{$\delta{f}_k$} & \cell{iter} & \cell{time[s]} &  \cell{$\delta{f}_k$} \\ \midrule
\cell{jester-1} & \cell{24983} & \cell{100} &\cell{  45} & \cell{11.60} &  \cell{1.75e-01} & \cell{  59} & \cell{11.35} &  \cell{1.75e-01}&  \cell{  36} & \cellb{4.78} &  \cell{1.75e-01} \\ 
\cell{jester-2} & \cell{23500} & \cell{100} &\cell{  41} & \cell{9.93} &  \cell{1.77e-01} & \cell{  57} & \cell{11.07} &  \cell{1.77e-01}&  \cell{  34} & \cellb{5.51} &  \cell{1.77e-01} \\ 
\cell{jester-3} & \cell{24938} & \cell{100} &\cell{  30} & \cell{5.15} &  \cell{9.04e-04} & \cell{  32} & \cell{4.94} &  \cell{9.66e-04}&  \cell{  25} & \cellb{2.13} &  \cell{9.26e-04} \\ 
\cell{jester-all} & \cell{73421} & \cell{100} &\cell{  48} & \cell{35.12} &  \cell{1.65e-01} & \cell{  57} & \cell{30.12} &  \cell{1.65e-01}&  \cell{  36} & \cellb{12.21} &  \cell{1.65e-01} \\ 
\cell{movie(s)} & \cell{668} & \cell{10325} &\cell{ 200} & \cell{16.69} &  \cell{1.64e-03} & \cell{  87} & \cellb{7.14} &  \cell{1.58e-03}&  \cell{ 200} & \cell{44.18} &  \cell{1.58e-03} \\ 
\cell{movie100k} & \cell{943} & \cell{1682} &\cell{ 200} & \cell{9.66} &  \cell{1.03e-02} & \cell{  84} & \cellb{4.87} &  \cell{1.00e-02}&  \cell{ 200} & \cell{15.63} &  \cell{1.00e-02} \\ 
\cell{movie1M} & \cell{6040} & \cell{3706} &\cell{  79} & \cell{41.71} &  \cell{1.18e-01} & \cell{  42} & \cellb{21.96} &  \cell{1.19e-01}&  \cell{  70} & \cell{49.84} &  \cell{1.18e-01} \\ 
\cell{movie10M} & \cell{69878} & \cell{10677} &\cell{  69} & \cell{109.02} &  \cell{2.14e-01} & \cell{  33} & \cellb{38.32} &  \cell{2.15e-01}&  \cell{  61} & \cell{40.48} &  \cell{2.14e-01} \\ 
\cell{movie20M} & \cell{138493} & \cell{26744} &\cell{  89} & \cell{307.22} &  \cell{2.30e-01} & \cell{  37} & \cellb{117.23} &  \cell{2.30e-01}&  \cell{  87} & \cell{133.86} &  \cell{2.30e-01} \\ 
\midrule
\multicolumn{3}{c}{} & \multicolumn{3}{c}{ Algorithm \ref{alg:A1} } & \multicolumn{3}{c}{ Algorithm \ref{alg:A2} }  &  \multicolumn{3}{c}{ LMaFit \cite{wen2012solving} } \\  \midrule
Name & $m$ & $n$  & \cell{rank} & \cell{$\delta{x}_k$} & \cell{NMAE} &  \cell{rank} & \cell{$\delta{x}_k$} & \cell{NMAE} &  \cell{rank} & \cell{$\delta{x}_k$} & \cell{NMAE} \\ \midrule
\cell{jester-1} & \cell{24983} & \cell{100} &\cell{80} & \cell{4.71e-01} & \cell{2.29e-02} & \cell{80} & \cell{4.71e-01} & \cell{2.36e-02} & \cell{80} & \cell{4.59e-01} & \cell{2.30e-02} \\ 
\cell{jester-2} & \cell{23500} & \cell{100} &\cell{80} & \cell{4.82e-01} & \cell{2.35e-02} & \cell{80} & \cell{4.82e-01} & \cell{2.42e-02} & \cell{80} & \cell{4.78e-01} & \cell{2.39e-02} \\ 
\cell{jester-3} & \cell{24938} & \cell{100} &\cell{80} & \cell{8.78e-04} & \cell{1.08e-05} & \cell{80} & \cell{8.78e-04} & \cell{7.22e-05} & \cell{80} & \cell{9.87e-04} & \cell{2.08e-05} \\ 
\cell{jester-all} & \cell{73421} & \cell{100} &\cell{80} & \cell{4.09e-01} & \cell{1.95e-02} & \cell{80} & \cell{4.09e+01} & \cell{2.05e-02} & \cell{80} & \cell{3.95e-01} & \cell{1.98e-02} \\ 
\cell{movie(s)} & \cell{668} & \cell{10325} &\cell{100} & \cell{1.36e-01} & \cell{3.76e-04} & \cell{100} & \cell{1.36e-01} & \cell{6.91e-04} & \cell{100} & \cell{1.36e-01} & \cell{3.97e-04} \\ 
\cell{movie100k} & \cell{943} & \cell{1682} &\cell{100} & \cell{1.10e-04} & \cell{5.24e-03} & \cell{100} & \cell{1.10e-04} & \cell{4.87e-03} & \cell{100} & \cell{1.00e-04} & \cell{4.64e-03} \\ 
\cell{movie1M} & \cell{6040} & \cell{3706} &\cell{100} & \cell{2.34e-01} & \cell{8.28e-02} & \cell{100} & \cell{2.34e-01} & \cell{8.44e-02} & \cell{100} & \cell{2.32e-01} & \cell{8.32e-02} \\ 
\cell{movie10M} & \cell{69878} & \cell{10677} &\cell{20} & \cell{5.86e-01} & \cell{1.34e-01} & \cell{20} & \cell{5.86e-01} & \cell{1.35e-01} & \cell{20} & \cell{5.81e-01} & \cell{1.34e-01} \\ 
\cell{movie20M} & \cell{138493} & \cell{26744} &\cell{10} & \cell{6.29e-01} & \cell{1.42e-01} & \cell{10} & \cell{6.29e-01} & \cell{1.42e-01} & \cell{10} & \cell{6.30e-01} & \cell{1.42e-01} \\ 
\bottomrule
\end{tabular}
\end{scriptsize}
\end{center}
\end{table}

In this test, since the data in ``movie10M'' and ``movie20M'' is sparse, we run the three algorithms on a MacBook laptop with a 2.6 GHz Intel Core i7 processor and 16GB memory.
We use \texttt{C-mex} routines in Matlab to compute $\Pc_{\Omega}(UV^{\top})$ in three algorithms to avoid forming $U^{\top}V$.
We terminates our algorithms based on the objective  obtained from LMaFit such that the three algorithms have similar objective values.

The result is summarized in Table \ref{tbl:mc_real_data2}, where we add a new measurement defined by $\delta{x}_k := \frac{1}{\vert\Omega\vert}\sum_{(i,j)\in\Omega}\left\vert \lfloor (U_kV_k^{\top})_{ij}\rfloor - B_{ij}\right\vert$ to measure the agreement ratio between the recovered matrix $M_k := U_kV_k^{\top}$ and the observed data $B$ projected onto $\Omega$.
Due to our stopping criterion, three  algorithms produce similar results in terms of the objective residuals, solution agreement, and NMAE. 
LMaFit works well on the Jester jokes dataset, but the computational time on these problems is relatively small.
Algorithm~\ref{alg:A2} works well on Movielen dataset, especially for Movie 10MB and Movie 20MB.
As mentioned previously, Algorithm~\ref{alg:A1} often achieve better solution in terms of accuracy if we run it long enough, while LMaFit and Algorithm~\ref{alg:A2} can be used to achieve a low or medium accurate solution for matrix completion.

\subsection{Robust low-rank matrix recovery}
We consider the following  nonsmooth problem in low-rank matrix recovery:
\begin{equation}\label{eq:rpca}
\min_{U,V} \set{ \Vert UV^{\top} - B\Vert_1 \ : \ U \in \R^{m\times r}, V\in\R^{n\times r}},
\end{equation}
where $\norm{Z}_1 := \sum_{ij}\abs{Z_{ij}}$ is the $\ell_1$-norm of $Z$.
This a low-rank matrix recovery problem with the $\ell_1$-norm, which can be referred to as a robust recovery as opposed to the standard square loss. 
This formulation is often used in background extraction, see, e.g., \cite{Shen2012}.

Clearly, we can solve \eqref{eq:rpca} using our ADMM-GN  scheme above, which can be written as
\begin{equation}\label{eq:gn_admm_for_rpca}
\left\{\begin{array}{lcl}
V_{k+1}^{\top} &:= &  U_k^{\dagger}(B + W_k - \Lambda_k), \vspace{1ex}\\
U_{k+1}    &:= &  U_k + \left(B + W_k - \Lambda_k - U_kV_{k+1}^{\top}\right)(V_k^{\dagger})^{\top}, \vspace{1ex}\\
W_{k+1}   &:= & \mathrm{prox}_{\rho^{-1}\norm{\cdot}_1}\left(U_{k+1}V_{k+1}^{\top} + \Lambda_k - B\right) \vspace{1ex}\\
\Lambda_{k+1} &:= & \Lambda_k + (U_{k+1}V_{k+1}^{\top} - W_{k+1}).
\end{array}\right.
\end{equation}
We apply this scheme to solve the \eqref{eq:rpca} using  video surveillance datasets at \url{http://perception.i2r.a-star.edu.sg/bk_model/bk_index.html}.
We implement \eqref{eq:gn_admm_for_rpca} in Algorithm \ref{alg:A2} and compare it with the augmented Lagrangian method proposed in \cite{Shen2012}, which we denote by L1-LMaFit. 
We use the same strategy as in L1-LMaFit to update the penalty parameter $\rho$, while using $U_0 := [\Id_r, \boldsymbol{0}_{(m-r)\times r}]$ and $V_0 := [\Id_r, \boldsymbol{0}_{(n-r)\times r}]$ as an initial point.
As suggested in \cite{Shen2012}, we choose the rank $r$ to be $r = 1$ when testing gray-scale video data.
As experienced, L1-LMaFit was based on alternating minimization idea, which can be saturated. Hence, we run both algorithm up to $100$ iterations to observe the outcome. 
The computational time and the relative objective value $\Vert U_kV_k^{\top}-B\Vert_1/\Vert B\Vert_1$ of these two algorithms are reported in Table \ref{tbl:rpca}.{}
\begin{table}[ht!]
\newcommand{\cell}[1]{{}#1{}}
\newcommand{\cells}[1]{{}#1{}}
\newcommand{\cellbf}[1]{{}{\color{blue}#1}{}}
\begin{center}
\caption{Summary of results of the two algorithms for video background extraction.}\label{tbl:rpca}
\begin{footnotesize}
\begin{tabular}{ c c c| rr| rr}
\toprule
\multicolumn{3}{c}{Video data} & \multicolumn{2}{c}{ Algorithm \ref{alg:A2} } & \multicolumn{2}{c}{ L1-LMaFit \cite{Shen2012}} \\ \midrule
\cells{Video} & \cells{Resolution} & \cells{\#Frames} & \cells{Time} & \cells{$\Vert U_kV_k^{\top}-B\Vert_1/\Vert B\Vert_1$} & \cells{Time} & \cells{$\Vert U_kV_k^{\top}-B\Vert_1/\Vert B\Vert_1$} \\ \midrule
\cells{Escalator} & \cell{$130\times 160$} & \cell{200} & \cell{12.48} & \cell{$9.434063 \times 10^{-2}$} & \cell{13.30} & \cell{$9.435117 \times 10^{-2}$} \\
\cells{Fountain} & \cell{$128\times 160$} & \cell{200} & \cell{13.27} & \cell{$4.197912 \times 10^{-2}$} & \cell{13.71} & \cell{$4.198963 \times 10^{-2}$} \\
\cells{Bootstrap} & \cell{$120\times 160$} & \cell{250} & \cell{15.76} & \cell{$13.103802 \times 10^{-2}$} & \cell{16.91} & \cell{$13.107209 \times 10^{-2}$} \\
\cells{Curtain} & \cell{$128\times 160$} & \cell{250} & \cell{18.43} & \cell{$2.965992 \times 10^{-2}$} & \cell{25.45} & \cell{$2.969248 \times 10^{-2}$} \\
\cells{Campus} & \cell{$128\times 160$} & \cell{300} & \cell{24.83} & \cell{$9.315523 \times 10^{-2}$} & \cell{30.10} & \cell{$9.316343 \times 10^{-2}$} \\
\cells{Hall} & \cell{$144\times 176$} & \cell{300} & \cell{31.63} & \cell{$5.708911 \times 10^{-2}$} & \cell{39.05} & \cell{$5.709121 \times 10^{-2}$} \\
\cells{ShoppingMall} & \cell{$256\times 320$} & \cell{350} & \cell{82.22} & \cell{$4.442732 \times 10^{-2}$} & \cell{85.48} & \cell{$4.442907 \times 10^{-2}$} \\
\cells{WaterSurface} & \cell{$128\times 160$} & \cell{350} & \cell{35.92} & \cell{$3.607625 \times 10^{-2}$} & \cell{40.25} & \cell{$3.607747 \times 10^{-2}$} \\
\bottomrule
\end{tabular}
\end{footnotesize}
\end{center}
\end{table}

We can observe from Table \ref{tbl:rpca} that the computational time in both algorithms is almost the same. 
This is consistent with our theoretical result, since the per-iteration complexity of the two algorithms is almost the same when we choose $r = 1$. However, Algorithm \ref{alg:A2} provides a slightly better objective value since it still improves the objective when running further compared  to L1-LMaFit.
Here, we use the full-step variant of Algorithm \ref{alg:A2}, a fast convergence guarantee can be achieved when a good initial point is provided. 
This remains unclear in L1-LMaFit \cite{Shen2012}. 
Unfortunately, global convergence of our variant as well as L1-LMaFit has not been known yet.

\section{Conclusions}\label{sec:conclude}
We have proposed a new Gauss-Newton scheme to approximate a stationary point of a class of low-rank matrix nonconvex optimization problems.
Our method features several advantages from classical Gauss-Newton (GN) method such as fast local convergence, achieving high accuracy solutions compared to the well-known alternating minimization algorithm (AMA).
We have proposed a linesearch GN algorithm and established its global and local convergence under standard assumptions. 
We have also specified this algorithm to the symmetric case, where AMA is not applicable.
Then, we have combined our GN scheme with the alternating direction method of multipliers (ADMM) to design a new ADMM-GN algorithm that has global convergence guarantee and low per-iteration complexity. 
Several numerical experiments have been presented to demonstrate the theory and show the advantages of nonconvex optimization approaches.
The theory and algorithms presented in this paper can be extended to different directions, including constrained low-rank matrix/tensor optimization.

\vskip 2mm
\noindent{\bf Acknowledgments.}
The author was partly supported by the Office of Naval Research under Grant No.
 ONR-N00014-20-1-2088.
 The author would like to thank Dr. Zheqi Zhang for proving some Matlab codes to conduct the experiments in Subsection~\ref{subsec:experiment0}.

\appendix
\section{The proof of technical results}
We provide the full proofs of all the technical results  in the main text. 

\subsection{The proof of Lemma \ref{le:gauss_newton_dir}: Closed form of Gauss-Newton direction}\label{apdx:le:gauss_newton_dir}
Le us define $x := [\vec{D_V^{\top}}, \vec{D_U}]$ and $b := [\mathrm{vec}(U^{\top}B), \mathrm{vec}(BV)]$.
Then, we can write \eqref{eq:opt_cond_subprob2} as
$\mathcal{B}x = b$, where $\mathcal{B} := \begin{bmatrix} \Id_n \otimes U^{\top}U & V\otimes U^{\top}\\ V^{\top}\otimes U & V^{\top}V\otimes \Id_m \end{bmatrix}$.
We can show that
\begin{equation*}
\mathcal{B} = \begin{bmatrix} (\Id_n\otimes U^{\top})(\Id_n\otimes U) & (\Id_n\otimes U^{\top})(V\otimes \Id_m) \\ (V^{\top}\otimes \Id_m)(\Id_n\otimes U) & (V^{\top}\otimes \Id_m)(V\otimes \Id_m) \end{bmatrix} = \begin{bmatrix}\Id_n\otimes U^{\top}\\  V^{\top}\otimes \Id_m\end{bmatrix}\begin{bmatrix}\Id_n\otimes U  & V\otimes \Id_m\end{bmatrix}.
\end{equation*}
By \cite[Fact. 7.4.24]{Bernstein2005}, we have $\rank{[\Id_n\otimes U,  V\otimes \Id_m]} \leq (m+n-r)r$. Hence,  $\rank{\mathcal{B}}$ in  \eqref{eq:opt_cond_subprob2} does not exceed $r(m+n-r) < r(m+n)$.

Next, we can rewrite $b = [(\Id_n\otimes U^{\top})\mathrm{vec}(B); (V^{\top}\otimes \Id_m)\mathrm{vec}(B)]$.
If we consider the extended matrix $\bar{\mathcal{B}} := [\mathcal{B}, b]$, then we can express it as
\begin{equation*}
\bar{\mathcal{B}} = \begin{bmatrix}\Id_n\otimes U^{\top}\\  V^{\top}\otimes \Id_m\end{bmatrix}\begin{bmatrix}\Id_n\otimes U  & V\otimes \Id_m & \mathrm{vec}(B)\end{bmatrix}.
\end{equation*}
This shows that $\rank{\bar{\mathcal{B}}} = \rank{\mathcal{B}}$. 
Hence, by the well-known consistency Rouch\'{e}--Capelli theorem, \eqref{eq:opt_cond_subprob2} has a solution.
 
Now, we find the closed form \eqref{eq:gn_dir}.
Since $\rank{U} = \rank{V} = r$, both $U^{\top}U$ and $V^{\top}V$ are invertible. 
Pre-multiplying the first equation of \eqref{eq:opt_cond_subprob2} by $(U^{\top}U)^{-1}$ and rearranging the result, we have 
\begin{equation}\label{eq:Dv}
D_V^{\top} = (U^{\top}U)^{-1}U^{\top}(Z - D_UV^{\top}).
\end{equation}
Substituting this expression into the second equation of \eqref{eq:opt_cond_subprob2} we get 
\begin{equation}\label{eq:Dv2}
(\Id - U(U^{\top}U)^{-1}U^{\top})D_VV^{\top}V = (\Id - U(U^{\top}U)^{-1}U^{\top})ZV.
\end{equation}
Using the definition of the projections $P_U$, $P_V$, $P^{\perp}_U$ and $P^{\perp}_V$, we have from \eqref{eq:Dv2} that $P^{\perp}_UD_VV^{\top}V = P^{\perp}_UZV$. 
Post-multiplying this expression by $(V^{\top}V)^{-1}$, we obtain 
\begin{equation}\label{eq:Dv3}
P^{\perp}_UD_V = P^{\perp}_UZV(V^{\top}V)^{-1}.
\end{equation}
Assume that $D_U := D_U^0 + U\hat{D}_r$, where $D_U^0$ is a given vector in the null space of $U^{\top}$, i.e., $U^{\top}D_U^0 = 0$, and $\hat{D}_r \in\R^{r\times r}$ is an arbitrary matrix.
Substituting this expression into \eqref{eq:Dv3} and noting that $P^{\perp}_UU = 0$, we obtain
\begin{equation*}
D_U^0 = D_U^0 - U(U^{\top}U)^{-1}U^{\top}D_U^0 + P^{\perp}_UU\hat{D}_r = P^{\perp}_UZV(V^{\top}V)^{-1}.
\end{equation*}
Hence, we finally get
\begin{equation*} 
D_U = P^{\perp}_UZV(V^{\top}V)^{-1} + U\hat{D}_r, ~~\text{for any}~\hat{D}_r\in\R^{r\times r},
\end{equation*}
which is exactly the first term in \eqref{eq:DuDv}.
Substituting this $D_U$ into \eqref{eq:Dv} to yield  the second term of \eqref{eq:DuDv} as
\begin{equation*}
\begin{array}{lcl}
D^{\top}_V &= & (U^{\top}U)^{-1}U^{\top}\left(Z - P^{\perp}_UZV(V^{\top}V)^{-1}V^{\top} - U\hat{D}_rV^{\top}\right) \vspace{1ex}\\
&= & (U^{\top}U)^{-1}U^{\top}Z - U^{\top}U)^{-1}U^{\top}P^{\perp}_UZV(V^{\top}V)^{-1}V^{\top} - \hat{D}_rV^{\top} \vspace{1ex}\\
&=& (U^{\top}U)^{-1}U^{\top}Z -  \hat{D}_rV^{\top}.
\end{array}
\end{equation*}
Since $\hat{D}_r$ is arbitrary in $\R^{r\times r}$, we choose $\hat{D}_r := \frac{1}{2}(U^{\top}U)^{-1}U^{\top}ZV(V^{\top}V)^{-1} \in\R^{r\times r}$. 
Substituting this choice  into \eqref{eq:DuDv}, we obtain 
\begin{equation*}
D_U = \big(\Id_m - (1/2)P_U\big)ZV(V^{\top}V)^{-1} \quad \text{and} \quad D_V^{\top} = (U^{\top}U)^{-1}U^{\top}Z\big(\Id_n - (1/2)P_V\big),
\end{equation*}
which is  \eqref{eq:gn_dir}. 
Hence, the solution set of \eqref{eq:opt_cond_subprob2} forms an $(r\times r)$-linear subspace.

Finally, let us denote the residual term in the objective of \eqref{eq:subprob2} by $R(D_U,D_V) := UD_V^{\top} + D_UV^{\top} - Z$.
Then, using the expression \eqref{eq:DuDv} we can easily show that 
\begin{equation*}
\begin{array}{lcl}
R(D_U,D_V) &= & U((U^{\top}U)^{-1}U^{\top}Z -  \hat{D}_rV^{\top}) + (P^{\perp}_UZV(V^{\top}V)^{-1} + U\hat{D}_r)V^{\top} - Z \vspace{1ex}\\
&=  & P_UZ + P^{\perp}_UZP_V - P_UZ - P^{\perp}_UZ = P^{\perp}_UZP_V - P^{\perp}_UZ.
\end{array}
\end{equation*}
Hence, we can write 
\begin{equation*} 
R(D_U,D_V) = -P^{\perp}_UZP^{\perp}_V \quad \text{and} \quad (1/2)\Vert R(D_U,D_V)\Vert_F^2 = (1/2)\Vert P^{\perp}_UZP^{\perp}_V\Vert_F^2.
\end{equation*}
The last term $(1/2)\Vert P^{\perp}_UZP^{\perp}_V\Vert_F^2$ is the optimal value of  \eqref{eq:subprob2}.
\Eproof

\subsection{The proof of Lemma \ref{le:descent_dir}: Descent property of GN algorithm}\label{apdx:le:descent_dir}
Let us define $U(\alpha) := U + \alpha D_U$ and $V(\alpha) := V + \alpha D_V$ for  $\alpha > 0$. Then
\begin{equation}\label{eq:lm23_proof_est1}
U(\alpha)V(\alpha)^{\top} = UV^{\top} + \alpha (UD_V^{\top} +  D_UV^{\top}) + \alpha^2D_UD_V^{\top}.
\end{equation}
Let $W := UD_V^{\top} +  D_UV^{\top}$ and  $r(\alpha) := \Vert U(\alpha)V(\alpha)^{\top} - UV^{\top} - Z\Vert_F^2$.
Using \eqref{eq:lm23_proof_est1} we have
\begin{equation*}
\begin{array}{lcl}
r(\alpha) &= & \Vert  \alpha (UD_V^{\top} +  D_UV^{\top}) + \alpha^2D_UD_V^{\top} - Z\Vert^2_F \vspace{1ex}\\
& = & \Vert Z \Vert_F^2 + \alpha^2\Vert W\Vert^2_F + \alpha^4\Vert D_UD_V^{\top}\Vert_F^2 + 2\alpha^3\iprods{W, D_UD_V^{\top}} - 2\alpha\iprods{W, Z} - 2\alpha^2\iprods{Z, D_UD_V^{\top}}.
\end{array}
\end{equation*}
Now, using the fact that 
\begin{align*}
\iprods{W, Z-W} &= \iprods{UD_V^{\top} +  D_UV^{\top}, P^{\perp}_UZP^{\perp}_V} = \trace{(D_VU^{\top} + VD_U^{\top})P^{\perp}_UZP^{\perp}_V} = 0, 
\end{align*}
we can further expand $r(\alpha)$ as
\begin{equation}\label{eq:lm23_proof_est2}
\begin{array}{lcl}
r(\alpha) &=&  \Vert Z \Vert_F^2 - \alpha(2-\alpha)\Vert W\Vert^2_F + \alpha^4\Vert D_UD_V^{\top}\Vert_F^2 \vspace{1ex}\\
&& - {~} 2\alpha^2(1-\alpha)\iprods{W, D_UD_V^{\top}} + 2\alpha^2\iprods{W-Z, D_UD_V^{\top}}.
\end{array}
\end{equation}
Using the pseudo-inverse of $U$ and $V$ and $(V^{\top})^{\dagger}U^{\dagger} = (UV^{\top})^{\dagger}$, we can show that
\begin{align*}
D_UD_V^{\top} 
& = \big(\Id_m -0.5P_U\big)Z(UV^{\top})^{\dagger}Z\big(\Id_n - 0.5P_V\big).
\end{align*}
From the optimality condition \eqref{eq:opt_cond} and the definition of $Z = -L_{\Phi}^{-1}\Ac^{*}\nabla{\phi}(\Ac(UV^{\top}) - B))$, we can show that $\nabla_{U}{\Phi}(U,V) = -L_{\Phi}U^{\top}Z$ and $\nabla_{V}{\Phi}(U,V) = -L_{\Phi}ZV$. 
However, since $D_U$ and $D_V$ are given by \eqref{eq:gn_dir}, we express
\begin{align*}
\left\{\begin{array}{lcl}
D_U     &= & -\frac{1}{L_{\Phi}}\big(P^{\perp}_U+\frac{1}{2}P_U\big)\nabla_{U}{\Phi}(U,V)(V^{\top}V)^{-1}, \vspace{1ex}\\
D_V^{\top} &= & -\frac{1}{L_{\Phi}}(U^{\top}U)^{-1}\nabla_{V}{\Phi}(U,V)\big(P^{\perp}_V + \frac{1}{2}P_V\big).
\end{array}\right.
\end{align*}
Using this expression, we can write $\nu := \norm{D_U}_F^2 + \norm{D_V}_F^2$ as 
\begin{equation*}
\nu = \frac{1}{L_{\Phi}^2}\Vert \big(P^{\perp}_U + \frac{1}{2}P_U\big)\nabla_{U}{\Phi}(U,V)(V^{\top}V)^{-1}\Vert_F^2 
 + \frac{1}{L_{\Phi}^2}\Vert (U^{\top}U)^{-1}\nabla_{V}{\Phi}(U,V)\big(P^{\perp}_V + \frac{1}{2}P_V\big)\Vert_F^2.\nonumber
\end{equation*}
Hence, we can estimate 
\begin{align}\label{eq:lm23_proof_est4}
\frac{\Vert \nabla_U{\Phi}(U,V)\Vert_F^2}{4L_{\Phi}^2(\sigma_{\max}(V))^4} + \frac{\Vert \nabla_V{\Phi}(U,V)\Vert_F^2}{4L_{\Phi}^2(\sigma_{\max}(U))^4} \leq \nu   \leq \frac{\Vert \nabla_U{\Phi}(U,V)\Vert_F^2}{L_{\Phi}^2(\sigma_{\min}(U))^4} + \frac{\Vert \nabla_V{\Phi}(U,V)\Vert_F^2}{L_{\Phi}^2(\sigma_{\min}(V))^4},
\end{align}
where $\sigma_{\min}(\cdot)$ and $\sigma_{\max}(\cdot)$ are the smallest and largest singular values of $(\cdot)$, respectively. 
Let $\sigma_{\max} := \max\set{\sigma_{\max}(U), \sigma_{\max}(V)}$ and $\sigma_{\min} := \min\set{\sigma_{\min}(U), \sigma_{\min}(V)}$. 
Using  $\Vert \nabla{\Phi}(U,V)\Vert_F^2 = \Vert \nabla_U{\Phi}(U,V)\Vert_F^2 + \Vert \nabla_V{\Phi}(U,V)\Vert_F^2$,  \eqref{eq:lm23_proof_est4} leads to
\begin{align}\label{eq:lm23_proof_est5}
\frac{\Vert \nabla{\Phi}(U,V)\Vert_F^2}{4L_{\Phi}^2\sigma_{\max}^4} \leq \nu  = \norm{D_U}_F^2 + \norm{D_V}_F^2 \leq \frac{\Vert \nabla{\Phi}(U,V)\Vert_F^2}{L_{\Phi}^2\sigma_{\min}^4}.
\end{align}
Next, using the orthonormality, we estimate $\norm{W}_F^2$ as follows:
\begin{equation}\label{eq:lm23_proof_est6}
\begin{array}{lcl}
\norm{W}_F^2 &= & \Vert UD_V^{\top} + D_UV^{\top}\Vert_F^2 = \Vert UD_V^{\top}\Vert_F^2 + \Vert D_UV^{\top}\Vert_F^2 \vspace{1ex}\\
&=&  \trace{D_V(U^{\top}U)D_V^{\top}} + \trace{D_U(V^{\top}V)D_U^{\top}} \vspace{1ex}\\
&\geq & (\sigma_{\min}(U))^2\Vert D_U\Vert_F^2 + (\sigma_{\min}(V))^2\Vert D_V\Vert_F^2 \vspace{1ex} \\ 
& \geq & \sigma_{\min}^2\nu.
\end{array}
\end{equation}
On the one hand, we estimate individually each term of the expression \eqref{eq:lm23_proof_est2} as follows:
\begin{align*}
\Vert D_UD_V^{\top}\Vert_F^2 &=  \trace{(D_UD_V^{\top})^{\top}(D_UD_V)}
 \leq \frac{1}{4}\left(\Vert D_U\Vert_F^2 + \Vert D_V\Vert_F^2\right)^2 = \frac{\nu^2}{4}.
\end{align*}
On the other hand, since $W - Z =  -P^{\perp}_UZP^{\perp}_V$ by Lemma \ref{le:gauss_newton_dir}, we can show that 
\begin{equation*}
\iprods{W-Z,D_UD_V^{\top}} = \iprods{P^{\perp}_UZP^{\perp}_V,D_UD_V^{\top}} = \trace{D_V^{\top}P^{\perp}_VZ^{\top}P^{\perp}_UD_U} \leq \Vert Z\Vert_F\Vert D_UD_V^{\top}\Vert_F.
\end{equation*}
In addition, $-\iprods{W, D_UD_V} \leq \Vert W \Vert_F\Vert D_UD_V^{\top}\Vert_F$.
Substituting these estimates into \eqref{eq:lm23_proof_est2} and using the fact that $2-\alpha \geq 1$ and $1-\alpha \leq 1$ we obtain
\begin{equation}\label{eq:lm23_proof_est7}
\arraycolsep=0.2em
\begin{array}{lcl}
r(\alpha) &\leq & \Vert Z \Vert_F^2 - \alpha\Vert W\Vert^2_F + \frac{\nu^2\alpha^4}{4} + 2\alpha^2\Vert W\Vert_F\Vert D_UD_V^{\top}\Vert_F + 2\alpha^2\Vert Z\Vert_F\Vert D_UD_V^{\top}\Vert_F \vspace{1ex}\\
&\leq & \Vert Z \Vert_F^2 - \frac{\alpha}{16}\Vert W\Vert^2_F - \frac{\alpha\nu\sigma_{\min}^2}{16}+ \frac{\nu^2\alpha^4}{4} - \frac{\alpha}{2}\Vert W \Vert_F^2 + 2\alpha^2\Vert W\Vert_F\Vert D_UD_V^{\top}\Vert_F \vspace{1ex}\\
&& - {~} \frac{3\alpha}{8}\Vert W\Vert_F^2 + 2\alpha^2\Vert Z\Vert_F\Vert D_UD_V^{\top}\Vert_F \vspace{1ex}\\
& = & \Vert Z \Vert_F^2 - \frac{\alpha\nu\sigma_{\min}^2}{16} - \frac{\alpha\nu}{16}\big( \sigma_{\min}^2 - 4\alpha^3\nu\big)_{[b_1]} 
- \frac{\alpha}{2}\Vert W\Vert_F \Big(\Vert W\Vert_F - 4\alpha\Vert D_UD_V^{\top}\Vert_F\Big)_{[b_2]} \vspace{1ex}\\
&&  - {~} \frac{\alpha}{8}\Big( 3\Vert W\Vert_F^2 -  16\alpha\norm{Z}_F\norm{D_UD_V^{\top}}_F\Big)_{[b_3]} \vspace{1ex}\\
&= & \Vert Z \Vert_F^2 - \frac{\alpha\nu\sigma_{\min}^2}{16} - \frac{\alpha \nu}{16}b_1 - \frac{\alpha}{2}\Vert W\Vert_Fb_2 - \frac{\alpha}{8}b_3.
\end{array}
\end{equation}
We estimate each term  in \eqref{eq:lm23_proof_est7}.
From \eqref{eq:lm23_proof_est6}, we can see that $W = 0$ implies $D_U = 0$ and $D_V = 0$, which is contradict to our assumption. Hence, $W \neq 0$.
First, we choose $\alpha \in (0, 1]$ such that 
\begin{equation}\label{eq:cond_on_alpha}
\sigma_{\min}^2 \geq \left(\frac{16\Vert Z\Vert_F}{3\Vert W\Vert_F}\right)^2\nu\alpha^2 \quad \text{and} \quad \sigma_{\min}^2 \geq 4\nu\alpha^2.
\end{equation}
Since $W\neq 0$, this condition allows us to compute $\alpha$ as
\begin{equation}\label{eq:cond_on_alpha2}
0 < \alpha \leq \frac{\sigma_{\min}}{2\sqrt{\nu}}\min\set{ 1, \frac{3\Vert W\Vert_F}{8\Vert Z\Vert_F}}.
\end{equation}
Under the second condition of \eqref{eq:cond_on_alpha} and $\alpha \in (0, 1]$, we have  $b_1 = \sigma_{\min}^2 - 4\nu\alpha^3 \geq  \sigma_{\min}^2 - 4\nu\alpha^2 \geq 0$. 
Next, since \eqref{eq:lm23_proof_est6} and the first condition in \eqref{eq:cond_on_alpha} we have $\sigma_{\min}^2\geq 4\alpha^2\nu$. 
Using \eqref{eq:lm23_proof_est6} we have $\Vert W\Vert_F^2 \geq \sigma_{\min}^2\nu \geq 4\alpha^2\nu^2 = 4\alpha^2\big(\norm{D_U}_F^2 + \norm{D_V}_F^2\big)^2 \geq 16\alpha^2\norm{D_UD_V^{\top}}_F^2$. 
Hence, $\Vert W\Vert_F \geq 4\alpha \norm{D_UD_V^{\top}}_F$.
This inequality leads to $b_2 =  \Vert W\Vert_F - 4\alpha\Vert D_UD_V^{\top}\Vert_F \geq 0$.

Now, using $\Vert W \Vert_F^2 \geq \sigma_{\min}^2\nu \geq  \left(\frac{16\Vert Z\Vert_F}{3\Vert W\Vert_F}\right)^2\alpha^2\nu^2$, we have $\Vert W\Vert_F \geq  \frac{16\Vert Z\Vert_F}{3\Vert W\Vert_F} \nu\alpha$. 
Therefore, we can estimate
\begin{equation*}
b_3 = 3\Vert W\Vert_F^2 - 16\alpha \Vert Z \Vert_F\Vert D_UD_V^{\top}\Vert_F \geq 3\alpha\Vert W \Vert_F \frac{16\Vert Z\Vert_F}{3\Vert W\Vert_F} - 16\alpha \Vert Z \Vert_F\Vert D_UD_V^{\top}\Vert_F = 0.
\end{equation*}
From \eqref{eq:lm23_proof_est5} we have $\sqrt{\nu} \leq \frac{\Vert\nabla{\Phi}(U,V)\Vert_F}{L_{\Phi}\sigma_{\min}^2}$, while from \eqref{eq:lm23_proof_est6} we have $\norm{W}_F \geq \sqrt{\nu}\sigma_{\min} \geq \frac{\sigma_{\min}\Vert\nabla{\Phi}(U,V)\Vert_F}{2L_{\Phi}\sigma_{\max}^2}$. 
Substituting these estimates into \eqref{eq:cond_on_alpha2} of $\alpha$ and using  $\norm{Z}_F =
 \frac{1}{L_{\Phi}}\norm{\Phi^{\prime}(UV^{\top})}_F$ we can lower estimate $\alpha$ as
\begin{equation}\label{eq:alpha_upper_bound}
0 < \alpha \leq \frac{\sigma_{\min}^3L_{\Phi}}{2\Vert\nabla{\Phi}(U,V)\Vert_F}\min\set{ 1, \frac{3\sigma_{\min}\Vert\nabla{\Phi}(U,V)\Vert_F}{16L_{\Phi}\Vert  \Phi^{\prime}(UV^{\top})\Vert_F\sigma_{\max}^2}}.
\end{equation}
Note that $\alpha \in (0, 1]$, we obtain from \eqref{eq:alpha_upper_bound} the update rule \eqref{eq:step_size_min}.

We finally estimate \eqref{eq:descent_property}. 
Since $\alpha$ satisfies \eqref{eq:step_size_min}, it follows from \eqref{eq:lm23_proof_est7} that 
\begin{equation*}
r(\alpha) \leq \Vert Z \Vert_F^2 - \frac{\alpha\nu\sigma_{\min}^2}{16} \overset{\tiny\eqref{eq:lm23_proof_est5}}{\leq} \Vert Z \Vert_F^2 - \frac{\alpha\sigma_{\min}^2}{64L_{\Phi}^2\sigma_{\max}^4}\Vert\nabla{\Phi}(U,V)\Vert^2
\end{equation*}
Substituting this inequality into \eqref{eq:surrogate} we obtain \eqref{eq:descent_property}.
\Eproof.

\subsection{The proof of Lemma \ref{le:full_rank_uv}: Full-rankness of iterates}\label{apdx:le:full_rank_uv}
Since $\rank{U} = \rank{V} = r$ by assumption, we have $\lambda_{\min}(U^{\top}U) > 0$ and $\lambda_{\min}(V^{\top}V) > 0$.
We consider $Q := (U^{\dagger})^{\top} = U(U^{\top}U)^{-1}$ and $S := (V^{\dagger})^{\top} = V(V^{\top}V)^{-1}$.  
We always have $U_{+}^{\top}(\lambda_{\max}(QQ^{\top})\Id - QQ^{\top})U_{+} \succeq 0$. 
This implies that $$\lambda_{\min}(U_{+}^{\top}U_{+})\lambda_{\max}(QQ^{\top}) \geq \lambda_{\min}((Q^{\top}U_{+})^{\top}(Q^{\top}U_{+})).$$ 
Clearly, since $Q^{\top} = U^{\dagger}$, we have $\lambda_{\max}(QQ^{\top}) = \lambda_{\min}^{-1}(U^{\top}U)$. 
Using this relation into the last inequality, we get
\vspace{-1ex}
\begin{equation}\label{eq:lm24_proof_est1}
\frac{\lambda_{\min}(U_{+}^{\top}U_{+})}{\lambda_{\min}(U^{\top}U)} \geq \lambda_{\min}((Q^{\top}U_{+})^{\top}(Q^{\top}U_{+})).
\vspace{-1ex}
\end{equation}
Hence, it is sufficient to show that $\lambda_{\min}((Q^{\top}U_{+})^{\top}(Q^{\top}U_{+}))  > 0$.
By Lemma \ref{le:gauss_newton_dir}, we have $U_{+} = U + \alpha D_U = U + \alpha(P^{\perp}_U + 0.5P_U)ZV(V^{\top}V)^{-1}$. 
Therefore,  we can compute $Q^{\top}U_{+} = \Id_m + 0.5\alpha H$, where $H := (U^{\top}U)^{-1}U^{\top}ZV(V^{\top}V)^{-1}$. 
Then, we estimate $\lambda_{\min}((Q^{\top}U_{+})^{\top}(Q^{\top}U_{+}))$ as follows:
\begin{equation*}
\arraycolsep=0.1em
\begin{array}{lcl}
\lambda_{\min}((Q^{\top}U_{+})^{\top}(Q^{\top}U_{+})) & = &  \lambda_{\min}\left(\Id + 0.5\alpha( H^{\top} + H) + \alpha^2H^{\top}H \right) \vspace{1ex} \\
&\overset{\tiny\text{\cite[9.13.6.]{Bernstein2005}}}{\geq} & 1 - 0.5\alpha \lambda_{\max}(H^{\top} + H) \overset{\tiny\text{\cite[5.11.25]{Bernstein2005}}}{\geq} 1 - \alpha \sigma_{\max}(H) \vspace{1ex}\\
& = & 1 - \alpha \sigma_{\max}\left( (U^{\top}U)^{-1}U^{\top}ZV(V^{\top}V)^{-1}\right) \vspace{1ex}\\
&\geq & 1 -  \frac{\alpha \sigma_{\max}(U^{\top}Z)}{\sigma_{\min}(U)^2\sigma_{\min}(V)} \geq 1 - \frac{\alpha \Vert U^{\top}Z \Vert_F}{\sigma_{\min}^3},
\end{array}
\end{equation*}
where $\sigma_{\min} = \min\set{\sigma_{\min}(U),\sigma_{\min}(V)}$ and $\Vert U^{\top}Z \Vert_F \geq \sigma_{\max}(U^{\top}Z)$.
We note that $\Vert\Phi(U,V)\Vert_F \geq \Vert \nabla_U{\Phi}(U,V)\Vert_F = L_{\Phi}\Vert U^{\top}Z\Vert$. 
Substituting this estimate into the last inequality and noting from \eqref{eq:step_size_min} that $\alpha \leq \frac{\sigma_{\min}^3L_{\Phi}}{2\Vert\nabla{\Phi}(U,V)\Vert_F}$, we obtain
\begin{align*}
\lambda_{\min}((Q^{\top}U_{+})^{\top}(Q^{\top}U_{+})) \geq 1  - \alpha\frac{\Vert \nabla{\Phi}(U, V)\Vert_F}{L_{\Phi}\sigma_{\min}^3} \geq 1 - \frac{L_{\Phi}}{2L_{\Phi}} = \frac{1}{2} > 0.
\end{align*}
Combining this estimate and \eqref{eq:lm24_proof_est1} we have $\lambda_{\min}(U_{+}^{\top}U_{+}) \geq 0.5\lambda_{\min}(U^{\top}U)$. Hence, we conclude that  $\rank{U_{+}} = r$. 
With a similar proof, we can show that $\rank{V_{+}} = r$.
\Eproof

\subsection{The proof of Theorem \ref{th:convergence_guarantee1}: Global convergence of  GN method}\label{apdx:th:convergence_guarantee1}
By Lemma  \ref{le:descent_dir},  we can see that the backtracking linesearch step at Step \ref{step:linesearch} of Algorithm \ref{alg:A1} is finite and $\alpha_k > 0$. 
The inequality \eqref{eq:ls_condition} guarantees that $\Phi(U_{k+1}, V_{k+1}) < \Phi(U_k, V_k)$. 
Hence, the sequence $\set{\Phi(U_k, V_k)}$ is decreasing and bounded from below by $\Phi^{\star}$. It converges to a limit point $\Phi^{*}$.
Now, using \eqref{eq:ls_condition} we obtain
\begin{equation*}
\sum_{k=0}^{n}\alpha_k\Vert\nabla{\Phi}(U_k,V_k)\Vert_F^2 \leq \Phi(U_0,V_0) - \Phi(U_{n+1},V_{n+1}) \leq \Phi(U_0,V_0) - \Phi^{\star} < +\infty.
\end{equation*}
Taking the limit in this inequality as $n\to\infty$, we obtain $\sum_{k=0}^{\infty}\alpha_k\Vert\nabla{\Phi}(U_k,V_k)\Vert_F^2~<+\infty$.
Consequently, $\lim_{k\to\infty}\alpha_k\Vert\nabla{\Phi}(U_k,V_k)\Vert_F^2 = 0$. 
This proves the first part \eqref{eq:sum_inf}. 

In order to prove the second part, we need to show that $\alpha_k \geq \alpha > 0$ for all~$k$ sufficiently large.
Indeed, by our assumption that $\set{[U_k,V_k]}$ is bounded.
Hence, $\Vert\Phi^{\prime}(U_k, V_k)\Vert_F \leq K_1 <+\infty$. 
Similarly, $\Vert\nabla{\Phi}(U_k, V_k) \Vert_F \leq K_2 <+\infty$ and $\max\set{\sigma_{\max}(U_k),\sigma_{\max}(V_k)} \leq K_3 <+\infty$. 
Using these arguments and condition \eqref{eq:bounded_assumption} into \eqref{eq:step_size_min}, we obtain 
\begin{equation*}
2\alpha_k \geq \underline{\alpha} \overset{\tiny\eqref{eq:step_size_min}}{\geq} 2\alpha := \min\set{1, \frac{L_{\Phi}\underline{\sigma}^3}{2K_2},\frac{3\underline{\sigma}^4}{32K_1K_3^2}} > 0.
\end{equation*}
Using this lower bound into \eqref{eq:sum_inf} we have $\lim\limits_{k\to\infty}\Vert\nabla{\Phi}(U_k,V_k)\Vert_F^2 \leq \alpha^{-1}\lim\limits_{k\to\infty}\alpha_k\Vert\nabla{\Phi}(U_k, V_k)\Vert_F^2 = 0$,
which implies \eqref{eq:lim_inf}.

By our assumption, $\set{X_k}$ generated by Algorithm \ref{alg:A1} is bounded. 
Hence, there exists a limit point $X_{\star} := [U_{\star}, V_{\star}]$.
Passing through the limit \eqref{eq:lim_inf} via subsequence, we can see that $\nabla{\Phi}(U_{\star},V_{\star}) = 0$,  and hence, $X_{\star}$ satisfies the optimality condition \eqref{eq:opt_cond}.
\Eproof

\subsection{The proof of Lemma \ref{le:descent2}: Descent property of $\Lc_{\rho}$}\label{apdx:le:descent2}
We first prove part (a).
Since $\set{[U_k,V_k]}$ is bounded by our assumption, and since $\lim\limits_{k\to\infty}\Vert W_k - \Ac(U_kV_k^{\top}) + B\Vert_F = 0$ due to part (b), the sequence $\set{W_k}$ is also bounded.

Now, we prove part (b) for \textbf{Option 1}.
First, since $[U_{k+1}, V_{k+1}]$ is updated by Step 5 of Algorithm \ref{alg:A2} that satisfies the backtracking linesearch condition  \eqref{eq:ls_cond2}, we have
\begin{equation}\label{eq:lm41_est1}
\mathcal{Q}_k( U_{k+1}, V_{k+1} ) \leq \mathcal{Q}_k(U_k, V_k) - 0.5c_1\alpha_k\Delta_k^2,
\end{equation}
where $\mathcal{Q}_k$ is defined by \eqref{eq:ls_cond2}, $E_k := \Ac(U_kV_k^{\top})-B +\rho^{-1}\Lambda_k$, and  $\Delta_k^2$ is 
\begin{equation}\label{eq:xi_k_term}
\Delta_k^2 := \Vert U_k^{\top}\Ac^{\ast}(E_k-W_k)\Vert_F^2 + \Vert \Ac^{\ast}(E_k-W_k)V_k \Vert_F^2.
\end{equation}
This condition implies
\begin{align}\label{eq:term1}
\Lc_{\rho}(U_{k+1},V_{k+1},W_k,\Lambda_k) \leq \Lc_{\rho}(U_{k},V_{k},W_k,\Lambda_k) - \frac{c_1\rho\alpha_k}{2}\Delta_k^2.
\end{align}
Second, we consider the objective function $h(W) := \phi(W) + (\rho/2)\norm{W - C_k}_F^2$ of \eqref{eq:aug_method2_b}, where $C_k := \Ac(U_{k+1}V_{k+1}) - B + \rho^{-1}\Lambda_k$. 
Since $h(\cdot)$ is strongly convex with the strong convexity parameter $\rho + \mu_{\phi}$, and $W_{k+1}$ is the optimal solution of $h$, we have
\begin{equation*}
h(W_{k+1}) \leq h(W_k) - ((\rho+\mu_{\phi})/2)\Vert W_{k+1} - W_k\Vert_F^2.
\end{equation*}
Using this inequality, and the definition of $h$ and $\Lc_{\rho}$, we can show that
\begin{equation}\label{eq:term2}
\Lc_{\rho}(U_{k+1},V_{k+1},W_{k+1},\Lambda_k) \leq \Lc_{\rho}(U_{k+1},V_{k+1},W_k,\Lambda_k) - \frac{(\rho+\mu_{\phi})}{2}\Vert W_{k+1} - W_k\Vert_F^2.
\end{equation}
In addition, since $\phi$ is $L_{\phi}$-smooth, we can write down the optimality condition of \eqref{eq:aug_method2_b} as $\nabla{\phi}(W_{k+1}) + \rho(W_{k+1} - C_k) = 0$. Using the definition of $C_k$ and \eqref{eq:aug_method2_c} we get $\Lambda_{k+1} = \nabla{\phi}(W_{k+1})$. Hence, we can derive
\begin{equation}\label{eq:term2b}
\Vert \Lambda_{k+1} - \Lambda_k\Vert_F = \Vert \nabla{\phi}(W_{k+1}) - \nabla{\phi}(W_k)\Vert_F \leq L_{\phi}\Vert W_{k+1} - W_k\Vert_F,
 \end{equation}
 which is the first inequality in \eqref{eq:bound_lambda}.
 The boundedness of $\set{\Lambda_k}$ also follows from the relation $\Lambda_{k+1} = \nabla{\phi}(W_{k+1})$ and the boundedness of $\set{W_k}$.
 
Third, since $\Lambda_k$ is updated by \eqref{eq:aug_method2_c}, using the definition of $\Lc_{\rho}$, it is easy to show that
\begin{equation}\label{eq:term3}
\arraycolsep=-0.1em
\begin{array}{lcl}
\Lc_{\rho}(U_{k+1},V_{k+1},W_{k+1},\Lambda_{k+1}) &= & \Lc_{\rho}(U_{k+1},V_{k+1},W_{k+1},\Lambda_k) + \rho^{-1}\Vert \Lambda_{k+1} - \Lambda_k\Vert_F^2 \vspace{1ex}\\
&\overset{\tiny\eqref{eq:term2b}}{\leq} & \Lc_{\rho}(U_{k+1},V_{k+1},W_{k+1},\Lambda_k) + \rho^{-1}L_{\phi}^2\Vert W_{k+1} - W_k\Vert_F^2.
\end{array}
\end{equation}
Summing up \eqref{eq:term1}, \eqref{eq:term2} and \eqref{eq:term3} we get \eqref{eq:descent2}.

Finally, we prove (b) for \textbf{Option 2}.
We consider the gradient step \eqref{eq:second_subprob_sol} instead of \eqref{eq:aug_method2_b}.
Using the optimality condition of \eqref{eq:aug_method2b} and \eqref{eq:aug_method2_c}, we can derive $\Lambda_{k+1} = \nabla{\phi}(W_k) + L_{\phi}(W_{k+1} - W_k)$. 
Using this relation and the Lipschitz continuity of $\nabla{\phi}$, we have
\begin{align}\label{eq:term6a}
\Vert\Lambda_{k+1} - \Lambda_k\Vert_F &= \Vert L_{\phi}(W_{k+1}-W_{k-1}) + \nabla{\phi}(W_k) - \nabla{\phi}(W_{k-1})\Vert_F\nonumber\\
&\leq L_{\phi}\big[ \Vert W_{k+1}-W_{k-1}\Vert_F + \Vert W_k-W_{k-1}\Vert_F\big],
\end{align}
which is exactly the second expression of \eqref{eq:bound_lambda}.
Using $\Lambda_{k+1} = \nabla{\phi}(W_k) + L_{\phi}(W_{k+1} - W_k)$, similar above, we can also show the boundedness of $\set{\Lambda_k}$.

Now, since we apply the gradient step to solve \eqref{eq:aug_method2_b}, with $h$ defined as in \eqref{eq:term2}, it is well-known that
\begin{equation*}
h(W_{k+1}) \leq h(W_k) - ((L_{\phi}+\rho)/2)\Vert W_{k+1}-W_k\Vert_F^2,
\end{equation*}
which implies 
\begin{equation}\label{eq:term6b}
\Lc_{\rho}(U_{k+1},V_{k+1},W_{k+1},\Lambda_k) \leq \Lc_{\rho}(U_{k+1},V_{k+1},W_k,\Lambda_k) - \frac{(\rho+ L_{\phi})}{2}\Vert W_{k+1} - W_k\Vert_F^2.
\end{equation}
Summing up \eqref{eq:term1}, \eqref{eq:term6b} and the first equality of \eqref{eq:term3} we obtain
\begin{align}\label{eq:term6c}
\Lc_{\rho}(U_{k+1},V_{k+1},W_{k+1},\Lambda_{k+1}) &= \Lc_{\rho}(U_k,V_k,W_k,\Lambda_k) - (1/2)c_1\rho\alpha_k\Delta^2_k - T_k,
\end{align}
where $T_k := \frac{(\rho+ L_{\phi})}{2}\Vert W_{k+1} - W_k\Vert_F^2 -  \rho^{-1}\Vert \Lambda_{k+1} - \Lambda_k\Vert_F^2$.
Finally, using \eqref{eq:term6a}, we can estimate $\Vert\Lambda_{k+1}-\Lambda_k\Vert_F$ as follows:
\begin{equation*}
\begin{array}{lcl}
\Vert\Lambda_{k+1}-\Lambda_k\Vert_F^2 &\leq & L_{\phi}^2\big[\Vert W_{k+1}-W_{k-1}\Vert_F + \Vert W_k-W_{k-1}\Vert_F\big]^2 \vspace{1ex}\\
&\leq & 2L_{\phi}^2\Vert W_{k+1}-W_k\Vert_F^2 + 4L_{\phi}^2\Vert W_k-W_{k-1}\Vert_F^2.
\end{array}
\end{equation*}
Hence, $T_k \geq (0.5(\rho+L_{\phi}) - 2\rho^{-1}L_{\phi}^2)\Vert W_{k+1}-W_k\Vert_F^2 - 4\rho^{-1}L_{\phi}^2\Vert W_k-W_{k-1}\Vert_F^2$.
Substituting this estimate of $T_k$ into \eqref{eq:term6c} we obtain \eqref{eq:descent2}.
\Eproof

\subsection{The proof of Theorem \ref{th:global_convergence2}: Global convergence of ADMM-GN}\label{apdx:th:global_convergence2}
We first prove for \textbf{Option 1}.
Let us define $\eta : = \rho^{-1}(\rho^2 + \mu_{\phi}\rho - 2L_{\phi}^2)$. 
Then, $\eta > 0$ if we choose $\rho > 0.5((\mu_{\phi} + 8L_{\phi}^2)^{1/2} + \mu_{\phi})$ as given by \eqref{eq:descent2b} in Lemma \ref{le:descent2}.
Hence,  the sequence $\set{\Lc_{\rho}(U_k,V_k,W_k, \Lambda_k)}$ is strictly decreasing, it is bounded from bellow due to Assumption A.\ref{as:A1} and the boundedness of $\set{(U_k, V_k, W_k, \Lambda_k]}$. 
It converges to a finite value $\Lc_{\rho}^{\star}$. In addition, \eqref{eq:descent2} implies
\begin{align}\label{lm41_est6}
\begin{array}{ll}
&\displaystyle\lim_{k\to\infty}\norm{W_{k+1} - W_k}_F = 0, \vspace{0.0ex}\\
&\displaystyle\lim_{k\to\infty} \alpha_k\big\Vert U_k^{\top}\Ac^{\ast}\big(\rho^{-1}\Lambda_k + \Ac(U_kV_k^{\top}) - B - W_k\big)\big\Vert_F^2 = 0,  \quad \text{and} \vspace{0.0ex}\\
&\displaystyle\lim_{k\to\infty} \alpha_k\big\Vert \Ac^{\ast}\big(\rho^{-1}\Lambda_k + \Ac(U_kV_k^{\top}) - B - W_k\big)V_k \big\Vert_F^2 = 0.
\end{array}
\end{align}
Under condition \eqref{eq:bounded_assumption}, similar to the proof of Theorem \ref{th:convergence_guarantee1} we can show that $\alpha_k \geq 0.5\underline{\alpha} > 0$ for $k$ sufficiently large.
Hence, the two last limits of \eqref{lm41_est6} imply 
\begin{equation}\label{eq:lm41_est6b}
\begin{array}{ll}
&\displaystyle\lim_{k\to\infty}  \big\Vert U_k^{\top}\Ac^{\ast}\big(\Lambda_k + \rho\big(\Ac(U_kV_k^{\top}) - B - W_k\big)\big)\big\Vert_F = 0 \quad \text{and} \vspace{0.5ex}\\
&\displaystyle\lim_{k\to\infty} \big\Vert \Ac^{\ast}\big(\Lambda_k + \rho\big(\Ac(U_kV_k^{\top}) - B - W_k\big)V_k \big\Vert_F = 0.
\end{array}
\end{equation}
On the other hand, using Lemma \ref{le:descent2}(a), \eqref{eq:aug_method2_c}, and the first limit in \eqref{lm41_est6}, we obtain
\begin{equation}\label{eq:lm41_est6c}
\arraycolsep=-0.1em
\begin{array}{lcl}
\lim_{k\to\infty}\Vert\Ac(U_{k+1}V_{k+1}^{\top}) - W_{k+1} - B\Vert_F &\overset{\tiny\eqref{eq:aug_method2_c}}{=} & \rho^{-1}\lim_{k\to\infty} \Vert \Lambda_{k+1} -\Lambda_k\Vert_F \vspace{0.5ex}\\
& \overset{\tiny\text{Lemma \ref{apdx:le:descent2}(a)}}{\leq} & \rho^{-1}L_{\phi}\lim_{k\to\infty}\norm{W_{k+1} - W_k}_F \overset{\tiny\eqref{lm41_est6}}{=} 0.
\end{array}
\end{equation}
We consider a convergent subsequence $\set{[U_{k_i},V_{k_i}]}_{i\in\mathbb{Nc}}$ with the limit $[U_{\ast}, V_{\ast}]$.
Then, the limit \eqref{eq:lm41_est6c} shows that the corresponding subsequence $\set{W_{k_i}}$ also converges to $W_{\ast}$ such that $W_{\ast} = \Ac(U_{\ast}V_{\ast}^{\top}) - B$, which is the last condition in \eqref{eq:kkt_cond}. 

Now, using the limit in \eqref{eq:lm41_est6c} and combining with the triangle inequality, we get
\begin{equation*}
\begin{array}{lcl}
\Vert U_k^{\top}\Ac^{\ast}(\Lambda_k)\Vert_F &\leq & \big\Vert U_k^{\top}\Ac^{\ast}\big(\Lambda_k + \rho\big(\Ac(U_kV_k^{\top}) - B - W_k\big)\big)\big\Vert_F \vspace{1ex}\\
&& + {~} \rho\big\Vert U_k^{\top}\Ac^{\ast}\big(\Ac(U_kV_k^{\top}) - B - W_k\big)\big\Vert_F \overset{\tiny\eqref{eq:lm41_est6b},\eqref{eq:lm41_est6c}}{\to} 0 \quad \text{as} \quad k_i\to\infty.
\end{array}
\end{equation*}
This implies $U_{\ast}^{\top}\Ac^{\ast}(\Lambda_{\ast}) = 0$ via subsequence. 
Similarly, we can also show that $\Ac^{\ast}(\Lambda_{\ast})V_{\ast} = 0$.
These are the second and the third conditions in \eqref{eq:kkt_cond}.
Finally, the first condition of \eqref{eq:kkt_cond} follows directly from the relation $\Lambda_k = \nabla{\phi}(W_k)$ as the optimality condition of \eqref{eq:aug_method2_b} by taking the limit via subsequence.

We have shown in the above steps that the limit point $(U_{\ast}, V_{\ast}, W_{\ast}, \Lambda_{\ast})$ satisfies the optimality condition \eqref{eq:kkt_cond} of \eqref{eq:constr_LMA_prob2}. By eliminating $\Lambda_{\ast}$ and $W_{\ast}$ in  \eqref{eq:kkt_cond} we obtain  \eqref{eq:opt_cond}, which shows that any limit point $[U_{\ast}, V_{\ast}]$ of $\set{[U_k, V_k]}$ is a stationary point of \eqref{eq:LMA_prob}.
The proof of \eqref{eq:norm_grad_limit} can be done as in Theorem \ref{th:convergence_guarantee1}.

We prove for \textbf{Option 2}.
We note that if $\rho > 3L_{\phi}$, then we can examine from \eqref{eq:descent2b} that $\eta_1 > \eta_0$. If we denote by $\Lc_k := \Lc_{\rho}(U_k,V_k,W_k,\Lambda_k)$ and $r_k := \Vert W_k - W_{k-1}\Vert_F$ for $k\geq 1$, then we can write \eqref{eq:descent2} as
\begin{equation*}
\Lc_{k+1} + \frac{\eta_0}{2}r_{k+1}^2 \leq \Lc_k + \frac{\eta_0}{2}r_k^2 - \frac{c_1\rho}{2}\Delta_k^2 - \frac{(\eta_1-\eta_0)}{2}r_{k+1}^2.
\end{equation*}
By induction, we can show from this inequality that $\sum_{k=0}^{\infty}\left( c_1\rho\Delta_k^2 + (\eta_1-\eta_0)r_{k+1}^2\right) = 0$, which implies \eqref{lm41_est6}.
With the same proof as in \textbf{Option 1} we obtain the same conclusions of the theorem as in \textbf{Option 1}.
\Eproof

\subsection{The proof of Theorem \ref{th:local_convergence}: Local convergence of GN method}\label{apdx:th:local_convergence}
Let us define $x := [\vec{U}, \vec{V^{\top}}] \in\R^{(m+n)r}$ the vecterization of $U$ and $V$, and $R(x) :=  \Ac(UV^{\top}) - B$  the residual term.
We can compute the Jacobian $J_R(x)$ of $R$ at $x$ as $J_R(x) = A[V \otimes \Id_m, \Id_n\otimes U^{\top}]\in\R^{l\times (m+n)r}$, where $A$ is the matrix form of the linear operator $\Ac$.
The objective function $\Phi(U,V)$ can be written as $\Phi(x) = \phi(R(x))$. Its gradient and Hessian are given by
\begin{equation}\label{eq:jac_hess}
\left\{\begin{array}{lcl}
\nabla{\Phi}(x)    &  = & J_R(x)^{\top}\nabla{\phi}(R(x)) \quad \text{and}\vspace{0.0ex}\\
\nabla^2{\Phi}(x) & = & J_R(x)^{\top}\nabla^2{\phi}(R(x))J_R(x) + \displaystyle\sum_{i=1}^l\frac{\partial{\phi(R(x))}}{\partial{R_i}}\nabla^2{R_i(x)}. 
\end{array}\right.
\end{equation}
First, we show that under Assumption A.\ref{as:A2}(a), $\nabla^2{\Phi}$ is also Lipschitz continuous in $\Nc(x_{\star})$ of $x_{\star} \in\Xc_{\star}$. 
Indeed, $\nabla^2R(x)$ is bounded in $\Nc(x_{\star})$ by $M_{R_i''}$, and $\nabla^2R_i(\cdot)$ is Lipschitz continuous with the Lipschitz constant $L_{R^{''}_i}$. 
In addition, $J_R(\cdot)$ is also bounded in $\Nc(x_{\star})$ by $M_{R'}$, and $R(\cdot)$ is also Lipschitz continuous with the Lipschitz constant $L_R$.
Since $\nabla^2\phi$ is Lipschitz continuous in $\Nc(R(x_{\star}))$, $\frac{\partial{\phi(R(x)}}{\partial{R_i}}$ is also bounded by $M^i_{\phi'}$, and Lipschitz continuous  in $\Nc(R(x_{\star}))$ with the Lipschitz constant $L^i_{\phi'}$.
Combining these statements and \eqref{eq:jac_hess}, we can show that for any $x, \hat{x}\in\Nc(x_{\star})$, the following estimate holds:
\begin{equation*}
\begin{array}{lcl}
\Vert\nabla^2{\Phi}(x) - \nabla^2{\Phi}(\hat{x})\Vert & \leq & \left\Vert J_R(x)^{\top}\nabla^2{\phi}(R(x))J_R(x) - J_R(\hat{x})^{\top}\nabla^2{\phi}(R(\hat{x}))J_R(\hat{x}) \right\Vert \vspace{1ex}\\
&& + {~} \left\Vert \sum_{i=1}^l\left[\frac{\partial{\phi(R(x))}}{\partial{R_i}}\nabla^2{R_i(x)} - \frac{\partial{\phi(R(\hat{x}))}}{\partial{R_i}}\nabla^2{R_i(\hat{x})}\right]\right\Vert \vspace{1ex}\\
& \leq & \left\Vert J_R(x)^{\top}\nabla^2{\phi}(R(x))\left(J_R(x) - J_R(\hat{x})\right) \right\Vert \vspace{1ex}\\
&& + {~} \left\Vert J_R(x)^{\top}\left( \nabla^2{\phi}(R(x)) - \nabla^2{\phi}(R(\hat{x}))\right)J_R(\hat{x}) \right\Vert  \vspace{1ex}\\
&& + {~} \left\Vert \left(J_R(x) - J_R(\hat{x})\right)^{\top}\nabla^2{\phi}(R(\hat{x})) J_R(\hat{x}) \right\Vert \vspace{1ex}\\
&& + {~} \sum_{i=1}^l\Big[ \Big\Vert \frac{\partial{\phi(R(x))}}{\partial{R_i}}\big(\nabla^2{R_i}(x) - \nabla^2{R_i}(\hat{x})\big)  \Big\Vert  \vspace{1ex}\\
&& + {~} \Big\Vert \left(\frac{\partial{\phi(R(x))}}{\partial{R_i}} - \frac{\partial{\phi(R(\hat{x}))}}{\partial{R_i}}\right)\nabla^2{R_i}(\hat{x})  \Big\Vert \Big] \vspace{1ex}\\
&\leq & \Big( 2M_{R'}M_{\phi''}L_{R'} + M_{R'}^2L_{\phi''} + \sum_{i=1}^l(M_{R_i^{''}}L^i_{\phi'} + L_{R_i^{''}}M^i_{\phi'})\Big)\Vert x -\hat{x}\Vert.
\end{array}
\end{equation*}
This inequality shows that $\nabla{\Phi}$ is Lipschitz continuous in $\Nc(x_{\star})$ with the Lipschitz constant $L_{\Phi^{''}} :=  2M_{R'}M_{\phi''}L_{R'} + M_{R'}^2L_{\phi''} + \sum_{i=1}^l(M_{R_i^{''}}L^i_{\phi'} + L_{R_i^{''}}M^i_{\phi'}) > 0$.

Next, we consider the GN direction $D_{X_k}$ in \eqref{eq:subprob2}. 
Let $d := [\vec{D_U}, \vec{D_V^{\top}}]$ and $H_0(x) := \begin{bmatrix}V^{\top}\otimes U & V^{\top}V \otimes \Id_m \\ \Id_n\otimes U^{\top}U & V\otimes U^{\top}\end{bmatrix}$.
Due to the full-rankness of $U$ and $V$, by using the result in \cite{campbell2009generalized} we can show that $H_0(x)^{\dagger}$ is bounded by $M_h$, i.e.:
\begin{equation}\label{eq:bound_H0dagger}
\Vert H_0(x)^{\dagger} \Vert \leq M_h < +\infty, \quad \forall x\in \Nc(x_{\star}),
\end{equation}
Moreover, we can see from \eqref{eq:opt_cond_subprob2} that $[\vec{ZV}, \vec{U^{\top}Z}] = -L_{\Phi}^{-1} J_R(x)^{\top}\nabla{\phi}(R(x))$. 
Hence, \eqref{eq:opt_cond_subprob2}  can be written as $H_0(x)d = -L_{\Phi}^{-1} J_R(x)^{\top}\nabla{\phi}(R(x))$, which implies $d = -L_{\Phi}^{-1}H_0(x)^{\dagger}\nabla{\Phi}(x)$. 
The full-step GN scheme becomes
\begin{equation}\label{eq:fs_gn_scheme0}
x_{+} = x  + d = x -  L_{\Phi}^{-1}H_0(x)^{\dagger}\nabla{\Phi}(x).
\end{equation}
We consider the residual term $r = x - x_{\star}$, where $x_{\star} := [\vec{U_{\star}}, \vec{V_{\star}^{\top}}]\in\Xc_{\star}$ is a given stationary point of \eqref{eq:LMA_prob}.
From \eqref{eq:fs_gn_scheme0} we can write
\begin{equation*}
\begin{array}{lcl}
r_{+} &= & x_{+} - x_{\star} = r - L_{\Phi}^{-1}H_0(x)^{\dagger}\nabla{\Phi}(x) \vspace{1ex}\\
& = & r - L_{\Phi}^{-1}H_0(x)^{\dagger}\left[\nabla{\Phi}(x) - \nabla{\Phi}(x_{\star})\right] \vspace{1ex}\\
&= & \left[\Id - L_{\Phi}^{-1}H_0(x)^{\dagger}\nabla^2\Phi(x_{\star})\right]r  \vspace{1ex}\\
&&  - {~}  L_{\Phi}^{-1}H_0(x)^{\dagger}\left[\int_0^1\left(\nabla^2{\Phi}(x_{\star} + \tau(x - x_{\star})) - \nabla^2\Phi(x_{\star})\right)(x - x_{\star})d\tau\right].
\end{array}
\end{equation*}
Using condition \eqref{eq:small_residual} and the Lipschitz continuity of $\nabla{\Phi}$, this expression leads to
\begin{equation}\label{eq:th32_proof2}
\begin{array}{lcl}
\Vert r_{+}\Vert &\leq &  \Vert \left(\Id - L_{\Phi}^{-1}H_0(x)^{\dagger}\nabla^2\Phi(x_{\star})\right)r\Vert \vspace{1ex}\\
&& + {~} L_{\Phi}^{-1}\Vert H_0(x)^{\dagger}\Vert\int_0^1\Vert \nabla^2{\Phi}(x_{\star} +\tau(x-x_{\star})) - \nabla^2\Phi(x_{\star})\Vert\Vert x - x_{\star}\Vert d\tau \vspace{1ex}\\
& \leq & \kappa(x_{\star})\Vert r\Vert + \frac{1}{2}L^{-1}L_{\Phi^{''}}\Vert H_0(x)^{\dagger}\Vert\Vert r \Vert^2 \vspace{1ex}\\
&\leq & \left( \bar{\kappa} + 0.5L_{\Phi}^{-1}L_{\Phi^{''}}K_h \Vert r\Vert \right)\Vert r\Vert.
\end{array}
\end{equation}
Since $r = x - x_{\star} = \vec{X - X_{\star}}$, we can write \eqref{eq:th32_proof2} as
\begin{equation*}
\Vert X_{+} - X_{\star}\Vert_F \leq \left( \bar{\kappa} + 0.5L_{\Phi}^{-1}L_{\Phi^{''}}K_h \Vert X-X_{\star}\Vert_F \right)\Vert X-X_{\star}\Vert_F,
\end{equation*}
which is exactly \eqref{eq:gn_local_est1} with $K_1 := L_{\Phi}^{-1}L_{\Phi^{''}}K_h > 0$.

Next, we  prove quadratic convergence of the full-step GN scheme.
Under Assumption \ref{as:A1}, it follows from \cite{golub1973differentiation} that  there exists a neighborhood $\Nc(x_{\star})$ of $x_{\star}$ such that $H_0(\cdot)^{\dagger}$ is Lipschitz continuous in $\Nc(x_{\star}$ with the Lipschitz constant $L_{H} > 0$.
Here, we use the same $\Nc(x_{\star}$ as in Assumption \ref{as:A2}. Otherwise, we can shrink it if necessary.
We consider the condition $H(X_{\star})^{\dagger}\nabla^2{\Phi}(X_{\star})  = L_{\Phi}\Id$.
Reforming this condition into vector form, we have $H_0(x_{\star})^{\dagger}\nabla^2{\Phi}(x_{\star})  = L_{\Phi}\Id$, which is equivalent to $\Id - L_{\Phi}^{-1}H(X_{\star})^{\dagger}\nabla^2{\Phi}(X_{\star}) = 0$.
Using the last condition, and the Lipschitz continuity of $H_0^{\dagger}(\cdot)$, we can show that
\begin{equation*}
\begin{array}{lcl}
S(x_{\star}) &:= & \Vert \left[\Id - L_{\Phi}^{-1}H_0(x)^{\dagger}\nabla^2{\Phi}(x_{\star})\right](x-x_{\star})\Vert \vspace{1ex}\\
&\leq & \Vert \left[\Id - L_{\Phi}^{-1}H_0(x_{\star})^{\dagger}\nabla^2{\Phi}(x_{\star})\right](x - x_{\star})\Vert + L_{\Phi}^{-1}\Vert \left(H_0(x)^{\dagger} - H_0(x_{\star})^{\dagger}\right)(x-x_{\star})\Vert  \vspace{1ex}\\
&\leq &  L_{\Phi}^{-1}\Vert H_0(x)^{\dagger} - H_0(x_{\star})^{\dagger}\Vert\Vert x-x_{\star}\Vert \vspace{1ex}\\
&\leq &  L_{\Phi}^{-1}L_{H}\Vert x - x_{\star}\Vert^2,~~~\forall x\in\Nc(x_{\star}).
\end{array}
\end{equation*}
Substituting this $S(x_{\star})$ estimate into \eqref{eq:th32_proof2} we get $\Vert r_{+}\Vert \leq L_{\Phi}^{-1}(L_{H} + 0.5L_{\Phi^{''}}K_h)\Vert r\Vert^2$, which is reformed into the matrix form as
\begin{equation*}
\Vert X_{+} - X_{\star}\Vert_F \leq 0.5K_2\Vert X- X_{\star}\Vert_F^2,~\forall X\in\Nc(X_{\star}),~\text{where}~K_2 := L_{\Phi}^{-1}\left(2L_{H} + L_{\Phi^{''}}K_h\right).
\end{equation*}
In order to guarantee the monotonicity of $\set{\Vert X-X_{\star}\Vert_F}$, we require $\Vert X_{+} - X_{\star}\Vert_F \leq 0.5K_1\Vert X - X_{\star}\Vert_F^2 <  \Vert X - X_{\star}\Vert_F$, which implies $ \Vert X - X_{\star}\Vert_F < 2K_2^{-1}$. 
Hence, if we choose $X_0\in\Nc(X_{\star})$ such that $\Vert X_0-X_{\star}\Vert_F < 2K_2^{-1}$, then $\Vert X_k-X_{\star}\Vert_F < 2K_2^{-1}$ for all $k \geq 0$ and $\set{\Vert X_k-X_{\star}\Vert_F}$ is monotone. 
Moreover, $\Vert X_{k+1}-X_{\star}\Vert_F \leq 0.5K_2\Vert X_k - X_{\star}\Vert_F^2$ shows that this sequence converges quadratically to zero. 
Hence, $\set{X_k}$ converges to $X_{\star}$ at a quadratic rate.
Here, we can easily check that $K_2 > K_1$.

Finally, if $\bar{\kappa} \in (0, 1)$, then for all $\geq 0$, the estimate \eqref{eq:gn_local_est1}  implies that 
\begin{equation*}
\Vert X_{k+1} - X_{\star}\Vert_F \leq \left(\bar{\kappa} + 0.5K_1\Vert X_k-X_{\star}\Vert_F\right)\Vert X_k-X_{\star}\Vert_F.
\end{equation*}
In order to guarantee $\Vert X_{k+1} -X_{\star}\Vert_F < \Vert X_k-X_{\star}\Vert_F$, we require $\bar{\kappa} + 0.5K_1\Vert X_k-X_{\star}\Vert_F < 1$, which leads to $\Vert X_k-X_{\star}\Vert_F < 2K_1^{-1}(1-\bar{\kappa})$. 
Hence, if we take $\bar{r}_0 < 2K_1^{-1}(1-\bar{\kappa})$, and choose $X_0\in\Nc(X_{\star})$ such that $\Vert X_0 - X_{\star}\Vert_F \leq \bar{r}_0$, then $\Vert X_k - X_{\star}\Vert_F \leq \bar{r}_0$ for all $k\geq 0$. 
In addition, we have $\Vert X_{k+1} - X_{\star}\Vert_F \leq \left(\bar{\kappa} + 0.5K_1\Vert X_k-X_{\star}\Vert_F\right)\Vert X_k-X_{\star}\Vert_F \leq (\bar{\kappa} + 0.5K_1\bar{r}_0)\Vert X_k - X_{\star}\Vert_F$, which shows that $\set{\Vert X_k-X_{\star}\Vert_F}$ converges to zero at a linear rate with the contraction factor $\omega := \bar{\kappa} + 0.5K_1\bar{r}_0 < 1$.
\Eproof

\bibliographystyle{plain}

\end{document}